\documentclass[11pt]{amsart}
\usepackage{amssymb}
\usepackage[all]{xy}
\usepackage{geometry}

\title{Character Sheaves and Depth-zero Representations}

\author{Anne-Marie~Aubert}
\address{Institut de Math\'ematiques de Jussieu, U.M.R. 7586 du C.N.R.S., Projet Formes Automorphes, Universit\'e Pierre et Marie Curie, F-75252 Paris Cedex 05, FRANCE} 
\email{aubert@math.jussieu.fr}

\author{Clifton~Cunningham}
\address{
Department of Mathematics, University of Calgary, 2500 University Drive N.W., Calgary Alberta T2N 1N4, CANADA}
\email{
cunning@math.ucalgary.ca}



\newtheorem{theorem}{Theorem}[section]
\newtheorem{proposition}[theorem]{Proposition}
\newtheorem{corollary}[theorem]{Corollary}
\newtheorem{lemma}[theorem]{Lemma}
\theoremstyle{definition}
\newtheorem{definition}[theorem]{Definition}
\newtheorem{example}[theorem]{Example}
\theoremstyle{remark}
\newtheorem{remark}[theorem]{Remark}
\numberwithin{equation}{section}

\newcommand{\ZZ}{{\mathbb{Z}}}
\newcommand{\QQ}{{\mathbb{Q}}}

\newcommand{\CC}{{\mathbb{C}}}
\newcommand{\EE}{{\bar{\QQ}_\ell}}
\newcommand{\FF}{{\mathbb{F}}}

\newcommand{\KK}{{\mathbb{K}}}
\newcommand{\RK}{{\mathfrak{o}_{\KK}}}
\newcommand{\PK}{{\mathfrak{p}_{\KK}}}
\newcommand{\kK}{{\Bbbk}}

\newcommand{\Knr}{{\mathbb{K}_1^{nr}}}
\newcommand{\Rnr}{{\mathfrak{o}_{\Knr}}}

\newcommand{\knr}{{\bar\FF_q}}

\newcommand{\Kq}{{\mathbb{K}_1}}
\newcommand{\Rq}{{\mathfrak{o}_{\Kq}}}

\newcommand{\kq}{{\FF_q}}

\newcommand{\Gal}[1]{{\operatorname{Gal}(#1)}}

\newcommand{\GKq}{{\Ksch{G}(\Kq)}}
\def\cf{{{\it cf.\ }}}
\def\ie{{{\it i.e.,\ }}}
\newcommand{\ceq}{{\ :=\ }}
\newcommand{\tq}{{\, \vert\, }}
\newcommand{\mutmut}{{\it mut.\,mut.}}
\newcommand{\Ksch}[1]{{{#1}}}
\newcommand{\Rsch}[1]{{{#1}}}
\newcommand{\ksch}[1]{{{#1}}}

\newcommand{\sfib}[1]{{\tilde{\ksch{#1}}}}
\newcommand{\quo}[1]{{\bar{\ksch{#1}}}}
\newcommand{\catM}{{\mathcal{M}}}
\newcommand{\catqA}{{\bar{\mathcal{A}}}}
\newcommand{\catqC}{{\bar{\mathcal{C}}}}
\newcommand{\catqD}{{\bar{\mathcal{D}}}}
\newcommand{\catK}{{\mathcal{K}}}
\newcommand{\catqK}{{\bar{\mathcal{K}}}}
\newcommand{\res}{{\,\operatorname{res}\,}}
\newcommand{\ind}{{\,\operatorname{ind}\,}}
\newcommand{\from}{{\ \leftarrow\ }}
\newcommand{\id}{{\, \operatorname{id}}}
\newcommand{\Hom}{{\operatorname{Hom}}}
\newcommand{\Aut}{{\operatorname{Aut}}}
\newcommand{\End}{{\operatorname{End}}}
\newcommand{\orbit}{{O}}
\newcommand{\iso}{{\ \cong\ }}
\newcommand{\proj}{{\operatorname{pr}}}
\newcommand{\conj}[2]{{\, ^{#1}\hskip-2pt{#2}}}
\newcommand{\Spec}[1]{{\operatorname{Spec}\left( #1 \right)}}

\newcommand{\Frob}{{\operatorname{Fr}}}
\newcommand{\frob}{{\operatorname{fr}}}

\newcommand{\chf}[1]{{\chi_{#1}}}
\newcommand{\trace}{{\operatorname{trace}\, }}
\newcommand{\cInd}{{\operatorname{cInd}}}

\newcommand{\Res}{{\operatorname{Res}}}

\newcommand{\er}{{\operatorname{er}}}
\newcommand{\GL}{{\operatorname{GL}}}
\newcommand{\SL}{{\operatorname{SL}}}
\newcommand{\Sp}{{\operatorname{Sp}}}
\newcommand{\SU}{{\operatorname{SU}}}

\newcommand{\fais}[1]{{\mathcal{#1}}}

\newcommand{\cind}{{\operatorname{cind}}}

\newcommand{\sgn}{{\operatorname{sgn}}}

\newcommand{\obj}{{\operatorname{obj}}}
\newcommand{\mor}{{\operatorname{mor}}}

\newcommand{\todo}[1]{\ \vspace{5mm}\par \noindent\framebox{\begin{minipage}[c]{0.95 \textwidth} \tt #1\end{minipage}} \vspace{5mm} \par}

\begin{document}

\begin{abstract}
In this paper we provide a geometric framework for the study of characters of depth-zero representations of unramified groups over local fields with finite residue fields which is built directly on Lusztig's theory of character sheaves for groups over finite fields and uses ideas due to Schneider-Stuhler. Specifically, we introduce a class of coefficient systems on Bruhat-Tits buildings of perverse sheaves sheaves on affine algebraic groups over an algebraic closure of a finite field, and to each supercuspidal depth-zero representation of an unramified $p$-adic group we associate a formal sum of these coefficient systems, called a model for the representation. Then, using a character formula due to Schneider-Stuhler and a fixed-point formula in etale cohomology we show that each model defines a distribution which coincides with the Harish-Chandra character of the corresponding representation, on the set of regular elliptic elements. The paper includes a detailed treatment of SL(2), Sp(4) and GL(n) as examples of the theory.
\end{abstract}

\maketitle


\section*{Introduction}

In their 1997 article in the Publications Math\'ematiques de l'Institut des Hautes \'Etudes Scientifiques, Peter Schneider and
Ulrich Stuhler defined a functor from the category of certain smooth representations over $\CC$ of a connected reductive $p$-adic group $\Ksch{G}(\QQ_p)$ to the category of $\Ksch{G}(\QQ_p)$-equivariant coefficient systems of vector spaces over $\CC$ and subsequently obtained a new formula for the characters of such representations on the set of regular elliptic elements of $\Ksch{G}(\QQ_p)$. 

In this article we restrict our attention to supercuspidal depth-zero representations and consider a related construction. We pass from $\QQ_p$ to an unramified closure $\QQ_p^{nr}$ and replace the category of vector spaces over $\CC$ by a triangulated category of $\ell$-adic sheaf complexes. This gives rise to a category of coefficient systems on the Bruhat-Tits building for $\Ksch{G}(\QQ_p^{nr})$ of $\ell$-adic sheaves, with $\ell \ne p$. Our category is equipped with parabolic restriction functors, an action of $\Ksch{G}(\QQ_p^{nr})$, a notion of parabolic induction and an action of Frobenius; we use all these to define \emph{Frobenius-stable admissible coefficient systems}. 

We then show that Frobenius-stable admissible coefficient systems are directly related to $\ell$-adic representations of $p$-adic groups. For example, we show there is a formal linear combination of Frobenius-stable admissible coefficient systems (an element of a Grothendieck group tensored with $\bar\QQ_\ell$) canonically associated to each supercuspidal depth-zero representation of $\Ksch{G}(\QQ_p)$; we call this a \emph{model for the representation}. We also show that each Frobenius-stable admissible coefficient system defines a distribution on the set of elliptic elements of $\Ksch{G}(\QQ_p)$. Finally, we show that the distribution associated to the model of a representation coincides with the character of the representation on the set of regular elliptic elements.

Although our focus here is on models for depth-zero supercuspidal representations, the distributions associated to admissible coefficient systems themselves are very interesting. In general, these distributions are neither orbital integrals nor characters of representations; however, they appear to generalize the distributions in \cite{W}.

We have recently found that admissible coefficient systems admit a geometric description: they may be interpreted as objects in a triangulated category of $\ell$-adic sheaf complexes on the \'etale site of a rigid analytic space associated to the group $\Ksch{G}_{\QQ_p}$. However, since the theory of derived categories of $\ell$-adic \'etale sheaves on rigid analytic spaces over $\QQ_p$ is, as far as the authors are aware, in some sense still under development, we have opted to restrict ourselves to $\ell$-adic \'etale sheaves on schemes over $\bar\FF_p$ for the moment. While the benefit of this decision is that we can provide a rigourous argument using ideas readily available in the literature, the cost of this decision is that several arguments in this paper are rather unpleasant due to the fact that we are essentially working with objects and morphisms defined by local data. The rigid analytic perspective is also the point of departure for expanding the scope of this paper to a larger class of admissible representations; in particular, we view the present paper as the depth-zero part of a naiscent theory involving $\ell$-adic sheaves on the \'etale site of a rigid analytic space which is compatible with the theory of character sheaves, via vanishing cycles functors, on the reductive quotients (over $\bar\FF_p$) of special fibres of affino\"{\i}d spaces formed from canonical integral models for parahoric subgroups.  In fact, that is where this story began, but the authors were surprised to find that much of the depth-zero story could be told without rigid analytic geometry. Hence this paper. 

\centerline{*\ *\ *}

We now describe the sections and principal results of this paper in more detail. 

In Sections~\ref{section: fundamental notions} through \ref{section: cuspidal coefficient systems}, $\KK$ denotes a field equipped with a non-trivial discrete valuation such that $\KK$ is strictly henselian and such that the residue field $\kK$ of $\KK$ is algebraically closed with non-zero characteristic. We let $\Ksch{G}$ be a connected reductive linear algebraic group over $\KK$ satisfying a hypothesis described in Section~\ref{subsection: stabilizers}. In Sections~\ref{subsection: integral models} through \ref{subsection: conjugation} we review certain basic constructions associated to parahoric subgroups of the group $\Ksch{G}(\KK)$. We denote facets of the extended Bruhat-Tits building $I(\Ksch{G},\KK)$ for $\Ksch{G}(\KK)$ by $i$, $j$, $k$ or $l$. For each such facet we consider a canonical integral model $\Rsch{G}_i$ such that $\Rsch{G}_i(\RK) = \Ksch{G}(\KK)_i$, where $\RK$ is the ring of integers in $\KK$. We are particularly interested in the maximal reductive quotient $\quo{G}_i$ of the special fibre of $\Rsch{G}_i$, which is a connected linear algebraic group over $\kK$ because of the conditions placed on $\Ksch{G}$ in Section~\ref{subsection: stabilizers}. After recalling some important facts concerning equivariant perverse sheaves in Section~\ref{subsection: equivariant perverse sheaves} we introduce cohomological parabolic induction functors on these reductive quotients in Section~\ref{subsection: parabolic induction on reductive quotients}. We then recall some important facts concerning character sheaves in Section~\ref{subsection: character sheaves}. In Section~\ref{subsection: categories} we introduce a new category, denoted $\catqD\Ksch{G}$, formed roughly by attaching the categories $D^b_c(\quo{G}_i;\EE)$ using the Bruhat order on facets of $I(\Ksch{G},\KK)$, where $D^b_c(\quo{G}_i;\EE)$ denotes the bounded derived category of constructible $\ell$-adic \'etale sheaves on $\quo{G}_i$. See Definition~\ref{definition: catqD} for the details. We end Section~\ref{section: fundamental notions} by defining an additive subcategory $\catqC\Ksch{G}$ (See Definition~\ref{definition: catqD}) of $\catqD\Ksch{G}$. Admissible coefficient systems are objects in this category with special properties.

In order to define admissible coefficient systems we first define cuspidal coefficient systems in Section~\ref{section: cuspidal coefficient systems}. We begin by defining a cohomological parabolic restriction functor $\res^{\Ksch{G}}_{\Ksch{P}}$ in Section~\ref{subsection: parabolic restriction}. Then, we describe an action of $\Ksch{G}(\KK)$ on the category $\catqC\Ksch{G}$ in Section~\ref{subsection: weakly-equivariant objects} and say that an object of $\catqC\Ksch{G}$ is \emph{weakly-equivariant} if its isomorphism class is fixed by the action of $\Ksch{G}(\KK)$. A \emph{cuspidal coefficient system} is, roughly, a weakly-equivariant simple object $C$ of $\catqC\Ksch{G}$ such that $\res^{\Ksch{G}}_{\Ksch{P}} C =0$ for every proper parabolic subgroup $\Ksch{P}$ of $\Ksch{G}$. Our first main result is Theorem~\ref{theorem: cuspidal}, which, together with Corollary~\ref{corollary: cuspidal}, provides a complete description of cuspidal coefficient systems in $\catqC\Ksch{G}$. We find that every cuspidal coefficient system may be produced, in a manner similar to compact induction, from some cuspidal character sheaf on the reductive quotient of the special fibre of the canonical integral scheme for a maximal parahoric subgroup of $\Ksch{G}(\KK)$. 

In Sections~\ref{section: admissible} we assume $\KK$ is a maximal unramified extension of a $p$-adic field. Note that such a field is strictly henselian and the residue field of that extension is an algebraic closure of a finite field. In Section~\ref{subsection: parabolic induction} we define a weakly-equivariant object $\ind^{\Ksch{G}}_{\Ksch{P}} A$ associated to any weakly-equivariant object $A$ in $\catqC\Ksch{L}$, where $\Ksch{L}$ is the Levi component for $\Ksch{P}$. Using this we define \emph{admissible coefficient systems} as those simple coefficient systems appearing in $\ind^{\Ksch{G}}_{\Ksch{P}}C$ for some parabolic subgroup $\Ksch{P}$ and some cuspidal coefficient system $C$ (see Definition~\ref{definition: admissible}).

In Sections Sections~\ref{section: frobenius} through \ref{section: examples} we fix a $p$-adic field $\Kq$ and let $\Knr$ denote a maximal unramified extension of $\Kq$. Thus, $\Knr$ plays the role of $\KK$ above. Let $\Ksch{G}_\Kq$ be a connected, quasi-split unramified linear algebraic group. Then $\Ksch{G}_\Kq \times_\Spec{\Kq} \Spec{\Knr}$ is a split connected reductive linear algebraic group over $\Knr$ and so we may let $\Ksch{G}_\Kq \times_\Spec{\Kq} \Spec{\Knr}$ play the role of $\Ksch{G}$ above.

The main idea of Section~\ref{section: frobenius} is to use the action of the Galois group $\Gal{\Knr/\Kq}$ on the extended Bruhat-Tits building $I(\Ksch{G},\Knr)$ to define a notion of (geometric) \emph{frobenius-stable} objects of $\catqC\Ksch{G}$; roughly, $C\in \obj \catqC\Ksch{G}$ is frobenius-stable if its isomorphism class is fixed by the action of Frobenius. See Proposition~\ref{proposition: frobenius} for details. In the rest of Section~\ref{section: frobenius} we briefly revisit the main ideas from Sections~\ref{section: cuspidal coefficient systems} and \ref{section: admissible} with this Galois action in mind. 

Section~\ref{section: representations} relates frobenius-stable admissible coefficient systems to depth-zero representations 
through the notion of a \emph{model for a representation}; a model is an element of the $\EE$-vector space obtained by
tensoring $\bar\QQ_\ell$ with the subgroup of the Grothendieck group for $\catqD\Ksch{G}$ generated by admissible coefficient systems (\cf Definition~\ref{definition: model}). Our second main result is Theorem~\ref{theorem: supercuspidal models}, which shows that every supercuspidal depth-zero representation admits a model.
In Section~\ref{subsection: characters} we use the character formula of \cite{SS} to associate a distribution to each admissible coefficient system. Our third main result is Theorem~\ref{theorem: character formula}, which shows that the distribution associated to a model of a depth-zero representation coincides with the character of the representation, in the sense of Harish-Chandra, on the set of regular elliptic elements of $\Ksch{G}(\Kq)$. In this way we suggest that the classification and character theory of depth-zero supercuspidal representations may be studied through the theory of admissible coefficient systems. Section~\ref{section: examples} applies the machinery of the paper to the groups of $p$-adic points on SL(2) and Sp(4) in order to illustrate this suggestion.

In summary, the main features of this paper are:
\begin{itemize}
\item
Theorem~\ref{theorem: cuspidal} and Corollary~\ref{corollary: cuspidal}, where cuspidal coefficient systems are classified;
\item
Theorem~\ref{theorem: supercuspidal models}, where models for supercuspidal depth-zero representations are
constructed;
\item
Theorem~\ref{theorem: character formula}, where we show that the distribution associated to a model of a representation equals the character of that representation on the set of regular elliptic elements;
\item
Section~\ref{section: examples}, where we give models for all supercuspidal depth-zero representations of $\SL(2)$, $\Sp(4)$ and $\GL(n)$.
\end{itemize}

\centerline{*\ *\ *}

A.-M.A. wishes to thank the University of Calgary for providing the opportunity for the authors to work together on this
project. She also wishes to thank the Pacific Institute of Mathematical Science and the Institute for Advanced Study
(Princeton) for their support. This work was presented by A.-M.A. at Ann Arbor at the invitation of Julee Kim and Gopal
Prasad and in Berlin (von Humboldt University) at the invitation of Volker Heiermann and Wilhem Zink, and also at the meeting on K-types for $p$-adic groups at M\"{u}nster organized by Guy Henniart, Philip Kutzko and Peter Schneider. She thanks Roger Carter, Ivan Fesenko and Michel Gros for their interest in this work. 

C.C. wishes to thank the \'Ecole Normale Sup\'erieure, the Centre National de la Recherche Scientifique, the Institut des Hautes \'Etudes Scientifiques and the National Science and Engineering Research Council (Canada) for support during various parts of the writing of this paper. He also wishes to thank Jonathan Korman for help and Sam Evans for his encouragement and hospitality. This paper was presented in abbreviated form by C.C. in the Summer of 2004 at Harvard at the invitation of Laurent Berger and then at a Fields Institute Workshop at the University of Ottawa at the invitation of Monica Nevins. 

Finally, both authors thank Peter Schneider for helpful conversations.

\medskip


\tableofcontents


\section{Fundamental Notions}\label{section: fundamental notions}


\subsection{Fields and algebraic groups}\label{subsection: fields and algebraic groups}

Let $\KK$ be a field equipped with a non-trivial discrete valuation, let $\RK$ be the ring of integers of $\KK$ and let $\kK$ be the residue field of $\RK$. We assume that $\KK$ is strictly henselian and that $\kK$ is algebraically closed with non-zero characteristic. Examples of such fields $\KK$ include $\QQ_p^{nr}$ (a maximal unramified extension of the field $\QQ_p$ of $p$-adic numbers) and $\bar{\mathbb{F}}_p((t))$ (formal Laurent series in $t$ with coefficients from an algebraic closure of a field with $p$ elements). Note that $\KK$ is neither complete nor locally compact.

Let $\Ksch{G}$ be a connected reductive linear algebraic group over $\KK$. We assume $\Ksch{G}$ splits over $\KK$. Let $\Ksch{G}(\KK)$ be the group of $\KK$-rational points on $\Ksch{G}$.\footnote{In Section~\ref{section: admissible} we apply the results of Sections~\ref{section: fundamental notions} and \ref{section: cuspidal coefficient systems} to the case when $\KK$ is an unramified closure of a local field with finite residue field. In Section~\ref{section: frobenius} we fix the local field, denoted $\Kq$, and consider a form $\Ksch{G}_\Kq$ for $\Ksch{G}$; thus, in Sections~\ref{section: frobenius}, \ref{section: representations} and \ref{section: examples}, $\Ksch{G}_\Kq$ is a connected reductive algebraic group over the local field $\Kq$ and $\Ksch{G}_\Kq$ splits over an unramified extension of $\Kq$.}


\subsection{Integral models}\label{subsection: integral models}

The enlarged Bruhat-Tits building for $\Ksch{G}(\KK)$ will be denoted $I(\Ksch{G},\KK)$. Recall that $I(\Ksch{G},\KK)$ is the product of the semi-simple Bruhat-Tits building for $\Ksch{G}(\KK)$ (the building for the derived group) by a real affine space. We denote polyfacets of $I(\Ksch{G},\KK)$ by $i$, $j$ or $k$ and refer to these as facets.

For each facet $i$ of $I(\Ksch{G},\KK)$, the parahoric subgroup 
of Bruhat-Tits will be denoted $\Ksch{G}(\KK)_i$ (\cf \cite[4.6.28]{BT2}). Let $\Rsch{G}_i$ denote the integral model of $\Ksch{G}$ associated to $i$ by \cite[5.1.30, 5.2.1]{BT2}; see also \cite[7.3.1]{Y1}. Thus, $\Rsch{G}_i$ is a smooth group scheme over $\RK$ equipped with an isomorphism between the generic fibre of $\Rsch{G}_i$ and $\Ksch{G}$ such that $\Rsch{G}_i(\RK)$ corresponds to $\Ksch{G}(\KK)_i$ under that isomorphism (\cf \cite[7.2]{Y1}). The special fibre of $\Rsch{G}_i$ will be denoted $\sfib{G}_i$; thus, $\sfib{G}_i = \Rsch{G}_i\times_{\Spec{\RK}}\Spec{\kK}$.  Then $\sfib{G}_i$ is a smooth connected affine group scheme over $\kK$ (\cf \cite[7.2]{Y1}). Although $\sfib{G}_i$ is reduced as a scheme, it need not be reductive as a group scheme; let ${\nu_i} : \sfib{G}_i \to \quo{G}_i$ be the maximal reductive quotient of $\sfib{G}_i$. Then $\quo{G}_i$ is a linear algebraic group over $\kK$ which is both connected and reductive. In fact,
$\quo{G}_i$ is a closed subscheme (over $\kK)$ of $\sfib{G}_i$, which is a closed subscheme of $\Rsch{G}_i$. In summary, we have the following commutative diagramme.
\begin{equation}
\xymatrix{
    \Ksch{G} \ar[r] \ar[d] & \Rsch{G}_i \ar[d] & \ar[l] \sfib{G}_i \ar[d] \ar[r]_{\nu_i} & \ar@{>.>}@/_1pc/[l] \quo{G}_i \\
    \Spec{\KK} \ar[r] & \Spec{\RK} & \ar[l] \Spec{\kK} \\
}
\end{equation}

Let $\rho_i : \Rsch{G}_i(\RK) \to \quo{G}_i(\kK)$ denote the composition of the group homomorphism $\Rsch{G}_i(\RK) \to
\Rsch{G}_i(\kK)$ defined by composition with the canonical map $\Spec{\kK} \to \Spec{\RK}$, the identification
$\Rsch{G}_i(\kK) = \sfib{G}_i(\kK)$, and the map of $\kK$-rational points $\sfib{G}_i(\kK) \to \quo{G}_i(\kK)$
induced from $\nu_i$. Observe that $\rho_i$ is a map of points; it is not a map of ringed spaces.

\begin{example}
Let $\Ksch{G} = \SL(2)_\KK$; thus, the global sections of this affine scheme are
\[
\mathcal{O}_{\Ksch{G}}(\Ksch{G}) 
= 
{\KK[X_{11}, X_{12}, X_{21},X_{22}]}/{(\det X -1)},
\]
where $\det X = X_{11}X_{22}-X_{12}X_{21} -1$. Let $i=(01)$ be the maximal facet (grande cellule) of the chamber corresponding to the Iwahori subgroup
\[
\Ksch{G}(\KK)_{(01)} = \left\{
    \begin{pmatrix}
        h_{11} & h_{12} \\
        h_{21} & h_{22}
    \end{pmatrix}
 \ \Big\vert
    \begin{array}{l}
        h_{11}, h_{12}, h_{22} \in\RK;\
        h_{21} \in \PK\\
        h_{11}h_{22} - h_{12}h_{21} = 1
    \end{array}
\right\}.
\]
In this case, $\Rsch{G}_{(01)}$ is the affine $\RK$-scheme with global sections
\[
\mathcal{O}_{\Rsch{G}_{(01)}}(\Rsch{G}_{(01)}) 
= 
{\RK[X_{11}, X_{12}, X_{21}, X_{22},X_{21}']}/{(\det X -1, X_{21}-\varpi X_{21}')},
\]
where $\varpi$ is a generator for $\PK$. (Of course, the scheme $\Rsch{G}_{(01)}$ is independent of this choice.) The isomorphism of the generic fibre of $\Rsch{G}_{(01)}$ with $\Ksch{G}$ is given by $X_{n m} \mapsto X_{n m}$, for $1\leq n,m \leq 2$. Since $\Rsch{G}_{(01)}$ is a group scheme, $\mathcal{O}_{\Rsch{G}_{(01)}}(\Rsch{G}_{(01)})$ is a Hopf algebra; the co-multiplication is given by $X_{n m} \mapsto \sum_{k} X_{nk}\otimes X_{km}$ and $X_{21}' \mapsto X_{21}'\otimes X_{12} + X_{22}\otimes X_{21}'$. The $\kK$-algebra of sections on the special fibre $\sfib{G}_{(01)}$ of $\Rsch{G}_{(01)}$ is
\begin{eqnarray*}
\mathcal{O}_{\Rsch{G}_{(01)}}(\Rsch{G}_{(01)})  \otimes \kK
    &=& {\kK[X_{11}, X_{12}, X_{22}, X_{21}']}/{(X_{11}X_{22}-1)},
\end{eqnarray*}
with co-multiplication given by
\begin{eqnarray*}
X_{11} &\mapsto& X_{11}\otimes X_{11}\\
X_{12} &\mapsto& X_{11}\otimes X_{12} + X_{12}\otimes X_{22}\\
X_{21}' &\mapsto& X_{21}'\otimes X_{11} + X_{22}\otimes X_{21}'\\
X_{22} &\mapsto& X_{22}\otimes X_{22}.
\end{eqnarray*}
Evidently, $\sfib{G}_{(01)}$ is not reductive. The reductive quotient of $\sfib{G}_{(01)}$ is $\GL(1)_\kK$, and the map $\nu_{(01)} : \sfib{G}_{(01)} \to \quo{G}_{(01)}$ is induced from the inclusion
\begin{eqnarray*}
{\kK[X_{11}, X_{22}]}/{(X_{11}X_{22}-1)} &\hookrightarrow&
{\kK[X_{11}, X_{12}, X_{22}, X_{21}']}/{(X_{11}X_{22}-1)}.
\end{eqnarray*}
In this example, $\rho_{(01)} : \Rsch{G}_{(01)}(\RK) \to \quo{G}_{(01)}(\kK)$ is given by
\[
\rho_{(01)}\begin{pmatrix}
        h_{11} & h_{12} \\
        h_{21} & h_{22}
    \end{pmatrix} = \bar{h}_{11},
\]
where $\bar{h}_{11}$ is the image of $h_{11}$ under the canonical map $\RK \to \kK$.
\end{example}

\begin{example}
Continuing with $\Ksch{G} = \SL(2)_\KK$, let $(0)$ and $(1)$ be the vertices in the closure of the facet considered above. Let $\Rsch{G}_{(0)}$ be the model for $\Ksch{G}$ with global section given by
\[
\mathcal{O}_{\Rsch{G}_{(0)}}(\Rsch{G}_{(0)}) 
= 
{\RK[X_{11}, X_{12}, X_{21},X_{22}]}/{(\det X -1)}.
\]
As above, the isomorphism of the generic fibre of $\Rsch{G}_{(0)}$ with $\Ksch{G}$ is the obvious one. On the other hand, $\Rsch{G}_{(1)}$ is the integral affine scheme with global sections
\[
\mathcal{O}_{\Rsch{G}_{(1)}}(\Rsch{G}_{(1)}) 
= 
\frac{\RK[X_{11}, X_{12}, X_{21},X_{22},X_{12}', X_{21}']}{(\det X - 1, \varpi X_{12}-X_{12}', X_{21} - \pi X_{21}')},
\]
where $\varpi$ is a uniformizer for $\KK$. (As above, the scheme $\Rsch{G}_{(1)}$ is independent of this choice.) The isomorphism from the generic fibre of $\Rsch{G}_{(1)}$ to $\Ksch{G}$ is determined by $X_{n m} \mapsto X_{n m}$ for $1 \leq n, m \leq 2$. In both cases, the special fibre is $\SL(2)_\kK$; since this group scheme is reductive, the maps $\nu_{(0)}$ and $\nu_{(1)}$ are identities.
\end{example}


\subsection{Stabilizers}\label{subsection: stabilizers}

When the time comes to relate admissible coefficient systems to characters of depth-zero representations of $p$-adic groups we will see that it is natural to study stabilizers of facets rather than parahoric subgroups. The stabiliser of any facet of $I(\Ksch{G},\KK)$ under the action of $\Ksch{G}(\KK)$ is a compact group (recall that $I(\Ksch{G},\KK)$ refers to the {\em enlarged} Bruhat-Tits building) which admits a canonical smooth integral model (\cf \cite[9.3.2]{Y1}). However, in general, the maximal reductive quotient of the special fibre of this integral model need not be connected, so we cannot use the theory of character sheaves as developed in \cite{CS}. We therefore now impose a condition on the groups $\Ksch{G}$ that we study: {we assume that the stabilizer of each facet in $I(\Ksch{G},\KK)$ is a parahoric subgroup, and further that the reductive quotient of the special fibre of the canonical integral model of that parahoric subgroup is connected, as an algebraic group over $\kK$; we also demand that the same property hold for all cuspidal Levi subgroups of $\Ksch{G}$ (\cf Definition~\ref{definition: cuspidal Levi})}. When $\Ksch{G}$ is simply connected, this condition is satisfied (\cf \cite[3.5.2]{T}); it also holds for general linear groups and symplectic groups. 

Remarkably, Lusztig has recently extended the definition of character sheaves to the disconnected group case so there is good reason to expect the main results of this paper can be extended, \mutmut, to a larger class of groups.


\subsection{Restriction between reductive quotients}\label{subsection: restriction between reductive quotients}

Let $i$ and $j$ be facets of $I(\Ksch{G},\KK)$ such that $i\leq j$ in the Bruhat order. Let $\Rsch{f}_{i\leq j} : \Rsch{G}_j \to \Rsch{G}_i$ be the morphism of group schemes over $\RK$ obtained by extending the identity morphism $\id_\Ksch{G}$ in the category of group schemes over $\RK$ (\cf \cite[6.2]{Lan1}). By restriction to special fibres, this defines a morphism $\sfib{f}_{i\leq j} : \sfib{G}_j \to \sfib{G}_i$ of group schemes over $\kK$, making the following diagramme commute.
\begin{equation}
	\xymatrix{
	\Rsch{G}_j \ar[r]^{\Rsch{f}_{i\leq j}} & \Rsch{G}_i \\
	\sfib{G}_j \ar[u] \ar[r]_{\sfib{f}_{i\leq j}} & \sfib{G}_i \ar[u]
	}
\end{equation}
In fact, this diagramme is cartesian. Let $\sfib{G}_{i\leq j}$ denote the schematic image of $\sfib{f}_{i\leq j}$ in $\sfib{G}_i$. Let $\quo{G}_{i\leq j}$ be the schematic image of $\sfib{G}_{i\leq j}$ under $\nu_i$ and let $\nu_{i\leq j}$ denote the restriction of $\nu_i$ to $\sfib{G}_{i\leq j}$. Next, let 
\begin{equation}\label{equation: iso thm}
	\sfib{f}_{i\leq j} = \ksch{h}_{i\leq j} \circ \ksch{g}_{i\leq j}
\end{equation}
be the factorization given by the Isomorphism Theorem. The kernel of $\sfib{f}_{i\leq j}$, which equals the kernel of $\ksch{g}_{i\leq j}$, is contained in the kernel of $\nu_j$; thus, $\ksch{g}_{i\leq j}$ factors through $\nu_j$ to give a map $\sfib{G}_{i\leq j} \to \quo{G}_j$; since the kernel of $\nu_{i\leq j}$ is contained in the kernel of this new map, it too factors, this time through $\nu_{i\leq j}$, thus defining $\ksch{r}_{i\leq j} : \quo{G}_{i\leq j} \to \quo{G}_j$. Notice that
\begin{equation}\label{equation: nuj}
\nu_j = \ksch{r}_{i\leq j} \circ \nu_{i\leq j} \circ \ksch{g}_{i\leq j}.
\end{equation}
Let $s_{i\leq j} : \quo{G}_{i\leq j} \to \quo{G}_i$ be the obvious inclusion; this is an affine closed immersion. By
\cite[9.22]{Lan1}, $\ksch{s}_{i\leq j} : \quo{G}_{i\leq j} \to \quo{G}_i$ is
a parabolic subgroup with Levi component $\quo{G}_j$ given by the reductive quotient map $\ksch{r}_{i\leq j} : \quo{G}_{i\leq j} \to \quo{G}_j$, which is a smooth projective map; we also have $\quo{G}_i = \sfib{G}_i \times_{\quo{G}_{i\leq j}} \sfib{G}_{i\leq j}$. In summary we have the following commutative diagramme, in which the square on the bottom right is cartesian.
\begin{equation}\label{equation: restriction diagramme}
\xymatrix{
	\sfib{G}_j \ar[rr]^{\sfib{f}_{i\leq j}} \ar[dr]^{\ksch{g}_{i\leq j}} \ar[d]_{\nu_j} && \sfib{G}_i \ar[d]^{\nu_i}\\
 	\quo{G}_j & \sfib{G}_{i\leq j} \ar[ur]^{\ksch{h}_{i\leq j}} \ar[d]^{\nu_{i\leq j}} & \quo{G}_i\\
	& \quo{G}_{i\leq j} \ar[ul]^{\ksch{r}_{i\leq j}} \ar[ur]_{\ksch{s}_{i\leq j}} & \\
}
\end{equation}

Recall the derived category $D^b_c(\quo{G}_i,\bar\QQ_\ell)$ of cohomologically bounded constructible $\ell$-adic sheaves with $\ell\ne p$, introduced in \cite[1.1.1-1.1.5]{D2} and \cite[2.2.9, 2.2.14, 2.2.18]{BBD} (\cf \cite[expos\'es VI, V, XV]{SGA5}). 
We will follow the notational conventions of \cite{BBD} regarding derived functors.

\begin{definition}\label{definition: restriction on quotients}
Let $i$ and $j$ be facets of $I(\Ksch{G},\KK)$ such that $i \leq j$. Define $\res_{i\leq j} :  D^b_c(\quo{G}_i;\EE) \to
 D^b_c(\quo{G}_j;\EE)$ by
 \begin{equation}
	\res_{i\leq j} = {\ksch{r}_{i\leq j}}_!\  \ksch{s}_{i\leq j}^*\ (d_{i\leq j}),
\end{equation}
where $(d_{i\leq j})$ denotes Tate twist by $d_{i\leq j} = \dim\ker\ksch{r}_{i\leq j}$.
\end{definition}

\begin{remark}\label{remark: specific}
Thus, $\res_{i\leq j } = \res^{\quo{G}_i}_{\quo{G}_{i\leq j}}$, where the right-hand side refers to the parabolic restriction functor defined in \cite[Sect.3.8]{CS}.  We will sometimes write $\res^{\quo{G}_i}_{\quo{G}_j}$ for $\res_{i\leq j}$ to emphasize the fact that it is a functor from $D^b_c(\quo{G}_i;\EE)$ to $D^b_c(\quo{G}_j;\EE)$. It must be understood that the definition of the functor makes reference to a \emph{specific} parabolic subgroup of $\quo{G}_i$ with Levi component $\quo{G}_j$.
Note also that $\res_{i\leq i}$ is an identity functor.
 \end{remark}


\begin{proposition}\label{proposition: transitive restriction on quotients}
If $i, j, k, l$ are facets of $I(\Ksch{G},\KK)$ such that $i\leq j\leq k\leq l$ then there are canonical isomorphisms of functors
\[
\res_{i\leq j\leq k} :\res_{j\leq k} \res_{i\leq j} \to \res_{i\leq k}
\]
such that the diagramme
\[
	\xymatrix{
	\ar[d]_{\res_{k\leq l} \res_{i\leq j\leq k} } 
\res_{k\leq l} \res_{j\leq k}\res_{i\leq j}  
	\ar[rrr]^{\res_{j\leq k\leq l} \res_{i\leq j} } 
&&& 
	\ar[d]^{\res_{i\leq j\leq l} } 
\res_{j\leq l} \res_{i\leq j} 
\\
\res_{k\leq l} \res_{i\leq k}  
\ar[rrr]_{\res_{i\leq k\leq l}} 
&&&
\res_{i\leq l}
	}
\]
commutes.
\end{proposition}
 
\begin{proof}
Observe that $i \leq j \leq k\leq l$ implies there is an apartment containing all of $i$, $j$, $k$ and $l$. 
We begin by defining $\res_{i\leq j\leq k}$. Let $\ksch{r}_{i\leq j\leq k} : \quo{G}_{i\leq k} \to \quo{G}_{j\leq k}$ and $\ksch{s}_{i\leq j\leq k} : \quo{G}_{i\leq k} \hookrightarrow \quo{G}_{i\leq j}$ be \emph{the} pull-back of $\ksch{r}_{i\leq j} : \quo{G}_{i\leq j} \to \quo{G}_j$ and $\ksch{s}_{j\leq k} : \quo{G}_{j\leq k} \to \quo{G}_k$ with domain $\quo{G}_{i\leq k}$; in particular,
\begin{equation}\label{equation: rijk}
  \ksch{r}_{i\leq k} = \ksch{r}_{j\leq k} \circ \ksch{r}_{i\leq j\leq k}
\end{equation}
and 
\begin{equation}\label{equation: sijk}
\ksch{s}_{i\leq k} = \ksch{s}_{i\leq j} \circ \ksch{s}_{i\leq j\leq k}.
\end{equation}
See \cite[Prop.9.22]{Lan1} for the existence of such a pull-back. The situation is summarized by the following diagramme, in which the square is cartesian and all triangles commute.
\begin{equation}\label{equation: trq.0}
	\xymatrix{
&& \ar@/_2pc/[ddll]_{\ksch{s}_{i\leq k}} \ar@{.>}[dl]^{\ksch{s}_{i\leq j\leq k}} \quo{G}_{i\leq k} \ar@{.>}[dr]_{\ksch{r}_{i\leq j\leq k}} \ar@/^2pc/[ddrr]^{\ksch{r}_{i\leq k}} && \\
& \ar[dl]^{\ksch{s}_{i\leq j}}  \quo{G}_{i\leq j} \ar[dr]_{\ksch{r}_{i\leq j}}  &&  \ar[dl]^{\ksch{s}_{j\leq k}} \quo{G}_{j\leq k} \ar[dr]_{\ksch{r}_{j\leq k}}  & \\	
\quo{G}_{i} && \quo{G}_j && \quo{G}_k \\
	}
\end{equation}
Observe that all maps $\ksch{r}_{\cdot}$ are smooth projective and all maps $\ksch{s}_{\cdot}$ are affine closed immersions. 

To define $\res_{i\leq j\leq k}$ we begin by observing that $d_{i\leq k} = d_{i\leq j}+ d_{j\leq k}$ and that Tate twists commute with everything below. Thus,
\begin{eqnarray*}
    \res_{j\leq k}
    \res_{i\leq j}
 &=& {\ksch{r}_{j\leq k}}_!\ \ksch{s}_{j\leq k}^*\ (d_{j\leq k})\ {\ksch{r}_{i\leq j}}_!\ \ksch{s}_{i\leq j}^*\ (d_{i\leq j})\\
 &=& {\ksch{r}_{j\leq k}}_!\ \ksch{s}_{j\leq k}^*\ {\ksch{r}_{i\leq j}}_!\ \ksch{s}_{i\leq j}^*\ (d_{i\leq k}).
\end{eqnarray*}
Applying the smooth base-change theorem for direct images with compact supports for $\ell$-adic sheaves (see \cite[Expos\'e XVII, \S 5.2]{SGA4}) to the cartesian square in Diagramme~\ref{equation: trq.0}, it follows that the base-change natural transformation
\begin{equation}\label{equation: trq.1}
\varphi_{i\leq j\leq k} : \ksch{s}_{j\leq k}^*\ {\ksch{r}_{i\leq j}}_! \to {\ksch{r}_{i\leq j\leq k}}_!\ \ksch{s}_{i\leq j\leq k}^*
\end{equation}
is an isomorphism of functors. Thus,
\begin{equation}\label{equation: trq.2}
	\xymatrix{
{\ksch{r}_{j\leq k}}_!\ \ksch{s}_{j\leq k}^*\ {\ksch{r}_{i\leq j}}_!\ \ksch{s}_{i\leq j}^*\ (d_{i\leq k}) 
	\ar[d]^{{\ksch{r}_{j\leq k}}_! \varphi_{i\leq j\leq k} \ksch{s}_{i\leq j}^*\ (d_{i\leq k}) }
	\\
{\ksch{r}_{j\leq k}}_!\ {\ksch{r}_{i\leq j\leq k}}_!\ \ksch{s}_{i\leq j\leq k}^*\ \ksch{s}_{i\leq j}^*\ (d_{i\leq k})
	}
\end{equation}
is a natural isomorphism. Let $\rho_{i\leq j\leq k} : {\ksch{r}_{j\leq k}}_!\ {\ksch{r}_{i\leq j\leq k}}_! \to {\ksch{r}_{i\leq k}}_!$ be the natural isomorphism determined by Equation~\ref{equation: rijk}; these isomorphisms satisfy a cocycle condition (see \cite[Expos\'e XVII, Thm~5.1.8(a)(i)]{SGA4}). Likewise, let $\sigma_{i\leq j\leq k} : \ksch{s}_{i\leq j\leq k}^*\ \ksch{s}_{j\leq k}^*  \to \ksch{s}_{i\leq k}^*$ be the natural isomorphism determined by Equation~\ref{equation: sijk}; these isomorphisms satisfy the analogous cocycle condition. Now, the following diagramme commutes.
\begin{equation}\label{equation: trq.3}
	\xymatrix{
{\ksch{r}_{j\leq k}}_!{\ksch{r}_{i\leq j\leq k}}_!\ksch{s}_{i\leq j\leq k}^*\ksch{s}_{i\leq j}^*(d_{i\leq k})
&& 
	\ar[ll]_{\rho_{i\leq j\leq k}\ksch{s}_{i\leq j\leq k}^*\ksch{s}_{i\leq j}^* (d_{i\leq k})}
 {\ksch{r}_{i\leq k}}_!\ksch{s}_{i\leq j\leq k}^*\ksch{s}_{i\leq j}^*(d_{i\leq k})
 \\
&& 
\\
	\ar[uu]^{{\ksch{r}_{j\leq k}}_!{\ksch{r}_{i\leq j\leq k}}_!\sigma_{i\leq j\leq k}(d_{i\leq k})} 
{\ksch{r}_{j\leq k}}_!{\ksch{r}_{i\leq j\leq k}}_!\ksch{s}_{i\leq k}^*(d_{i\leq k}) 
&&  
{\ksch{r}_{i\leq k}}_!\ksch{s}_{i\leq k}^*(d_{i\leq k}) 
 	\ar[uu]_{{\ksch{r}_{i\leq k}}_!\sigma_{i\leq j\leq k}(d_{i\leq k})}
 	\ar[ll]^{\rho_{i\leq j\leq k}\ksch{s}_{j\leq k}^* (d_{i\leq k})}
\\
	}
\end{equation}
Define $\res_{i\leq j\leq k}$ by composing Diagramme~\ref{equation: trq.2} with Diagramme~\ref{equation: trq.3} in the obvious manner. It is clearly a natural isomorphism as it is defined by composing natural isomorphisms.

Having defined $\res_{i\leq j\leq k}$ we now turn to the remaining part of Proposition~\ref{proposition: transitive restriction on quotients}. Using the same procedure as above, let $\ksch{r}_{j\leq k\leq l} : \quo{G}_{j\leq l} \to \quo{G}_{j\leq k}$ and $\ksch{s}_{j\leq k\leq l} : \quo{G}_{j\leq l} \hookrightarrow \quo{G}_{j\leq k}$ be the pull-back of $\ksch{r}_{j\leq k} : \quo{G}_{j\leq k} \to \quo{G}_k$ and $\ksch{s}_{k\leq l} : \quo{G}_{k\leq l} \to \quo{G}_l$ with domain $\quo{G}_{j \leq l}$; in particular, $\ksch{r}_{j\leq l} = \ksch{r}_{k\leq l} \circ \ksch{r}_{j\leq k\leq l}$ and  $\ksch{s}_{j\leq l} = \ksch{s}_{j\leq k} \circ \ksch{s}_{j\leq k\leq l}$.  Likewise, let $\ksch{r}_{i\leq j\leq k\leq l} : \quo{G}_{i\leq l} \to \quo{G}_{j\leq l}$ and $\ksch{s}_{i\leq j\leq k\leq l} : \quo{G}_{i\leq l} \hookrightarrow \quo{G}_{i\leq k}$ be the pull-back of $\ksch{r}_{i\leq j\leq k} : \quo{G}_{i\leq k} \to \quo{G}_{j\leq k}$ and $\ksch{s}_{j\leq k\leq l} : \quo{G}_{j\leq l} \to \quo{G}_{j\leq k}$ with domain $\quo{G}_{i\leq l}$;  in particular, $\ksch{r}_{i\leq l} = \ksch{r}_{j\leq k} \circ \ksch{r}_{i\leq j\leq k} \circ \ksch{r}_{i\leq j\leq k\leq l} $
and $\ksch{s}_{i\leq k} = \ksch{s}_{i\leq j} \circ \ksch{s}_{i\leq j\leq k} \circ \ksch{s}_{i\leq j\leq k}$. The situation is summarized by the top part (the upper twelve arrows) of Diagramme~\ref{equation: trq.4}, in which all squares are cartesian. As above, observe that all maps $\ksch{r}_{\cdot}$ are smooth projective and all maps $\ksch{s}_{\cdot}$ are affine closed immersions. The bottom part (the lower six arrows) of Diagramme~\ref{equation: trq.4} is obtained by pushing-out, which is possible exactly because the maps $\ksch{r}_{\cdot}$ are smooth projective and all maps $\ksch{s}_{\cdot}$ are affine closed immersions!
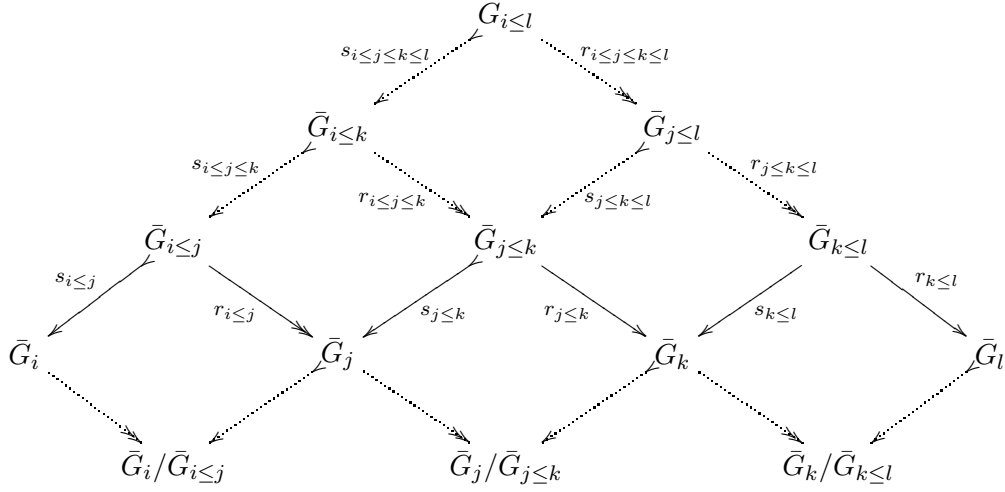
\begin{equation}\label{equation: trq.4}
	\xymatrix{
&&& \ar@{>.>}[dl]_{\ksch{s}_{i\leq j\leq k\leq l}} \quo{G}_{i\leq l} \ar@{.>>}[dr]^{\ksch{r}_{i\leq j\leq k\leq l}}  &&&\\
&&  \ar@{>.>}[dl]_{\ksch{s}_{i\leq j\leq k}} \quo{G}_{i\leq k} \ar@{.>>}[dr]_{\ksch{r}_{i\leq j\leq k}}  && 
\ar@{>.>}[dl]^{\ksch{s}_{j\leq k\leq l}} \quo{G}_{j\leq l} \ar@{.>>}[dr]^{\ksch{r}_{j\leq k\leq l}}  &&\\
& \ar@{>->}[dl]_{\ksch{s}_{i\leq j}}  \quo{G}_{i\leq j} \ar@{->>}[dr]_{\ksch{r}_{i\leq j}}  &&  \ar@{>->}[dl]^{\ksch{s}_{j\leq k}} \quo{G}_{j\leq k} \ar[dr]_{\ksch{r}_{j\leq k}}  && \ar[dl]^{\ksch{s}_{k\leq l}} \quo{G}_{k\leq l} \ar[dr]^{\ksch{r}_{k\leq l}} & \\	
\quo{G}_{i} \ar@{.>>}[dr] && \ar@{>.>}[dl] \quo{G}_j \ar@{.>>}[dr] && \ar@{>.>}[dl] \quo{G}_k \ar@{.>>}[dr] && \ar@{>.>}[dl] \quo{G}_l \\
& \quo{G}_{i}/\quo{G}_{i\leq j} && \quo{G}_{j}/\quo{G}_{j\leq k}  && \quo{G}_{k}/\quo{G}_{k\leq l} & \\
	}
\end{equation}
The cocycle relation for the restriction functors is obtained by repeated application of \cite[Expos\'e XVII, Thm~4.4]{SGA4} to Diagramme~\ref{equation: trq.4}.
\end{proof}


\subsection{Parabolic restriction on the level of reductive quotients}\label{subsection: parabolic restriction on quotients}

In Section~\ref{subsection: parabolic restriction} we will need the following consequence of Proposition~\ref{proposition: transitive restriction on quotients}. Let $\Ksch{P} \subseteq \Ksch{G}$ be a parabolic subgroup with reductive quotient $\Ksch{L}$.  We fix an imbedding of buildings $I(\Ksch{L},\KK) \to I(\Ksch{G},\KK)$ (\cf \cite{Lan2}). We will write $i_\Ksch{G}$ for the image of a facet $i$ of $I(\Ksch{L},\KK)$ under this embedding. Let $i$ be any facet of $I(\Ksch{L},\KK)$. Then $\quo{L}_i$ is a Levi subgroup of $\quo{G}_{i_\Ksch{G}}$. By \cite[9.22]{Lan1} there is a unique facet $i_\Ksch{P}$ in $I(\Ksch{G},\KK)$ such that $i_\Ksch{G} \leq i_\Ksch{P}$ and $\quo{L}_i = \quo{G}_{i_\Ksch{P}}$ and $\quo{G}_{i_\Ksch{G} \leq i_\Ksch{P}}$ is the schematic intersection of $\quo{G}_{i_\Ksch{G}}$ with $\Ksch{P}$ in $\Rsch{G}_{i_\Ksch{G}}$. (See Section~\ref{subsection: restriction between reductive quotients} for the definition of $\quo{G}_{i_\Ksch{G} \leq i_\Ksch{P}}$.) 

\begin{lemma}\label{lemma: parabolic restriction on quotients}
Let $\Ksch{P}$ be a parabolic subgroup of $\Ksch{G}$ with levi component $\Ksch{L}$. With notation as above, there is an isomorphism of functors
\[
\res^{\Ksch{P}}_{i\leq j} :  \res_{j_\Ksch{P}\leq i_\Ksch{P}} \res_{i_\Ksch{G}\leq i_\Ksch{P}} \to \res_{j_\Ksch{G}\leq j_\Ksch{P}} \res_{i_\ksch{G}\leq j_\Ksch{G}} 
\]
such that
\[
	\xymatrix{
	\ar[rrrr]^{ \res_{i_\Ksch{P}\leq j_\Ksch{P} \leq k_\Ksch{P}} \res_{i_\Ksch{G}\leq i_\Ksch{P}}}
\res_{j_\Ksch{P}\leq k_\Ksch{P}} \res_{i_\ksch{P}\leq j_\ksch{P}} \res_{i_\Ksch{G}\leq i_\Ksch{P}} 
	\ar[d]^{\res_{j_\Ksch{P}\leq k_\Ksch{P}} \res^{\Ksch{P}}_{i\leq j}}
&&&&
\res_{i_\ksch{P}\leq k_\Ksch{P}} \res_{i_\Ksch{G}\leq i_\Ksch{P}} 
	\ar[dd]_{\res^{\Ksch{P}}_{i\leq k}}
\\ 
\res_{j_\Ksch{P}\leq k_\Ksch{P}} \res_{j_\ksch{G}\leq j_\ksch{P}} \res_{i_\Ksch{G}\leq j_\Ksch{G}} 
	\ar[d]^{\res^{\Ksch{P}}_{j\leq k} \res_{j_\Ksch{P}\leq k_\Ksch{P}}}
&&&&
\\ 
	\ar[rrrr]_{\res_{k_\Ksch{G}\leq k_\Ksch{P}} \res_{i_\Ksch{G}\leq j_\Ksch{G}\leq k_\ksch{G}}} 
\res_{k_\Ksch{G}\leq k_\Ksch{P}} \res_{j_\ksch{G}\leq k_\ksch{G}} \res_{i_\Ksch{G}\leq j_\Ksch{G}} 
&&&&
\res_{k_\Ksch{G}\leq k_\Ksch{P}} \res_{i_\Ksch{G}\leq k_\ksch{G}} 
\\
	}
\]
commutes for all facets $i$, $j$ and $k$ of $I(\Ksch{L},\KK)$ such that $i\leq j\leq k$.
\end{lemma}

\begin{proof}
As the notation perhaps suggests, the natural transformation $\res^{\Ksch{P}}_{i\leq j}$ is defined using Definition~\ref{definition: restriction on quotients}; specifically,
\begin{equation}
\res^{\Ksch{P}}_{i\leq j} \ceq \res_{i_\Ksch{G}\leq j_\Ksch{G}\leq j_\ksch{P}}^{-1} \circ \res_{i_\Ksch{G}\leq i_\Ksch{P}\leq j_\Ksch{P}} .
\end{equation}
This is clearly an isomorphism of functors. The property appearing in Lemma~\ref{lemma: parabolic restriction on quotients} follows from Proposition~\ref{proposition: transitive restriction on quotients}.
\end{proof}


\subsection{Conjugation}\label{subsection: conjugation}

Let $\Ksch{m} : \Ksch{G}\times\Ksch{G} \to \Ksch{G}$ be conjugation over $\KK$. Recall that the Bruhat-Tits building
$I(\Ksch{G},\KK)$ is equipped with an action of $\Ksch{G}(\KK)$ which we indicate by
\begin{eqnarray*}
\Ksch{G}(\KK) \times I(\Ksch{G},\KK) &\to& I(\Ksch{G},\KK)\\
    (g,i) &\mapsto& gi.
\end{eqnarray*}
We will also write $ig$ for $g^{-1}i$. 

Fix an element $g$ of $\Ksch{G}(\KK)$ and let $\Ksch{m}(g) : \Ksch{G} \to \Ksch{G}$ be the morphism given by $\Ksch{m}(g)(h) = \Ksch{m}(g,h)$ for $h \in \Ksch{G}(\KK)$ (recall that $\kK$ is algebraically closed, see Section~\ref{subsection: fields and algebraic groups}). Fix a facet $i$ and recall that $\Rsch{G}_i$ and $\Rsch{G}_{gi}$ are smooth integral models of $\Ksch{G}$. Since $\Ksch{m}(g)(\Rsch{G}_i(\RK)) = \Rsch{G}_{gi}(\RK)$, it follows from the Extension Principle (\cf \cite[1.7]{BT2}) that the isomorphism $\Ksch{m}(g):\Ksch{G} \to \Ksch{G}$ of group schemes over $\KK$ extends to an isomorphism $\Rsch{m}(g)_{i}: \Rsch{G}_i \to \Rsch{G}_{gi}$ of group schemes over $\RK$. Restricting to special fibres, $\Rsch{m}(g)_i$ defines an isomorphism $\quo{m}(g)_i : \quo{G}_i \to \quo{G}_{gi}$ of reductive quotients.
Restricting to case when $g$ in an element of $\Rsch{G}_i(\RK)$ we obtain a family of isomorphisms 
\[
\quo{m}(g)_i : \quo{G}_i \to \quo{G}_i
\] 
which together define conjugation $\quo{m}_i : \quo{G}_i \times \quo{G}_i \to \quo{G}_i$ on the level of reductive quotients. If $h$ is an element of $\Rsch{G}_i(\RK)$ then $\quo{m}_i(\rho_i(h)) = \quo{m}(h)_i$, with $\rho_i$ as
defined in Section~\ref{subsection: integral models}.

\begin{lemma}\label{lemma: conjugation and restriction}
Let $g$ be an element of $\Ksch{G}(\KK)$ and let $i,j$ be facets of $I(\Ksch{G},\KK)$ with $i\leq j$. Then there is an isomorphism of functors in $D^b_c(\quo{G}_{g i};\EE)$
\[
\res^g_{i\leq j} : \res_{i\leq j}\ \quo{m}(g)_{i}^* \iso \quo{m}(g)_{j}^*\ \res_{gi\leq gj}
\] 
such that
\[
	\xymatrix{
	\ar[rrrr]^{\res_{i\leq j\leq k}\quo{m}(g)_i^*}   
\res_{j\leq k} \res_{i\leq j}\ \quo{m}(g)_i^* 	
	\ar[d]^{\res_{j\leq k} \res^g_{i\leq j}}
&&&&
\res_{i\leq k}\  \quo{m}(g)_i^*\  
	\ar[dd]_{\res^g_{i\leq k}} 
\\ 
\res_{j\leq k}\ \quo{m}(g)_j^*\ \res_{gi\leq gj} 	
	\ar[d]^{\res^g_{j\leq k} \res_{gi\leq gj}}
&&&&
\\ 
\quo{m}(g)_k^*\ \res_{gj\leq gk} \res_{gi\leq gj}   
		\ar[rrrr]_{\quo{m}(g)_k^*\ \res_{gi\leq gj\leq gk} }
&&&&
\quo{m}(g)_k^*\ \res_{gi\leq gk}\\
	}
\]
commutes for all facets $i$, $j$ and $k$ of $I(\Ksch{G},\KK)$ such that $i\leq j\leq k$.
\end{lemma}

\begin{proof}
Observe that 
\begin{equation}\label{equation: cr.0}
\res_{i\leq j}\ {\quo{m}(g)_i}^*= {\ksch{r}_{i\leq j}}_!\ {\ksch{s}_{i\leq j}}^*\ {\quo{m}(g)_i}^*\ (d_{i\leq j}),
\end{equation}
by Definition~\ref{definition: restriction on quotients}. Now, consider the following commutative diagramme, where $\quo{m}(g)_{i\leq j}$ is the isomorphism of special fibres obtained by restricting $\quo{m}(g)_j$ to $\quo{G}_{i\leq j}$.
\begin{equation}\label{equation: cr.1}
	\xymatrix{
\quo{G}_i \ar[d]_{\quo{m}(g)_i} 	& \ar[l]_{\ksch{s}_{i\leq j}} \quo{G}_{i\leq j} \ar[r]^{\ksch{r}_{i\leq j}} \ar[d]^{\quo{m}(g)_{i\leq j}}  & \quo{G}_j \ar[d]^{\quo{m}(g)_j} \\
\quo{G}_{gi} & \ar[l]^{\ksch{s}_{gi\leq gj}} \quo{G}_{gi\leq gj} \ar[r]_{\ksch{r}_{gi\leq gj}}  &  \quo{G}_{gj}
	}
\end{equation}
Since the left-hand square in Diagramme~\ref{equation: cr.1} commutes (by construction) we have the following natural isomorphisms in $D(\quo{G}_{gi};\EE)$.
\begin{equation}\label{equation: cr.2}
{\ksch{s}_{i\leq j}}^*\ {\quo{m}(g)_i}^*  
{\iso} 
 (\quo{m}(g)_i \circ s_{i\leq j})^*  
{=}
(\ksch{s}_{gi\leq gj} \circ \quo{m}(g)_{i\leq j})^*  
{\iso}
{\quo{m}(g)_{i\leq j}}^*\ {\ksch{s}_{gi\leq gj}}^*
\end{equation}
Applying the smooth base-change theorem for direct images with compact support for $\ell$-adic sheaves (see \cite[Expos\'e XVII, \S 5.2]{SGA4}, see also \cite[Thm~6.3 (c)]{E}) to the right-hand square in Diagramme~\ref{equation: cr.1} (which is indeed cartesian) it follows that the base-change natural transformation
\begin{equation}\label{equation: cr.3}
\varphi^g_{i\leq j} :   \quo{m}(g)_j^*\ {\ksch{r}_{gi\leq gj}}_! \to {\ksch{r}_{i\leq j}}_!\ \quo{m}(g)_{i\leq j}^*
\end{equation}
is an isomorphism of functors. Since 
\begin{equation}\label{equation: cr.4}
 {\ksch{r}_{gi\leq gj}}_!\ {\ksch{s}_{gi\leq gj}}^*\ (d_{gi\leq gj}) = \res_{gi\leq gj},
\end{equation}
by Definition~\ref{definition: restriction on quotients} again, we define $\res^g_{i\leq j}$ by composing the isomorphisms appearing in Equations~\ref{equation: cr.0}, \ref{equation: cr.2}, \ref{equation: cr.3} and \ref{equation: cr.4} in the obvious manner. The property appearing in Lemma~\ref{lemma: conjugation and restriction} is now a direct result Proposition~\ref{proposition: transitive restriction on quotients} (which in turn follows from \cite[Expos\'e XVII, \S 5]{SGA4}).
\end{proof}


\subsection{Equivariant perverse sheaves}\label{subsection: equivariant perverse sheaves}

In this Subsection we discuss a fundamental result concerning the category of equivariant perverse sheaves which plays a key role in the definition of parabolic induction \emph{as a functor}. 

Let $\ksch{X}$ be an algebraic variety over $\kK$ and let $\catM\ksch{X}$ denote the category of perverse sheaves on $\ksch{X}$.  Let $\ksch{m} : \ksch{P}\times \ksch{X} \to \ksch{X}$ be an action of a {connected} algebraic group on $\ksch{X}$ over $\kK$. Recall from \cite[\S0]{L0} that a perverse sheaf $F$ on $\ksch{X}$ is {equivariant} if there is an isomorphism 
	\begin{equation}\label{equation: equivariant}
	\mu_{F} : \ksch{m}^* F \to \proj^* F
	\end{equation}
in $D^b_c(\ksch{P}\times\ksch{X};\EE)$ such that $\ksch{e}^*\mu_F = \id_F$, where  $\ksch{e} : \ksch{X} \to \ksch{P}\times\ksch{X}$ is defined by $x \mapsto (1,x)$ and $\proj : \ksch{P}\times \ksch{X} \to \ksch{X}$ is projection onto the second component. As observed in \cite[\S0]{L0}, if $F$ is an equivariant perverse sheaf, then there is exactly one such isomorphism $\mu_F$.

Fix $h\in \ksch{P}(\kK)$ and let $\ksch{e}_x : \ksch{X} \to \ksch{P} \times \ksch{X}$ be the morphism determined by $\ksch{e}_x(y) = (x,y)$ (recall that $\kK$ is algebraically closed, see Section~\ref{subsection: fields and algebraic groups}). Define
	\begin{equation}\label{equation: equivariance}
	\mu_{F}(x) = \ksch{e}_{x^{-1}}^*\ \mu_{F}.
	\end{equation}
Set $\ksch{m}(x^{-1}) \ceq \ksch{m}\circ \ksch{e}_{x^{-1}}$. Then $\proj \circ \ksch{e}_{x^{-1}} = \id$ and Equation~\ref{equation: equivariance} defines a family of isomorphisms
	\begin{equation}\label{equation: iii reductive}
	\forall x\in \ksch{P}(\kK),\qquad \mu_{F}(x): \ksch{m}(x^{-1})^*\ F \to F,
	\end{equation}
such that $\mu_{F}(1) = \id_{F}$ and
\begin{equation}\label{equation: condition 3 reductive}
	\xymatrix{
\ksch{m}(x^{-1})^*\ \ksch{m}(y^{-1})^*\ F
\ar[rrr]^{\ksch{m}(x^{-1})^*\mu_F(y)}
\ar[d]
&&&
\ksch{m}(x^{-1})^* F
\ar[d]^{\mu_F(x)}
\\
\ksch{m}((xy)^{-1})^*\ F
\ar[rrr]^{\mu_F(xy)}
&&&
F
	}
\end{equation}
commutes, for all $h,h' \in \ksch{P}(\kK)$, where the left-hand arrow refers to the inverse of the isomorphism in $D^b_c(\ksch{P},\EE)$ determined by the isomorphism of functors $\ksch{m}((xy)^{-1})^*\ \to \ksch{m}(x^{-1})^*\ \ksch{m}(y^{-1})^*$. (We will use Equations~\ref{equation: equivariance}, \ref{equation: iii reductive} and \ref{equation: condition 3 reductive} in Section~\ref{subsection: weakly-equivariant objects}.)

We will also say that a morphism $\phi : F_1\to F_2$ in $\catM\ksch{X}$ is \emph{equivariant} if $F_1$ and $F_2$ are equivariant perverse sheaves and the following diagramme commutes.
\[
\xymatrix{
	\ar[d]_{\mu_{F_1}} \ksch{m}^*F_1 \ar[r]^{\ksch{m}^*\phi} & \ar[d]^{\mu_{F_2}} \ksch{m}^*F_2\\
	\proj^*F_1 \ar[r]^{\proj^*\phi} & \proj^*F
}
\]
Note that this definition makes implicit use of the uniqueness of the isomorphisms $\ksch{m}^* F_1 \to \proj^* F_1$ and $\ksch{m}^* F_2 \to \proj^* F_2$ as above. Since $\id_F$ is equivariant if $F$ is equivariant and since the composition of equivariant morphisms is equivariant, it follows that equivariant perverse sheaves define a category, with morphisms as above, henceforth denoted $\catM_\ksch{P}\ksch{X}$.

\begin{proposition}\label{proposition: torsor}
Let $\ksch{f} : \ksch{X} \to \ksch{Y}$ be a locally trivial principal fibre space with group $\ksch{P}$ and suppose $\ksch{P}$ is connected. Let $F_\ksch{X}$ be a perverse sheaf on $\ksch{X}$. Then $F_\ksch{X}$ is equivariant if and only if $F_\ksch{X} \iso \ksch{f}^*[\dim\ksch{P}] F_\ksch{Y}$ for some perverse sheaf $F_\ksch{Y}$ on $\ksch{Y}$.
\end{proposition}

\begin{proof}
(This result is presented in \cite[1.9.3]{CS}.) By the definition of a locally trivial principal fibre space there is an open covering $\ksch{Y} = \cup_n \ksch{Y}_n$ and isomorphisms $\ksch{t}_n : \ksch{f}^{-1}\ksch{Y}_n \to \ksch{P}\times \ksch{Y}_n$ such that  $\ksch{f}$ is given locally by $\ksch{f}_n = \proj \circ \ksch{t}_n$ (so $\ksch{f}_n : \ksch{f}^{-1}\ksch{Y}_n \to \ksch{Y}_n$) and the action $\ksch{m}_\ksch{X}$ of $\ksch{P}$ on $\ksch{X}$ is given locally by $\ksch{t}_n \circ \ksch{m}_n =  (\ksch{m}_\ksch{P} \times \id) \circ (\id\times\ksch{t}_n)$ (so $\ksch{m}_n : \ksch{P}\times\ksch{f}^{-1}\ksch{Y}_n \to \ksch{f}^{-1}\ksch{Y}_n$). To simplify notation slightly, let $\ksch{X}_n$ denote $\ksch{f}^{-1}\ksch{Y}_n$; also, let $\ksch{j}_n : \ksch{Y}_n \to \ksch{Y}$ and $\ksch{i}_n : \ksch{X}_n \to \ksch{X}$ denote inclusions.

First, suppose $F_\ksch{Y}$ is a perverse sheaf on $\ksch{Y}$. Let $F_\ksch{X} = \ksch{f}^*[\dim\ksch{P}] F_\ksch{Y}$. Since $\ksch{f}$ is smooth with fibres isomorphic to $\ksch{P}$ (so the relative dimension of $\ksch{f}$ is $\dim\ksch{P}$) and since $\ksch{P}$ is geometrically connected (recall that $\kK$ is algebraically closed), it follows from \cite[Prop~4.2.5]{BBD} that $F_\ksch{X}$ is a perverse sheaf. To show that $F_\ksch{X}$ is equivariant we must find an isomorphism $\mu : \ksch{m}^* F_\ksch{X} \to \proj^* F_\ksch{X}$ in $D^b_c(\ksch{P}\times\ksch{X};\EE)$ such that $\ksch{e}^*\mu = \id_{F_\ksch{X}}$, where $\ksch{e} : \ksch{X} \to \ksch{P}\times\ksch{X}$ is the section defined by $x \mapsto (1,x)$. To see this, consider the restriction of $\ksch{m}^* F_\ksch{X}$ to $\ksch{P}\times\ksch{X}_n$:
\begin{eqnarray*}
(\ksch{m}^* F_\ksch{X})\vert_{\ksch{P}\times \ksch{X}_n}
	&=& (\id\times \ksch{i}_n)^* \ksch{m}^* \ksch{f}^* F_\ksch{Y} [\dim\ksch{P}]\\
	&\iso& (\ksch{f}\circ \ksch{m} \circ \id\times \ksch{i}_n)^* F_\ksch{Y} [\dim\ksch{P}]\\
	&=& (\ksch{f} \circ \ksch{i}_n \circ \ksch{m}_n)^* F_\ksch{Y} [\dim\ksch{P}].
\end{eqnarray*}
On the other hand, the restriction of $\proj^* F_\ksch{X}$ to $\ksch{P}\times\ksch{X}_n$ is 
\begin{eqnarray*}
(\proj^* F_\ksch{X})\vert_{\ksch{P}\times \ksch{X}_n}
	&=& (\id\times \ksch{i}_n)^* \proj^* \ksch{f}^* F_\ksch{Y} [\dim\ksch{P}]\\
	&\iso& (\ksch{f}\circ \proj \circ \id\times \ksch{i}_n)^* F_\ksch{Y} [\dim\ksch{P}]\\
	&=& (\ksch{f} \circ \ksch{i}_n \circ \proj)^* F_\ksch{Y} [\dim\ksch{P}].
\end{eqnarray*}
Since $\ksch{f} \circ \ksch{i}_n \circ \ksch{m}_n = \ksch{f} \circ \ksch{i}_n \circ \proj$, we have $(\ksch{m}^* {F_\ksch{X}})\vert_{\ksch{P}\times \ksch{X}_n} \iso (\proj^* F_\ksch{X})\vert_{\ksch{P}\times \ksch{X}_n}$. Since $\cup_n \ksch{P}\times\ksch{X}_n$ is an open cover for $\ksch{P}\times\ksch{X}$, this gives the isomorphism we seek.

Next, suppose $F_\ksch{X} \in \obj\catM\ksch{X}$ is equivariant; thus, $F_\ksch{X} \in \obj\catM_\ksch{P}\ksch{X}$. Let $F_n$ denote the restriction of $F_\ksch{X}$ to $\ksch{X}_n$. Since $\id \times \ksch{i}_n$ satisfies the hypotheses of \cite[Prop~4.2.5]{BBD}, it follows that $F_n$ is a perverse sheaf on $\ksch{X}_n$. Recall that $\ksch{m}: \ksch{P}\times\ksch{X} \to \ksch{X}$ is given locally by $\ksch{m}_n : \ksch{P}\times \ksch{X}_n \to \ksch{X}_n$, as above. Restricting the isomorphism $\mu_{F_\ksch{X}} : \ksch{m}^* F_\ksch{X} \to \proj^* F_\ksch{X}$ to $\ksch{P}\times\ksch{X}_n$ yields the isomorphism $\mu_{F_n} : \ksch{m}_n^* F_n \to \proj^* F_n$. Thus, $F_n$ is an equivariant perverse sheaf. Now, let $\ksch{v}_n : \ksch{Y}_n \to \ksch{X}_n$ be the section of $\ksch{f}_n : \ksch{X}_n \to \ksch{Y}_n$ corresponding to $1\in \ksch{P}(\kK)$ (so $\ksch{v}_n$ is the unique morphism of varieties such that $(\ksch{t}_n \circ \ksch{v}_n) (y) = (1,y)$). Define
\[
F'_n \ceq \ksch{v}_n^* F_n [-\dim\ksch{P}].
\]
Then $F'_n \in \obj D^b_c(\ksch{Y}_n;\EE)$. By standard glueing arguments, the collection of $F'_n \in\obj D^b_c(\ksch{Y}_n;\EE)$, as $\ksch{Y}_n$ ranges over the open cover of $\ksch{Y}$ fixed above, uniquely determines an object $F_\ksch{Y}$ of $D^b_c(\ksch{Y};\EE)$.

It remains to be shown that $F_\ksch{Y}$ is a perverse sheaf. Again, we work locally. For each such $n$,
\begin{eqnarray*}
\ksch{f}_n^*[\dim\ksch{P}] F'_n 
	&=& \ksch{f}_n^*[\dim\ksch{P}] \ksch{v}_n^* F_n [-\dim\ksch{P}]\\
	&\iso& (\ksch{v}_n\circ \ksch{f}_n)^* F_n.
\end{eqnarray*}
Let $\ksch{u}_n : \ksch{X}_n \to \ksch{P}\times\ksch{X}_n$ be the section of $\ksch{m}_n : \ksch{P}\times\ksch{X}_n \to \ksch{X}_n$ corresponding to $1$ (so $\ksch{u}_n$ is the unique morphism of varieties such that $(\id\times\ksch{t}_n)\circ \ksch{u}_n \circ \ksch{t}_n^{-1} (h,y) = (h,1,y)$). The domain of $\ksch{u}_n^*\mu_{F_n}$ is $\ksch{u}_n^*\ksch{m}_n^* F_n \iso (\ksch{m}_n \circ \ksch{u}_n)^* F_n = F_n$, since $\ksch{u}_n$ is a section of $\ksch{m}_n$; the codomain of $\ksch{u}_n^* \mu_{F_n}$ is $\ksch{u}_n^* \proj^* F_n \iso (\proj\circ \ksch{u}_n)^*F_n$. Since $\proj\circ \ksch{u}_n = \ksch{v}_n \circ \ksch{f}_n$, it follows that 
\[
\ksch{u}_n^* \mu_{F_n} : F_n \to \ksch{f}_n^*[\dim\ksch{P}] F'_n
\]
is an isomorphism in $D^b_c(\ksch{X}_n;\EE)$. By \cite[Prop~4.2.5]{BBD} and the fact that $\catM\ksch{X}_n$ is stable in $D^b_c(\ksch{X}_n;\EE)$ under isomorphisms, it follows that $F'_n \in \obj \catM\ksch{Y}_n$. By standard glueing arguments, the collection of isomorphisms $\ksch{u}_n^* \mu_{F_n} \in \mor\catM\ksch{X}_n$, as $\ksch{X}_n$ ranges over the open cover of $\ksch{X}$ fixed above, uniquely determines an isomorphism in $\Hom_{D^b_c(\ksch{X};\EE)}(F_\ksch{X},\ksch{f}_n^*[\dim\ksch{P}] F_\ksch{Y})$.  From \cite[Prop~4.2.5]{BBD}  it follows that $F_\ksch{Y}\in \obj\catM\ksch{Y}$. This completes the proof of Proposition~\ref{proposition: torsor}.
\end{proof}

\begin{proposition}\label{proposition: more torsor}
Let $\ksch{f} : \ksch{X} \to \ksch{Y}$ be a locally trivial principal fibre space with group $\ksch{P}$ and suppose $\ksch{P}$ is connected. Then $\ksch{f}^*[\dim\ksch{P}] : \catM\ksch{Y} \to \catM_\ksch{P}\ksch{X}$ is an equivalence of categories and $\catM_\ksch{P}\ksch{X}$ is a thick subcategory of $\catM\ksch{X}$.
\end{proposition}

\begin{proof}
By \cite[Prop~4.2.5]{BBD} we know that $\ksch{f}^*[\dim\ksch{P}] : \catM\ksch{Y} \to \catM\ksch{X}$ is full and faithful. Let $F_\ksch{Y}$ be a perverse sheaf on $\ksch{Y}$. From the proof of Proposition~\ref{proposition: torsor} we have seen that $F_\ksch{X} \ceq \ksch{f}^*[\dim\ksch{P}] F_\ksch{Y}$ is an equivariant perverse sheaf on $\ksch{X}$ and that $\ksch{f}^*[\dim\ksch{P}] \phi$ is an equivariant morphism in $\catM\ksch{X}$ for each morphism $\phi$ in $\catM\ksch{Y}$. Thus, $\ksch{f}^*[\dim\ksch{P}]$ is a functor from $\catM\ksch{Y}$ to $\catM_\ksch{P}\ksch{X}$. Thus, $\ksch{f}^*[\dim\ksch{P}] : \catM\ksch{Y} \to \catM_\ksch{P}\ksch{X}$ is full and faithful. Proposition~\ref{proposition: torsor} tells us that this functor is essentially surjective. Thus, $\ksch{f}^*[\dim\ksch{P}]$ is an equivalence. The last clause of Proposition~\ref{proposition: more torsor} follows from \cite[4.2.6]{BBD}.
\end{proof}


\subsection{Parabolic induction on reductive quotients}\label{subsection: parabolic induction on reductive quotients}

Let $\ksch{H}$ be a reductive algebraic group over $\kK$ and let $\ksch{P}$ be a parabolic subgroup of $\ksch{H}$ with levi component $\ksch{L}$ and unipotent radical $\ksch{U}$. Let $\ksch{r} : \ksch{P}\to \ksch{L}$ denote the reductive quotient map and let $\ksch{s} : \ksch{P} \to \ksch{H}$ be inclusion. Equip $\ksch{X} \ceq \ksch{H}\times\ksch{P}$ with the $\ksch{P}$-action defined by $p\cdot(g,h) = (pg^{-1},php^{-1})$ and let $\ksch{Y}$ denote the quotient by this action. (This is a variety!)
Consider the diagramme
\begin{equation}\label{equation: parabolic induction on reductive quotients}
\xymatrix{
\ksch{L} & \ar[l]_{\ksch{a}} \ksch{X} \ar[r]^{\ksch{b}} & \ksch{Y} \ar[r]^{\ksch{c}}& \ksch{H} 
}
\end{equation}
where $\ksch{a}(g,h) = \ksch{r}(h)$, $\ksch{b}(g,h) = [g,h]$ and $\ksch{c}[g,h] = ghg^{-1}$. Observe that $\ksch{a}$ is smooth with connected fibres of equal dimension $\dim\ksch{H} + \dim\ksch{U}$, which is therefore the relative dimension of $\ksch{a}$. Observe that $\ksch{b}$ is a locally trivial principal fibre bundle with group $\ksch{P}$, which is connected. Observe also that $\ksch{c}$ is proper. Let $F$ be an equivariant perverse sheaf on $\ksch{L}$ with respect to conjugation $\ksch{m}_\ksch{L}$. It follows from \cite[Prop~4.2.5]{BBD} that $\ksch{a}^*[\dim\ksch{a}] F$ is a perverse sheaf on $\ksch{X}$. Moreover, since $\ksch{a}$ is $\ksch{P}$-equivariant (with respect to the action on $\ksch{X}$ defined above and the action $p\cdot l = \ksch{r}(p)l \ksch{r}(p)^{-1}$ on $\ksch{L}$) and since $F$ is $\ksch{P}$-equivariant with respect to the action just defined on $\ksch{L}$, it follows that $\ksch{a}^*[\dim\ksch{a}] F$ is a $\ksch{P}$-equivariant perverse sheaf on $\ksch{X}$. Let 
\begin{equation}
F_\ksch{X} = \ksch{a}^*[\dim\ksch{a}] F. 
\end{equation}
Since $\ksch{b} : \ksch{X} \to \ksch{Y}$ is a locally trivial principle fibre space with group $\ksch{P}$, and since $\ksch{P}$ is connected, it follows from Proposition~\ref{proposition: torsor} that there is a some perverse sheaf $F_\ksch{Y}$ on $\ksch{Y}$ such that 
\begin{equation}
F_\ksch{X} = \ksch{b}^*[\dim\ksch{P}] F_\ksch{Y}.
\end{equation}
Note that Proposition~\ref{proposition: torsor} tells us exactly how to construct the perverse sheaf $F_\ksch{Y}$. Define 
\begin{equation}
\ind^\ksch{H}_\ksch{P} F \ceq {\ksch{c}}_* F_\ksch{Y}
\end{equation}
(Since $\ksch{c}$ is proper, this is equal to ${\ksch{c}}_! F_\ksch{Y}$.) 
Notice that, \textit{a priori}, $\ind^\ksch{H}_\ksch{P} F$ is an object of $D^b_c(\ksch{H};\EE)$; we do not claim that this is a perverse sheaf.

Next, let $\phi$ be a morphism of $\ksch{P}$-equivariant perverse sheaves on $\ksch{L}$. Then $\ksch{a}^*[\dim\ksch{H}] \phi$ is a morphism of $\ksch{P}$-equivariant perverse sheaves on $\ksch{X}$ (by \cite[Prop~4.2.5]{BBD} and arguments as above). Using the equivalence of categories in Proposition~\ref{proposition: torsor} again, there is a unique (given the choices made  above) morphism $\phi_\ksch{Y}$ of perverse sheaves on $\ksch{Y}$ such that 
\[
\phi_\ksch{X} = \ksch{b}^*[\dim\ksch{P}] \phi_\ksch{Y}.
\] 
Let $\phi_\ksch{X}$ be that morphism. Define 
\[
\ind^\ksch{H}_\ksch{P} \phi \ceq {\ksch{c}}_* \phi_\ksch{Y}.
\] 
Notice that $\ind^\ksch{H}_\ksch{P} \phi$ is an morphism $D^b_c(\ksch{H};\EE)$. Thus, we have defined a \emph{functor}
\[
\ind^\ksch{H}_\ksch{P} : \catM_\ksch{L} \ksch{L} \to D^b_c(\ksch{H};\EE).
\]

In this paper we will use the above construction with $\ksch{H}= \quo{G}_i$, $\ksch{P}=\quo{G}_{i\leq j}$ (see Section~\ref{subsection: restriction between reductive quotients}) and $\ksch{L} = \quo{G}_j$, where $i$ and $j$ are facets of $I(\Ksch{G},\KK)$ and $i \leq j$ in the Bruhat order. In that case we will denote the functor $\ind^{\quo{G}_i}_{\quo{G}_{i\leq j}}$ by $\ind_{i\leq j}$.

\begin{lemma}\label{lemma: conjugation and induction}
Let $g$ be an element of $\Ksch{G}(\KK)$ and let $i,j$ be facets of $I(\Ksch{G},\KK)$ such that $i\leq j$. There is an isomorphism of functors
\[
\ind_{i\leq j}\ \quo{m}(g)_j^*\  \iso \quo{m}(g)_i^*\ \ind_{g i \leq g j}
\]
in the category of equivariant perverse sheaves on $\quo{G}_j$.
\end{lemma}

\begin{proof}
Work locally and use the construction appearing in the proof of Proposition~\ref{proposition: torsor}.
\end{proof}


\subsection{Character Sheaves}\label{subsection: character sheaves}

Let $\ksch{H}$ be a connected reductive algebraic group over $\kK$. We recall that a character sheaf is an irreducible perverse sheaf satisfying any one (and hence all) of the conditions appearing in \cite[Prop.2.9]{CS}. We also remind the reader that an irreducible perverse sheaf on $\ksch{H}$ is \emph{admissible} if it is an irreducible component of $\ind^\ksch{H}_\ksch{P} F$ for some parabolic subgroup $\ksch{P}$ and some cuspidal perverse sheaf $F$ on the levi component $\ksch{L}$ of $\ksch{P}$ (see \cite[7.1.10]{CS})  (Cuspidal perverse sheaves are defined in \cite[7.1.1]{CS}.)  In \cite[Th.23.1]{CS} it is shown that, under some extremely mild conditions on $\ksch{H}$ (which are satisfied if $p \geq 7$, for example), these two classes of perverse sheaves coincide.

Regarding parabolic induction as defined in Section~\ref{subsection: parabolic induction on reductive quotients}, if $F$ is a character sheaf on $\ksch{L}$ then $F$ is equivariant for the conjugation action (see \cite[Prop.2.18]{CS}), in which case $\ind^\ksch{H}_\ksch{P} F$ is defined. Moreover, in \cite[Prop.4.8]{CS} it is shown that if $F$ is a character sheaf on $\ksch{L}$ then $\ind^\ksch{H}_\ksch{P} F$ is a finite direct sum of character sheaves on $\ksch{H}$, and therefore equivariant for the conjugation action on $\ksch{H}$. Thus, if $F$ is a finite direct sum of character sheaves on $\ksch{L}$ then $\ind^\ksch{H}_\ksch{P} F$ is a finite direct sum of character sheaves on $\ksch{H}$. 
From our treatment of parabolic induction as a functor, we see further that if $\phi$ is an isomorphism in the category of equivariant perverse sheaves on $\ksch{L}$ and the domain and codomain of $\phi$ are finite direct sums of character sheaves on $\ksch{L}$ then $\ind^\ksch{H}_\ksch{P} \phi$ is an isomorphism in the category of equivariant perverse sheaves on $\ksch{H}$. We will use this fact in Section~\ref{subsection: parabolic induction}.

In this paper we will use these facts with $\ksch{H}= \quo{G}_i$ and 
$\ksch{P} = \quo{G}_{i\leq j}$ where $i$ and $j$ are facets of $I(\Ksch{G},\KK)$ 
with $i\leq j$.


\subsection{Categories}\label{subsection: categories}

We may now introduce the main categories appearing in this paper. 

\begin{definition}\label{definition: catqD}
Let $\catqD\Ksch{G}$ denote the following category.
\begin{enumerate}
\item[\bf obj:]
An object $A$ of $\catqD\Ksch{G}$ is a family of objects
\[
\left\{ A_i \in  D^b_c(\quo{G}_i;\EE) \tq i\text{\ facet of }
I(\Ksch{G},\KK)\right\},
\]
equipped with a family of isomorphisms
\[
\left\{
A_{i\leq j}\in\Hom_{D^b_c(\quo{G}_j;\EE)}(\res_{i\leq j} A_{i},A_{j})
\tq i \leq j \text{ in } I(\Ksch{G},\KK)\right\},
\]
such that $A_{i\leq i} = \id_{A_i}$ for each facet $i$ of $I(\Ksch{G},\KK)$, and such that the diagramme
\[
	\xymatrix{
\res_{j\leq k} \res_{i\leq j} A_i 
	\ar[rr]^{\res_{j\leq k} A_{i\leq j}} 
	\ar[d]_{\res_{i\leq j\leq k} A_i} 
&&
\res_{j\leq k} A_j 
	\ar[d]^{A_{j\leq k}}
	\\
\res_{i\leq k} A_i 
	\ar[rr]_{A_{i\leq k}} 
&&
A_k
	}
\]
is commutative for each triplet $i,j,k$ of facets of $I(\Ksch{G},\KK)$ such that $i\leq j \leq k$. The isomorphism of sheaves appearing on the left-hand side of this diagramme is determined by the isomorphism of functors appearing in Proposition~\ref{proposition: transitive restriction on quotients}.
\item[\bf mor:]
A morphism $\phi \in \Hom_{\catqD\Ksch{G}}(A,B)$ in the category $\catqD\Ksch{G}$ is a family
\[
\left\{ \phi_i \in \Hom_{ D^b_c(\quo{G}_i;\EE)}(A_i,B_i) \tq i\text{\
facet of } I(\Ksch{G},\KK)\right\},
\]
such that the diagramme
\[
\xymatrix{
\res_{i\leq j} A_i 
	\ar[rr]^{\res_{i\leq j} \phi_i} 
	\ar[d]_{A_{i\leq j}}
&& 
\res_{i\leq j} B_i 
	\ar[d]^{B_{i\leq j}}
\\
A_j 
	\ar[rr]_{\phi_j} 
&& 
B_j
	}
\]
is commutative for each pair $i,j$ of facets of $I(\Ksch{G},\KK)$ such that $i \leq j$. 
\item[\bf com:]
If $u$ and $v$ are morphisms in $\catqD\Ksch{G}$ then the composition $u\circ v$ is defined in $\catqD\Ksch{G}$ by $(u\circ v)_i = u_i\circ v_i$ for each facet $i$ of $I(\Ksch{G},\KK)$. 
\item[\bf id:]
For any object $A$, the identity $\id_A : A\to A$ is defined by $(\id_A)_i = \id_{A_i}$ for each facet $i$ of $I(\Ksch{G},\KK)$.  
\end{enumerate}
Let $\catqC\Ksch{G}$ denote the full subcategory of $\catqD\Ksch{G}$ consisting of objects $A\in \obj\catqD\Ksch{G}$ such that $A_i$ is a perverse sheaf for $\quo{G}_i$ for each facet $i$ of $I(\Ksch{G},\KK)$.
A \emph{coefficient system (for $\Ksch{G}$)} is an object of $\catqC\Ksch{G}$.
\end{definition}


\section{Cuspidal coefficient systems}\label{section: cuspidal coefficient systems}

Let $\KK$, $\RK$ and $\kK$ be as in Section~\ref{section: fundamental notions}. Likewise, let $\Ksch{G}$ be a connected reductive linear algebraic group over $\KK$ satisfying the conditions of Section~\ref{section: fundamental notions} (and in particular, Section~\ref{subsection: stabilizers}).


\subsection{Parabolic restriction}\label{subsection: parabolic restriction}

Let $\Ksch{P} \subseteq \Ksch{G}$ be a parabolic subgroup with reductive quotient $\Ksch{L}$. Recall the notation of Section~\ref{subsection: parabolic restriction on quotients}. 

\begin{proposition}\label{proposition: parabolic restriction}
There is a canonical functor $\res^{\Ksch{G}}_{\Ksch{P}} : \catqD\Ksch{G} \to \catqD\Ksch{L}$ such that for each facet $i$ of $I(\Ksch{L},\KK)$, $(\res^{\Ksch{G}}_{\Ksch{P}} A)_i=\res_{i_\Ksch{G} \leq i_\Ksch{P} } A_{i_\Ksch{G}}$, considered as an object in $D^b_c(\quo{L}_i,\EE)$.
\end{proposition}

\begin{proof}
The functor is defined as follows. Let $i$ and $j$ be facets of $I(\Ksch{L},\KK)$ with $i \leq j$. For any object $A$ in $\catqD\Ksch{G}$, define 
\[
(\res^{\Ksch{G}}_{\Ksch{P}} A)_i\ceq \res_{i_\Ksch{G} \leq i_\Ksch{P} } A_{i_\Ksch{G}},
\] 
considered as an object in $D^b_c(\quo{L}_i,\EE)$ using the identification $\quo{G}_{i_\Ksch{P}} =\quo{L}_i$; also define 
\[
(\res^{\Ksch{G}}_{\Ksch{P}} A)_{i\leq j}\ceq \res_{j_\Ksch{G} \leq j_\Ksch{P} } A_{i_\Ksch{G} \leq j_\Ksch{G} }  \circ \res^{\Ksch{P}}_{i\leq j} A_{i_\Ksch{G}},
\] 
likewise considered as a morphism in $D^b_c(\quo{L}_j,\EE)$ (see Lemma~\ref{lemma: parabolic restriction on quotients} for the definition of the natural transformation $\res^{\Ksch{P}}_{i\leq j}$). For any morphism $\phi$ in $\catqD\Ksch{G}$, define 
\[
(\res^{\Ksch{G}}_{\Ksch{P}} \phi)_i\ceq \res_{i_\Ksch{G} \leq i_\Ksch{P} } \phi_{i_\Ksch{G}},
\] 
considered as a morphism in $D^b_c(\quo{L}_i,\EE)$.

We must verify that $\res^{\Ksch{G}}_{\Ksch{P}} A$ is an object of $\catqD\Ksch{L}$ (\cf Definition~\ref{definition: catqD}(obj)). Using the definition above it follows that
\[
(\res^{\Ksch{G}}_{\Ksch{P}} A)_{i\leq i} = \res_{i_\Ksch{G} \leq i_\Ksch{P} } A_{i_\Ksch{G} \leq i_\Ksch{G} }  \circ \res^{\Ksch{P}}_{i\leq i}  A_{i_\Ksch{G}},
\]
for all facets $i$ of $I(\Ksch{L},\KK)$. By Definition~\ref{definition: catqD}(obj), $A_{i_\Ksch{G} \leq i_\Ksch{G} }= \id_{A_{i_\Ksch{G}}}$. From the definition of $\res^{\Ksch{P}}_{i\leq j}$ appearing in Lemma~\ref{lemma: parabolic restriction on quotients} we see that $\res^{\Ksch{P}}_{i\leq i} = \id_{\res_{i_\Ksch{G} \leq i_\Ksch{P}}}$. Thus, \[(\res^{\Ksch{G}}_{\Ksch{P}} A)_{i\leq i } = \id_{(\res^{\Ksch{G}}_{\Ksch{P}} A)_i}.\]
Having shown that $\res^{\Ksch{G}}_{\Ksch{P}} A$ satisfies the first condition set out in Definition~\ref{definition: catqD}(obj), we now turn to the second part of Definition~\ref{definition: catqD}(obj). Suppose $i$, $j$ and $k$ are facets of $I(\Ksch{L},\KK)$ with $i \leq j \leq k$. To show that $\res^{\Ksch{G}}_{\Ksch{P}} A$ satisfies the second condition appearing in Definition~\ref{definition: catqD}(obj) we must show that the following diagramme commutes.
\begin{equation}\label{equation: lpr.0}
	\xymatrix{
	\ar[d]_{\res_{i\leq j\leq k} (\res^{\Ksch{G}}_{\Ksch{P}} A)_{i\leq j}} 
	\res_{j\leq k} \res_{i\leq j} (\res^{\Ksch{G}}_{\Ksch{P}} A)_{i} 
	\ar[rrr]^{\res_{j\leq k} (\res^{\Ksch{G}}_{\Ksch{P}} A)_{i\leq j} } 
&&& 
	\ar[d]^{(\res^{\Ksch{G}}_{\Ksch{P}} A)_{j\leq k}} 
	\res_{j\leq k} (\res^{\Ksch{G}}_{\Ksch{P}} A)_{j}
\\
\res_{i\leq k}  (\res^{\Ksch{G}}_{\Ksch{P}} A)_{i} 
	\ar[rrr]_{(\res^{\Ksch{G}}_{\Ksch{P}} A)_{i\leq k}}
&&&  
(\res^{\Ksch{G}}_{\Ksch{P}} A)_{k} 
}
\end{equation}
To do this, we begin be recalling (from Section~\ref{subsection: parabolic restriction on quotients}) that $\quo{L}_i = \quo{G}_{i_\Ksch{P}}$ (likewise, $\quo{L}_j = \quo{G}_{j_\Ksch{P}}$ and $\quo{L}_k = \quo{G}_{k_\Ksch{P}}$). Together with the definition of $(\res^{\Ksch{G}}_{\Ksch{P}} A)_i$ given in above, the top left-hand corner of Diagramme~\ref{equation: lpr.0} may be re-written as follows:
\[
\res_{j\leq k} \res_{i\leq j} (\res^{\Ksch{G}}_{\Ksch{P}} A)_{i} 
= 
\res_{j_\Ksch{P}\leq k_\Ksch{P}} \res_{i_\Ksch{P}\leq j_\Ksch{P}} \res_{i_\Ksch{G} \leq i_\Ksch{P}} A_{i_\Ksch{G}},
\]
and likewise for all the corners of Diagramme~\ref{equation: lpr.0}. Now, to show that Diagramme~\ref{equation: lpr.0} commutes, consider the diagramme below, in which the outer square is exactly Diagramme~\ref{equation: lpr.0}. (To save space we write $\res^{i_\Ksch{G}}_{j_\Ksch{G}}$ for $\res_{i_\Ksch{G}\leq j_\Ksch{G}}$, etc...)
\[
	\xymatrix{
	\ar[ddd]
\res^{j_\Ksch{P}}_{k_\Ksch{P}} \res^{i_\Ksch{P}}_{j_\Ksch{P}} \res^{i_\Ksch{G}}_{i_\Ksch{P}} A_{i_\Ksch{G}} 
	\ar[rrr] 
	\ar[dr]_{4.}
&&& 
	\ar[ddd]
\res^{j_\Ksch{P}}_{k_\Ksch{P}} \res^{j_\Ksch{G}}_{j_\Ksch{P}} A_{j_\Ksch{G}} 
	\ar[dl]_{2.} 
\\
&	
\res^{k_\Ksch{G}}_{k_\Ksch{P}} \res^{j_\Ksch{G}}_{k_\Ksch{G}} \res^{i_\Ksch{G}}_{j_\Ksch{G}} A_{i_\Ksch{G}}
	\ar[r]
	\ar[d] 
& 
\res^{k_\Ksch{G}}_{k_\Ksch{P}} \res^{j_\Ksch{G}}_{k_\Ksch{G}} A_{j_\Ksch{G}} 
	\ar[d]
& 
\\
&	
\res^{k_\Ksch{G}}_{k_\Ksch{P}} \res^{i_\Ksch{G}}_{k_\Ksch{G}} A_{i_\Ksch{G}} 
	\ar[r]
& 
\res^{k_\Ksch{G}}_{k_\Ksch{P}} A_{k_\Ksch{G}} 
& 
\\
	\ar[ur]_{3.} 
\res^{i_\Ksch{P}}_{k_\Ksch{P}}  \res^{i_\Ksch{G}}_{i_\Ksch{P}} A_{i_\Ksch{G}} 
	\ar[rrr]
&&&  
\res^{k_\Ksch{G}}_{k_\Ksch{P}} A_{k_\Ksch{G}} 
	\ar[ul]^{1.} 
}
\]
The inner square is the result of applying the functor $\res^{k_\Ksch{G}}_{k_\Ksch{P}}$ to the relevant form of the commuting square appearing in Definition~\ref{definition: catqD}(obj), and therefore commutes, since $A$ is an object of $\catqD\Ksch{G}$. The arrow marked $1.$ is the identity. The arrow marked $2.$ is $\res^{\Ksch{P}}_{j \leq k} A_{j_\Ksch{G}}$, so the right-hand square commutes by virtue of the definition of $(\res^{\Ksch{G}}_{\Ksch{P}} A)_{j\leq k}$; likewise, the arrow marked $3.$ is $\res^{\Ksch{P}}_{i \leq k} A_{i_\Ksch{G}}$, so the bottom square commutes by virtue of the definition of $(\res^{\Ksch{G}}_{\Ksch{P}} A)_{i\leq k}$. The arrow marked $4.$ is $\res^{\Ksch{P}}_{j \leq k} \res_{i_\Ksch{G}\leq j_\Ksch{G}} A_{i_\Ksch{G}} \circ \res_{j\leq k} \res^{\Ksch{P}}_{i \leq j} A_{i_\Ksch{G}}$ and the top and left-hand squares commute by Lemma~\ref{lemma: parabolic restriction on quotients}. This, the outer square in the diagramme above commutes. Therefore, Diagramme~\ref{equation: lpr.0} commutes. This concludes the demonstration that $\res^{\Ksch{G}}_{\Ksch{P}}A$ is an object in $\catqD\Ksch{L}$.

Suppose $\phi : A \to B$ is a morphism in $\catqD\Ksch{G}$ and let $i$ and $j$ be facets of $I(\Ksch{L},\KK)$ such that $i\leq j$. In order to show that $\res^{\Ksch{G}}_{\Ksch{P}}\phi$ is a morphism in $\catqD\Ksch{L}$ we must show that the following diagramme commutes.
\begin{equation}\label{equation: lpr.1}
\xymatrix{
\res_{i\leq j} (\res^{\Ksch{G}}_{\Ksch{P}} A)_i 
	\ar[rrr]^{\res_{i\leq j} (\res^{\Ksch{G}}_{\Ksch{P}} \phi)_i} 
	\ar[d]_{(\res^{\Ksch{G}}_{\Ksch{P}} A)_{i\leq j}}
&&& 
\res_{i\leq j} (\res^{\Ksch{G}}_{\Ksch{P}} B)_i 
	\ar[d]^{(\res^{\Ksch{G}}_{\Ksch{P}} B)_{i\leq j}}
\\
(\res^{\Ksch{G}}_{\Ksch{P}} A)_j 
	\ar[rrr]_{(\res^{\Ksch{G}}_{\Ksch{P}} \phi)_j} 
&&& 
(\res^{\Ksch{G}}_{\Ksch{P}} B)_j
	}
\end{equation}
As above, we begin by translating this from a statement about morphisms of sheaves on the reductive quotients of parahoric subgroups of $\Ksch{L}$ to a statement about morphisms of sheaves on the reductive quotients of parahoric subgroups of $\Ksch{G}$, where we will see that the diagramme commutes. Thus, for example, the top arrow becomes
\[
{\res_{i\leq j} (\res^{\Ksch{G}}_{\Ksch{P}} \phi)_i}
	= \res_{i_\Ksch{P}\leq j_\Ksch{P}} \res_{i_\Ksch{G}\leq i_\Ksch{P}} \phi_{i_\Ksch{G}}.
\]
Now, consider the diagramme below, in which the outer square is exactly Diagramme~\ref{equation: lpr.1}. (To save space we write $\res^{i_\Ksch{G}}_{j_\Ksch{G}}$ for $\res_{i_\Ksch{G}\leq j_\Ksch{G}}$, etc..., as above.)
\[
\xymatrix{
\res^{i_\Ksch{P}}_{j_\Ksch{P}} \res^{i_\Ksch{G}}_{i_\Ksch{P}} A_{i_\Ksch{G}} 
	\ar[rrr]
	\ar[dd]
	\ar[dr]^{3.}
&&& 
\res^{i_\Ksch{P}}_{j_\Ksch{P}} \res^{i_\Ksch{G}}_{i_\Ksch{P}} B_{i_\Ksch{G}} 
	\ar[dd]
	\ar[dl]^{4.}
\\
& 
	\ar[dl]^{2.}
\res^{j_\Ksch{G}}_{j_\Ksch{P}} \res^{i_\Ksch{G}}_{j_\Ksch{G}}  A_{i_\Ksch{G}} 
	\ar[r]^{0.}
&
\res^{j_\Ksch{G}}_{j_\Ksch{P}} \res^{i_\Ksch{G}}_{j_\Ksch{G}} B_{i_\Ksch{G}} 
	\ar[dr]^{1.}
&
\\
\res^{j_\Ksch{G}}_{j_\Ksch{P}} A_{j_\Ksch{G}} 
	\ar[rrr]
&&& 
\res^{j_\Ksch{G}}_{j_\Ksch{P}} B_{j_\Ksch{G}}
	}
\]
The arrow marked $0.$ is $\res_{j_\Ksch{G}\leq j_\Ksch{P}} \res_{i_\Ksch{G}\leq j_\Ksch{G}}  \phi_{i_\Ksch{G}}$, the arrow marked $1.$ is $\res_{j_\Ksch{G}\leq j_\Ksch{P}} B_{i_\Ksch{G}\leq j_\Ksch{G}}$ and the arrow marked $2.$ is $\res_{j_\Ksch{G}\leq j_\Ksch{P}} A_{i_\Ksch{G}\leq j_\Ksch{G}}$; thus, the bottom square is the result of applying the functor 
$\res_{j_\Ksch{G}\leq j_\Ksch{P}}$ to the relevant form of the commuting square appearing in Definition~\ref{definition: catqD}(mor), and therefore commutes since $\phi$ is a morphism in $\catqD\Ksch{G}$. The arrow marked $3.$ is $\res^\Ksch{P}_{i\leq j} A_{i_\Ksch{G}}$ and the arrow marked $3.$ is $\res^\Ksch{P}_{i\leq j} B_{i_\Ksch{G}}$, so the upper square commutes because $\res^\Ksch{P}_{i\leq j}$ is a natural transformation. The left-hand triangle commutes by virtue of the definition of $(\res^{\Ksch{G}}_{\Ksch{P}} A)_{i\leq j}$ and likewise the right-hand triangle commutes by virtue of the definition of $(\res^{\Ksch{G}}_{\Ksch{P}}B)_{i\leq j}$. Therefore, the outer square commutes. This concludes the demonstration that $\res^{\Ksch{G}}_{\Ksch{P}} : \catqD\Ksch{G} \to \catqD\Ksch{L}$ is a functor.
\end{proof}

Let $\Ksch{P} \to \Ksch{G}$ be a parabolic subgroup containing Borel $\Ksch{B}$ and with Levi component (\ie maximal
reductive quotient) $\Ksch{P} \to \Ksch{L}$. Let $\Ksch{Q}\to \Ksch{L}$ be a parabolic subgroup containing $\Ksch{B}\cap \Ksch{L}$ with Levi component $\Ksch{Q} \to \Ksch{M}$. Let $ \Ksch{P}\hookleftarrow \Ksch{R} \to \Ksch{Q}$
be a pull-back of $\Ksch{P} \to \Ksch{L} \from \Ksch{Q}$ in the category of group $\KK$-schemes. Then $\Ksch{R} \to \Ksch{G}$ is a parabolic subgroup with Levi component $\Ksch{R} \to \Ksch{L}$ and $\Ksch{R} = \Ksch{Q}\Ksch{U}$, where $\Ksch{U}$ is the kernel of $\Ksch{P} \to \Ksch{L}$.


\begin{proposition}\label{proposition: transitive restriction}
With notation as above, $\res^{\Ksch{L}}_{\Ksch{Q}}\ \res^{\Ksch{G}}_{\Ksch{P}}  \cong
\res^{\Ksch{G}}_{\Ksch{R}}$.
\end{proposition}

\begin{proof}
Proposition~\ref{proposition: transitive restriction} is a consequence of Proposition~\ref{proposition: parabolic restriction} and Proposition~\ref{proposition: transitive restriction on quotients}.
\end{proof}

\begin{remark}\label{remark: transitive restriction}
It was very important for us to keep track of the isomorphisms appearing in Proposition~\ref{proposition: transitive restriction on quotients} in order to have a good definition of $\catqD\Ksch{G}$ and in order for us to define parabolic restriction above. However, it is not important for us to keep track of the isomorphism in Proposition~\ref{proposition: transitive restriction}, as we will see below.
\end{remark}


\subsection{Weakly-equivariant objects}\label{subsection: weakly-equivariant objects}

In this Section we make extensive use of ideas and notation introduced in Section~\ref{subsection: conjugation}. 

\begin{proposition}\label{proposition: conjugation}
For each $g\in \Ksch{G}(\KK)$ there is a canonical functor from $\catqD\Ksch{G}$ to $\catqD\Ksch{G}$
such that the image of $A\in \obj\catqD$ under this functor is an object $\conj{g}{A}$ given by
$\conj{g}{A}_i \ceq \quo{m}(g^{-1})_i^*\ A_{g^{-1}i}$ for each facet $i$ of $I(\Ksch{G},\KK)$.
We shall denote this functor by $\quo{m}(g)^* : \catqD\Ksch{G} \to \catqD\Ksch{G}$
\end{proposition}

\begin{proof}
Fix $g\in \Ksch{G}(\KK)$. For any $A\in \obj \catqD\Ksch{G}$, define $\conj{g}{A}\in \obj \catqD\Ksch{G}$ as follows: for each pair of facets $i,j$ of $I(\Ksch{G},\KK)$ with $i \leq j$,
\[
\conj{g}{A}_i \ceq \quo{m}(g^{-1})_i^*\ A_{g^{-1}i}
\]
and
\[
{\conj{g}{A}}_{i\leq j} \ceq \quo{m}(g^{-1})_{j}^* A_{g^{-1}i\leq g^{-1}j} \circ \res^{g^{-1}}_{g^{-1}i\leq g^{-1}j} A_{g^{-1}i}.
\]
(Here,  $\res^{g^{-1}}_{g^{-1}i\leq g^{-1}j}$ refers to the natural transformation introduced in Lemma~\ref{lemma: conjugation and restriction}.)  
Likewise, for any $\phi\in \Hom_{\catqD\Ksch{G}}(A,B)$ we define $\conj{g}{\phi}\in \Hom_{\catqD\Ksch{G}}(\conj{g}{A},\conj{g}{B})$ by
\[
\conj{g}{\phi}_i \ceq \quo{m}(g^{-1})_i^*\ \phi_{g^{-1}i} ,
\]
for each facet $i$ of $I(\Ksch{G},\KK)$. 

From the proof of Lemma~\ref{lemma: conjugation and restriction}, we see that $\res_{g^{-1}i\leq g^{-1}i} A_{g^{-1}i} = \id_{\quo{m}(g^{-1})_i^*\ A_{g^{-1}i}}$. From Definition~\ref{definition: catqD}(obj) we see that $A_{g^{-1}i\leq g^{-1}i} = \id_{A_{g^{-1}i}}$. Thus, using Proposition~\ref{proposition: conjugation} we have
\begin{eqnarray*}
{\conj{g}{A}}_{i\leq i} 
	&=& \quo{m}(g^{-1})_{i}^*\ A_{g^{-1}i\leq g^{-1}i} \circ \res^{g^{-1}}_{g^{-1}i\leq g^{-1}i} A_{g^{-1}i}\\
	&=& \quo{m}(g^{-1})_{i}^*\ \id_{A_{g^{-1}i}} \circ \id_{\quo{m}(g^{-1})_i^*\ A_{g^{-1}i}}\\
	&=& \id_{\quo{m}(g^{-1})_{i}^*\ A_{g^{-1}i}} \circ \id_{\quo{m}(g^{-1})_i^*\ A_{g^{-1}i}}\\
	&=& \id_{\conj{g}{A}_i} \circ \id_{\conj{g}{A}_i}\\
	&=& \id_{\conj{g}{A}_i}.
\end{eqnarray*}
Having shown that $\conj{g}{A}$ satisfies the first condition set out in Definition~\ref{definition: catqD}(obj), we now turn to the second part of Definition~\ref{definition: catqD}(obj). Suppose $i$, $j$ and $k$ are facets of $I(\Ksch{G},\KK)$ with $i \leq j \leq k$; suppose also that $g\in \Ksch{G}(\KK)$ as above. To show that $\conj{g}{A}$ satisfies the second condition appearing in Definition~\ref{definition: catqD}(obj) we must show that the following diagramme commutes.
\begin{equation}\label{equation: c.0}
	\xymatrix{
	\ar[d]_{\res_{i\leq j\leq k} \conj{g}{A}_{i\leq j}} 
	\res_{j\leq k} \res_{i\leq j} \conj{g}{A}_{i} 
	\ar[rrr]^{\res_{j\leq k} \conj{g}{A}_{i\leq j} } 
&&& 
	\ar[d]^{\conj{g}{A}_{j\leq k}} 
	\res_{j\leq k} \conj{g}{A}_{j}
\\
\res_{i\leq k} \conj{g}{A}_{i} 
	\ar[rrr]_{\conj{g}{A}_{i\leq k}}
&&&  
\conj{g}{A}_{k} 
}
\end{equation}
Consider the following diagramme, in which the outer square is exactly Diagramme~\ref{equation: c.0}. (To save space we have written $\res^{gi}_{gj}$ for $\res_{gi\leq gj}$, etc...)
\[
	\xymatrix{
	\ar[ddd]
\res^{j}_{k}\res_{i\leq j}\quo{m}(g)_i^* A_{gi}
	\ar[dr]^{4.}
	\ar[rrr]
&&&
	\ar[ddd]
\res^{j}_{k} \quo{m}(g)_j^* A_{gj}
	\ar[dl]^{2.}
\\
& 
\quo{m}(g)_k^* \res^{gj}_{gk} \res^{gi}_{gj} A_{gi}
	\ar[r]
	\ar[d]
&
\quo{m}(g)_k^* \res^{gj}_{gk} A_{gj}
	\ar[d]
&
\\
&
\quo{m}(g)_k^* \res^{gi}_{gk} A_{gi}
	\ar[r]
&
\quo{m}(g)_k^* A_{gk}
&
\\
	\ar[ur]^{3.}
\res^{i}_{k} \quo{m}(g)_i^* A_{gi}
	\ar[rrr]
&&&
	\ar[ul]^{1.}
\quo{m}(g)_k^* A_{gk}
	}
\]
The inner square is the result of applying the functor $\quo{m}(g)_k^*$ to the relevant form of the commuting square appearing in Definition~\ref{definition: catqD}(obj), and therefore commutes. The arrow marked $1.$ is the identity. The arrow marked $2.$ is $\res^{g}_{j \leq k} A_{gj}$, so the right-hand square commutes by virtue of the definition of $\conj{g}{A}_{j\leq k}$; likewise, the arrow marked $3.$ is $\res^{g}_{i \leq k} A_{gi}$, so the bottom square commutes by virtue of the definition of $\conj{g}{A}_{i\leq k}$. The arrow marked $4.$ is $\res^{g}_{j \leq k} \res_{gi\leq gj} A_{gi} \circ \res_{j\leq k} \res^{g}_{i \leq j} A_{gi}$ and the top and left-hand squares commute by Lemma~\ref{lemma: conjugation}. This concludes the demonstration that $\conj{g}{A}$ is an object in $\catqD\Ksch{G}$.

Suppose $\phi : A \to B$ is a morphism in $\catqD\Ksch{G}$ and let $i$ and $j$ be facets of $I(\Ksch{G},\KK)$ such that $i\leq j$; also, fix $g\in \Ksch{G}(\KK)$. In order to show that $\conj{g}{\phi}$ is a morphism in $\catqD\Ksch{G}$ we must show that the following diagramme commutes.
\begin{equation}\label{equation: c.1}
\xymatrix{
\res_{i\leq j} \conj{g}{A}_i 
	\ar[rrr]^{\res_{i\leq j} \conj{g}{\phi}_i} 
	\ar[d]_{\conj{g}{A}_{i\leq j}}
&&& 
\res_{i\leq j} \conj{g}{B}_i 
	\ar[d]^{\conj{g}{B}_{i\leq j}}
\\
\conj{g}{A}_j 
	\ar[rrr]_{\conj{g}{\phi}_j} 
&&& 
\conj{g}{B}_j
	}
\end{equation}
To do this, consider the diagramme below, in which the outer square is exactly Diagramme~\ref{equation: c.1}. (To save space we have written $\res^{gi}_{gj}$ for $\res_{gi\leq gj}$, etc... , as above.)
\[
\xymatrix{
\res^{i}_{j} \quo{m}(g)_i^*  A_{gi}
	\ar[rrr] 
	\ar[dd]
	\ar[dr]^{3.}
&&& 
\res^{i}_{j} \quo{m}(g)_i^*  B_{gi} 
	\ar[dd]
	\ar[dl]^{4.}
\\
& 
\quo{m}(g)_j^* \res^{gi}_{gj}  A_{gi} 
	\ar[r]^{0.}
	\ar[dl]^{2.}
&
\quo{m}(g)_j^* \res^{gi}_{gj}  B_{gi} 
	\ar[dr]^{1.}
&
\\
\quo{m}(g)_j^*  A_{gj}
	\ar[rrr]
&&& 
\quo{m}(g)_j^*  B_{gj}
	}
\]
The arrow marked $0.$ is $\quo{m}(g)_j^* \res_{gi\leq gj} \phi_{gi}$, the arrow marked $1.$ is $\quo{m}(g)_j^* B_{gi\leq gj}$ and the arrow marked $2.$ is $\quo{m}(g)_j^* A_{gi\leq gj}$; thus, the bottom square is the result of applying the functor 
$\quo{m}(g)_j^*$ to the relevant form of the commuting square appearing in Definition~\ref{definition: catqD}(mor), and therefore commutes since $\phi$ is a morphism in $\catqD\Ksch{G}$. The arrow marked $3.$ is $\res^{g}_{i\leq j} A_{gi}$ and the arrow marked $3.$ is $\res^{g}_{i\leq j} B_{gi}$, so the upper square commutes because $\res^{g}_{i\leq j} $ is a natural transformation. The left-hand triangle commutes by virtue of the definition of $\conj{g}{A}_{i\leq j}$ and likewise the right-hand triangle commutes by virtue of the definition of $\conj{g}{B}_{i\leq j}$. Therefore, the outer square commutes. This concludes the demonstration that $\conj{g}{\phi}$ is a morphism in $\catqD\Ksch{G}$.
\end{proof}

\begin{lemma}\label{lemma: conjugation}
Let $A$ be an object in $\catqD\Ksch{G}$ and let $g$, $h$ be elements of $\Ksch{G}(\KK)$. Then $\conj{g}{(\conj{h}{A})} \iso \conj{gh}{A}$ in $\catqD\Ksch{G}$.
\end{lemma}

\begin{proof}
By Proposition~\ref{proposition: conjugation}, for each facet $i$ of $I(\ksch{G},\KK)$,
\[
\conj{g}{(\conj{h}{A})}_i =  \quo{m}(g^{-1})_i^*\ \quo{m}(h^{-1})_{g^{-1}i}^*\ A_{h^{-1}g^{-1}i}
\] 
and 
\[
\conj{gh}{A}_i = \quo{m}((gh)^{-1})_i^* A_{(gh)^{-1}i} = \quo{m}(h^{-1}g^{-1})_i^* A_{h^{-1}g^{-1}i}.
\] 
Let 
\[
\phi_i : \quo{m}(h^{-1}g^{-1})_i^* A_{h^{-1}g^{-1}i} \to \quo{m}(g^{-1})_i^*\ \quo{m}(h^{-1})_{g^{-1}i}^*\ A_{h^{-1}g^{-1}i}
\]
be the canonical isomorphism. To prove the Lemma~\ref{lemma: conjugation} we show that $\phi \ceq (\phi_i)_i$ is a morphism in $\catqD\Ksch{G}$; thus, we show that $\phi$ satisfies the condition of Definition~\ref{definition: catqD}(mor). This follows from Lemma~\ref{lemma: conjugation and restriction}.
\end{proof}

Recall the hypothesis of Section~\ref{subsection: stabilizers}. The fact that $\quo{G}_i$ is a \emph{connected} linear algebraic group over an algebraically closed field allows us to apply the theory of character sheaves from \cite{CS} to $\quo{G}_i$. Suppose $A$ is an object of $\catqD\Ksch{G}$ such that $A_i$ is a finite direct sum of character sheaves for $\quo{G}_i$, for each facet $i$ of $I(\Ksch{G},\KK)$. By \cite[2.18]{CS}, $A_i$ is an equivariant perverse sheaf.  For each $x\in \quo{G}_i(\kK)$, let $\mu_{A_i}(x) \in \mor D^b_c(\quo{G}_i;\EE)$ be the isomorphism defined by Equation~\ref{equation: equivariance} in Section~\ref{subsection: equivariant perverse sheaves}.

\begin{definition}\label{definition: weakly-equivariant}
An object $A \in \obj\catqC\Ksch{G}$ is \emph{weakly-equivariant} if the following conditions are met.
\begin{enumerate}
\item[(a)]
For each facet $i$ of $I(\Ksch{G},\KK)$, the perverse sheaf $A_i$ is equivariant. 
\item[(b)]
There is a family of isomorphisms
\[	\mu_A = \left\{ \mu_A(g) \in \Hom_{\catqD\Ksch{G}}(\conj{g}{A},A) \tq g\in\Ksch{G}(\KK)\right\} \]
such that $\mu_A(1) = \id_{A}$ and the diagramme
\[
	\xymatrix{
	\ar[d] \conj{g\, }{(\conj{h}{A})} \ar[r]^{\conj{g}{\mu_A(h)}} & \conj{g}{A} \ar[d]^{\mu_A(g)}\\
	\conj{gh}{A}  \ar[r]_{\mu_A(gh)} & A\\
	}
\]
commutes, for all $g$ and $h$ in  $\Ksch{G}(\KK)$. The arrow appearing on the left-hand side of this diagramme is the isomorphism of Lemma~\ref{lemma: conjugation}.
\item[(c)]
For each facet $i$ of $I(\Ksch{G},\KK)$ and for each $g\in \Rsch{G}_i(\RK)$, 
\[
\mu_{A_i}(\rho_i(g)) = \mu_A(g)_i.
\] 
\end{enumerate}
A morphism $\phi :  A\to B$ of weakly-equivariant objects of $\catqC\Ksch{G}$ is itself \emph{weakly-equivariant} if the diagramme
\begin{equation}\label{diagramme: weakly-equivariant morphism}
	\xymatrix{
\conj{g}{A} \ar[d]_{\mu_A(g)} \ar[r]^{\conj{g}{\phi}}& \conj{g}{B} \ar[d]^{\mu_B(g)}\\
A \ar[r]_{\phi}& B
	}
\end{equation}
commutes for all $g\in \Ksch{G}(\KK)$. 
\end{definition}

\begin{lemma}\label{lemma: iso}
Suppose $A$ is weakly-equivariant. For each $g\in \Ksch{G}(\KK)$, the morphism $\mu_A(g) : \conj{g}{A} \to A$ is an isomorphism in $\catqD\Ksch{G}$, and 
\[
\mu_A(g)^{-1} = \conj{g}{\mu_A(g^{-1})}.
\]
\end{lemma}

\begin{proof}
Using Definition~\ref{definition: weakly-equivariant}(b), we have $\id_{A} = \mu_A(g) \circ \conj{g}{\mu_A(g^{-1})}$. Using the functorality of conjugation yields
\begin{eqnarray*}
\conj{g^{-1}}{(\conj{g}{\mu_A(g^{-1})}\circ \mu_A(g))}
    &=& \conj{1}{\mu_A(g^{-1})}\circ \conj{g^{-1}}{\mu_A(g)}\\
    &=& \mu_A(g^{-1} g)\\
    &=& \id_{A}.
\end{eqnarray*}
Thus,
\begin{eqnarray*}
\conj{g}{\mu_A(g^{-1})}\circ \mu_A(g) &=& \conj{g}{\id_A} \\
&=& \id_{\conj{g}{A}}.
\end{eqnarray*}
It follows that $\mu_A(g)^{-1} = \conj{g}{\mu_A(g^{-1})}$, as
promised.
\end{proof}


\subsection{Cuspidal coefficient systems}\label{subsection: cuspidal}

Recall the definition of the additive category $\catqC\Ksch{G}$ from Section~\ref{subsection: categories}. In particular, recall that objects of $\catqC\Ksch{G}$ are called coefficient systems for $\Ksch{G}$.

\begin{definition}\label{definition: cuspidal}
A non-zero coefficient system $C$ for $\Ksch{G}$ is \emph{cuspidal} if it satisfies the following conditions:
	\begin{enumerate}
	\item[(a)] $C_i$ is a finite direct sum of character sheaves for $\quo{G}_i$, or $0$, for each facet $i$ of $I(\Ksch{G},\KK)$.
	\item[(b)] $C$ is weakly-equivariant (see Definition~\ref{definition: weakly-equivariant}).
	\item[(c)] If $C = A \oplus B$ in $\catqC\Ksch{G}$ and $A$ and $B$ are weakly-equivariant, then $A=0$ or $B=0$. 
	\item[(d)] $\res^{\Ksch{G}}_{\Ksch{P}} C =0$ for each proper parabolic Levi subgroup $\Ksch{L}\subset \Ksch{G}$. 
	\end{enumerate}
Let $\catqA^{(0)}\Ksch{G}$ denote the set of cuspidal coefficient systems for $\Ksch{G}$.
\end{definition}

In this Section we give a complete description of the isomorphism classes in $\catqA^{(0)}\Ksch{G}$.

\begin{proposition}\label{proposition: cuspidal}
Let $i_0$ be a vertex of $I(\Ksch{G},\KK)$ and let $F$ be a cuspidal character sheaf for $\quo{G}_{i_0}$. There is a cuspidal coefficient system $C$ for $\Ksch{G}$ such that $C_{i_0} = F$ and $C_i=0$ unless $i$ is in the $\Ksch{G}(\KK)$-orbit of $i_0$; moreover, up to a weakly-equivariant isomorphism, $C$ is the unique cuspidal coefficient system for $\Ksch{G}$ with these properties.
\end{proposition}

\begin{proof}
We begin by showing existence of $C\in \obj\catqC\Ksch{G}$ with the properties claimed above.  Denote the $\Ksch{G}(\KK)$-orbit of the vertex $i_0$ in $I(\Ksch{G},\KK)$ by $\orbit(i_0)$. Consider the function $\Ksch{G}(\KK) \to \orbit(i_0)$ given by $g \mapsto g i_0$ and let $i \mapsto w_i$ denote a normalised section of that function; thus, $w_{i_0} = 1_{\Ksch{G}(\KK)}$ and for each $i\in \orbit(i_0)$ the element $w_i$ of $\Ksch{G}(\KK)$ satisfies that $w_i i_0 = i$. For each facet $i$ of $I(\Ksch{G},\KK)$, define
\begin{equation}\label{equation: i cuspidal}
C_i \ceq
\begin{cases}
    \quo{m}(w_i^{-1})_i^*\ F ,& i\in \orbit(i_0),\\
    0,   & \text{otherwise}.
\end{cases}
\end{equation}
If $i$ and $j$ are facets of $I(\Ksch{G},\KK)$ and $i\leq j$ then we define
\begin{equation}\label{equation: ii cuspidal}
C_{i\leq j} = 
\begin{cases}
    \id_{C_i},& i=j,\\
    0,   & \text{otherwise}.
\end{cases}
\end{equation}
We will show that Equations~\ref{equation: i cuspidal} and \ref{equation: ii cuspidal} define an object of $\catqC\Ksch{G}$. If $i$ and $j$ are not facets in $\orbit(i_0)$, or if $i=j$, then the diagramme in Definition~\ref{definition: catqD}(obj) is commutative for trivial reasons. Thus, we suppose now that $i$ or $j$ is contained in $\orbit(i_0)$ and $i<j$. Since $i_0$ is a vertex, it follows that $i \in \orbit(i_0)$ and $j \not\in\orbit(i_0)$; therefore $C_j =0$. Since $C_i$ is a cuspidal character sheaf (or $0$) for each such facet $i$, and since $\quo{G}_j$ is a proper Levi subgroup of $\quo{G}_i$, it follows that $\res_{i\leq j} C_i =0$. Thus, in all cases, $C_{i\leq i} = \id_{C_i}$ and $\res_{i\leq j\leq k} C_i \circ C_{i\leq k} = C_{j\leq k}\circ\res_{j\leq k} C_{i\leq j}$ for each triplet $i,j,k$ of facets of $I(\Ksch{G},\KK)$ such that $i\leq j \leq k$.  It is clear from these definitions that $C$ is an object of $\catqD\Ksch{G}$ (see Definition~\ref{definition: catqD}) and that $C_{i_0} = F$. Since $C_i$ is a perverse sheaf for every facet $i$, it follows immediately that $C$ is an object of $\catqC\Ksch{G}$ (see Definition~\ref{definition: catqD}).  We must now show that $C$ is cuspidal (see Definition~\ref{definition: cuspidal}). 

	It is clear that $C$ satisfies the condition appearing in Definition~\ref{definition: cuspidal}(a). In order to demonstrate Definition~\ref{definition: cuspidal}(b) we define a family $\mu_C$ of isomorphisms in $\catqD\Ksch{G}$ satisfying the conditions of Definition~\ref{definition: weakly-equivariant}. First, recall that a cuspidal character sheaf is strongly cuspidal (\cf\cite[7.1.6]{CS}). Using \cite[7.1.1]{CS} and \cite[7.1.5]{CS}, observe that the strongly cuspidal perverse sheaf $F$ is
$\quo{G}_i$-equivariant. For each $x\in \quo{G}_{i_0}(\kK)$, let
\[
\mu_{F}(x) : \quo{m}_{i_0}(x^{-1})^*F \to F
\]
be the isomorphism as in Equation~\ref{equation: iii reductive}. For each $g\in\Ksch{G}(\KK)$ and for each facet $i$ of
$I(\Ksch{G},\KK)$ in the $\Ksch{G}(\KK)$-orbit of $i_0$, define
	\begin{equation}\label{equation: k1}
	k_{i,g} \ceq w_i^{-1} g w_{ig}.
	\end{equation}
Then $k_{i,g}$ is an element of $\Rsch{G}_{i_0}(\RK)$, as we now show. By definition, $i = w_i i_0$; thus, $g^{-1} i = g^{-1}(w_i i_0) = (g^{-1} w_i) i_0$. On the other hand, $g^{-1} i$ is a facet of $\orbit(i_0)$ implies $g^{-1} i = w_{g^{-1} i} i_0$. Comparing these last two equations it follows that $w_i^{-1}g w_{ig} \in \Ksch{G}(\KK)_{i_0}$. (Note that here we use the assumption on $\Ksch{G}$ described in Section~\ref{subsection: stabilizers}.) Now, set $\bar{k}_{i,g} = \rho_{i_0}(k_{i,g})$ (\cf Section~\ref{subsection: integral models}). 

In order to define $\mu_C(g)_i : \conj{g}{C}_i \to C_i$ we first suppose $i\subset \orbit(i_0)$. Then
\begin{eqnarray*}
\conj{g}{C}_i 
	&=& \quo{m}(g^{-1})^*\ C_{ig}\\
	&=& \quo{m}(g^{-1})^*\ \conj{w_{ig}}{F}\\
	&=& \quo{m}(g^{-1})^*\ \quo{m}(w_{ig}^{-1})^*\ F.
\end{eqnarray*}
Now, let 
\begin{equation}
\quo{m}(g^{-1})^*\ \quo{m}(w_{ig}^{-1})^*\ F \iso \quo{m}(w_{ig}^{-1} g^{-1})_{i_0}^*\ F
\end{equation}
be the canonical isomorphism and note also that
\begin{eqnarray*}
\quo{m}(w_{ig}^{-1} g^{-1})_{i_0}^*\ F
	&=& \quo{m}(w_{ig}^{-1} g^{-1} w_i w_{i}^{-1})_{i_0}^*\ F\\
	&=& \quo{m}(k_{i,g}^{-1} w_{i}^{-1})_{i_0}^*\ F.
\end{eqnarray*}
Let
\begin{equation}
\quo{m}(k_{i,g}^{-1} w_i^{-1})_{i_0}^* \ F \iso \quo{m}(w_i^{-1})_{i_0}^* \quo{m}(k_{i,g}^{-1})_{i_0}^*\ F
\end{equation}
be the canonical isomorphism and consider the isomorphism
\begin{equation}
\quo{m}(w_i^{-1})_{i_0}^* \mu_F(\bar{k}_{i,g}) : \quo{m}(w_i^{-1})_{i_0}^* \quo{m}(k_{i,g}^{-1})_{i_0}^*\ F \to \quo{m}(w_i^{-1})_{i_0}^* F.
\end{equation}
Since $\conj{w_i}{F} = C_i$, we let $\mu_C(g)_i : \conj{g}{C}_i \to C_i$ be the composition of the isomorphisms above, when $i\subset \orbit(i_0)$. Otherwise, $\mu_C(g)_i  =0$.

It is clear that the family of morphisms $\mu_C(g)_i\in\mor D^b_c(\quo{G}_i;\EE)$ defined above, as $i$ ranges over all facets of $I(\Ksch{G},\KK)$, defines a morphism $\mu_C(g)$ in $\catqD\Ksch{G}$ (\cf Definition~\ref{definition: catqD}(mor)), since
	\begin{equation}
	\mu_C(g)_j \circ \conj{g}{C_{i\leq j}} = C_{i\leq j} \circ \res_{i\leq j} \mu_C(g)_i,
	\end{equation}
for $i\leq j$ in $I(\Ksch{G},\KK)$. We now show that the family of isomorphisms $\mu_C(g)\in \mor\catqD\Ksch{G}$, as $g$ ranges over $\Ksch{G}(\KK)$, denoted $\mu_C$, satisfies the conditions of Definition~\ref{definition: weakly-equivariant}.  If $i$ is not contained in the $\Ksch{G}(\KK)$-orbit of $i_0$ then these conditions are trivial. We now suppose, therefore, that $i$ is a facet of $\orbit(i_0)$, whence $\mu_C(g)_i = \conj{w_i}{\mu_{F}(\bar{k}_{i,g}^{-1})}$. If $g=1$, we have
$k_{i,g} = k_{i,1} = w_i^{-1} w_i = 1$ and $\mu_{F}(1) = \id_{F}$, so 
\begin{eqnarray*}
\mu_C(1)_i 
	&=& \conj{w_i}{\mu_F(1)}_i\\
	&=& \conj{w_i}{\id_{F}}_i\\
	&=&  \id_{\conj{w_i}{F}_i}\\
	&=& \id_{C_i}.
\end{eqnarray*}
Applying the functor $\quo{m}(w_i^{-1})_{i_0}^*$ to Diagramme~\ref{equation: condition 3 reductive} with $x=\bar{k}_{i,g}$ and $y= \bar{k}_{ig,h}$ (see Equation~\ref{equation: k1}) yields the following commutative diagramme.
\begin{equation}\label{equation: ca.0}
	\xymatrix{
\quo{m}(w_i^{-1})_{i_0}^*\ \ksch{m}(\bar{k}_{i,g}^{-1})^*\ \quo{m}(\bar{k}_{ig,h}^{-1})^*\ F
\ar[rrrr]^{\quo{m}(w_i^{-1})_{i_0}^*\ \ksch{m}(\bar{k}_{i,g}^{-1})^*\mu_F(\bar{k}_{ig,h})}
\ar[d]
&&&&
\quo{m}(w_i^{-1})_{i_0}^*\ \ksch{m}(\bar{k}_{i,g}^{-1})^* F
\ar[d]_{\quo{m}(w_i^{-1})_{i_0}^*\ \mu_F(\bar{k}_{i,g})}
\\
\quo{m}(w_i^{-1})_{i_0}^*\ \ksch{m}((\bar{k}_{i,g}\bar{k}_{ig,h})^{-1})^*\ F
\ar[rrrr]^{\quo{m}(w_i^{-1})_{i_0}^*\ \mu_F(\bar{k}_{i,g}\bar{k}_{ig,h})}
&&&&
\quo{m}(w_i^{-1})_{i_0}^*\ F
	}
\end{equation}
The top arrow is $\conj{w_i\bar{k}_{i,g}}{\mu_{F}(\bar{k}_{ig,h})}$ while the right-hand side arrow is $\conj{w_i}{\mu_{F}(\bar{k}_{i,g})}$; thus, the clockwise path is
\begin{eqnarray*}
 \conj{w_i}{\mu_{F}(\bar{k}_{i,g}) \circ \conj{w_i\bar{k}_{i,g}}{\mu_{F}(\bar{k}_{ig,h})}}
    &=& \conj{w_i}{\mu_{F}(\bar{k}_{i,g})} \circ \conj{g w_{ig}}{\mu_{F}(\bar{k}_{ig,h})}\\
    &=& \conj{w_i}{\mu_{F}(\bar{k}_{i,g})} \circ \quo{m}(g^{-1})^*\ \conj{w_{ig}}{\mu_{F}(\bar{k}_{ig,h})}\\
    &=& \mu_C(g)_i \circ \quo{m}(g^{-1})^*\ \mu_C(h)_{ig}\\
    &=& \mu_C(g)_i \circ \conj{g}{\mu_C(h)}_i.
\end{eqnarray*}
On the other hand, the left-hand arrow is the canonical isomorphism, while the bottom arrow is $\conj{w_i}{\mu_{F}(\bar{k}_{i,g} \bar{k}_{ig,h})}$, and
\begin{eqnarray*}
\conj{w_i}{\mu_{F}(\bar{k}_{i,g} \bar{k}_{ig,h})}
    &=& \conj{w_i}{\mu_{F}(\rho_{i_0}(k_{i,g} k_{ig,h}))}\\
    &=& \conj{w_i}{\mu_{F}(\rho_{i_0}(w_i^{-1} g w_{ig} w_{ig}^{-1} h w_{igh}))}\\
    &=& \conj{w_i}{\mu_{F}(\rho_{i_0}(w_i^{-1} gh w_{igh}))}\\
    &=& \conj{w_i}{\mu_{F}(\bar{k}_{i,gh})} \\
	&=& \mu_C(gh)_i.
\end{eqnarray*}
Thus, Diagramme~\ref{equation: ca.0} gives us the condition appearing in Definition~\ref{definition: cuspidal}(b).

	To show that $C$ satisfies the condition appearing in Definition~\ref{definition: cuspidal}(c), suppose $C = A \oplus B$ in $\catqC\Ksch{G}$ and that $A$ and $B$ are weakly-equivariant. Then, for each facet $i$, $C_i = A_i \oplus B_i$ in the category of perverse sheaves for $\quo{G}_i$. Since $C_i=0$ unless $i$ is in the $\Ksch{G}(\KK)$-orbit of $i_0$, we have $A_i=0$ and $B_i=0$ unless $i$ is in the $\Ksch{G}(\KK)$-orbit of $i_0$. Since $F$ is a character sheaf, it is irreducible, so $C_{i_0} = A_{i_0} \oplus B_{i_0}$ implies $A_{i_0} =0$ or $B_{i_0}$ (recall that $C_{i_0} = F$). Without loss of generality, suppose $B_{i_0} =0$. Since $B$ is weakly-equivariant, this implies $B_i =0$ for each facet $i$ in the $\Ksch{G}(\KK)$-orbit of $i_0$. Since $B_i=0$ when $i$ is not in the $\Ksch{G}(\KK)$-orbit of $i_0$ also, it follows that $B=0$. Thus, $C$ satisfies Definition~\ref{definition: cuspidal}(c). 

	 We now consider the condition of Definition~\ref{definition: cuspidal}(d). Let $\Ksch{P}$ be a proper parabolic subgroup $\Ksch{P}$ of $\Ksch{G}$ and let $\Ksch{L}$ be the reductive quotient of $\Ksch{P}$. Let $i$ be a facet of $I(\Ksch{L},\KK)$. From Proposition~\ref{proposition: parabolic restriction}, we see that $(\res^{\Ksch{G}}_{\Ksch{P}} C)_i$ may be indentified with $\res_{i_\Ksch{G}\leq i_\Ksch{P}} A_{i_\Ksch{G}}$. Since $\Ksch{P}$ is proper, $i_\Ksch{P}$ is strictly greater than $i_\Ksch{G}$; thus, $\quo{G}_{i_\Ksch{P}}$ is a proper levi subgroup of $\quo{G}_{i_\Ksch{G}}$. Since $F$ is a cuspidal character sheaf (and therefore strongly cuspidal) and $A$ is weakly-equivariant, it follows that $A_{i_\Ksch{G}}$ is either strongly cuspidal or $0$; in either case, $i_\Ksch{G} \lneq i_\Ksch{P}$ implies $\res_{i_\Ksch{G}\leq i_\Ksch{P}} A_{i_\Ksch{G}} =0$. Thus, $(\res^{\Ksch{G}}_{\Ksch{P}} C)_i =0$. Since $i$ was an arbitrary facet of $I(\Ksch{L},\KK)$, it follows that $\res^{\Ksch{G}}_{\Ksch{P}} C =0$. Thus, $C$ satisfies Definition~\ref{definition: cuspidal}(d).

	We now show uniqueness. Let $A$ and $B$ be cuspidal objects of $\catqC\Ksch{G}$ such that $A_{i_0} = F = B_{i_0}$ and $A_i=0 =B_i$ unless $i$ is in the $\Ksch{G}(\KK)$-orbit of $i_0$. Note that $\mu_A(k)_{i_0} = \mu_{F}(\rho_{i_0}(k)) = \mu_B(k)_{i_0}$ for $k\in \Rsch{G}_{i_0}(\RK)$ by Definition~\ref{definition: cuspidal}(b). Now, $A_{i_0} = B_{i_0}$ implies $A_{w_i^{-1} i} = B_{w_i^{-1}i}$, hence $\quo{m}(w_i^{-1})_i^*\ A_{w_i^{-1}i} = \quo{m}(w_i^{-1})_i^*\ B_{w_i^{-1}i}$, so $\conj{w_i}{A}_i = \conj{w_i}{B}_i$ for all facets $i$ in the $\Ksch{G}(\KK)$-orbit of $i_0$. Define $\phi : A\to B$ by
\begin{equation}\label{equation: u1}
\phi_i = \begin{cases}
    \mu_B(w_i)_i \circ \conj{w_i}{\mu_A(w_i^{-1})}_i & i \subset \orbit(i_0)\\
    \id_0 & \text{otherwise}.
    \end{cases}
\end{equation}
We will now show that  $\phi$ is an isomorphism in $\catqC\Ksch{G}$. Suppose $i$ is a facet in the $\Ksch{G}(\KK)$-orbit of $i_0$. First, using Equation~\ref{equation: k1} and Definition~\ref{definition: weakly-equivariant}(b) (twice), we have
\begin{eqnarray*}
\mu_A(g)
 &=& \mu_A(w_i k_{i,g} w_{ig}^{-1})\\
 &=& \mu_A(w_i) \circ \conj{w_i}{\mu_A(k_{i,g} w_{ig}^{-1})}\\
 &=& \mu_A(w_i) \circ \conj{w_i}{\left( \mu_A(k_{i,g}\circ \conj{k_{i,g}}{\mu_A(w_{ig}^{-1})}\right)}\\
 &=& \mu_A(w_i) \circ \conj{w_i}{\mu_A(k_{i,g})} \circ \conj{w_ik_{i,g}}{\mu_A(w_{ig}^{-1})}\\
 &=& \mu_A(w_i) \circ \conj{w_i}{\mu_A(k_{i,g})} \circ \conj{g w_{ig}}{\mu_A(w_{ig}^{-1})}.
\end{eqnarray*}
Thus,
\begin{eqnarray*}
\mu_A(g)_i
    &=& \mu_A(w_i)_i \circ \conj{w_i}{\mu_A(k_{i,g})}_i \circ \conj{g w_{ig}}{\mu_A(w_{ig}^{-1})}_i\\
    &=& \mu_A(w_i)_i \circ \quo{m}(w_i^{-1})^*\ \mu_A(k_{i,g})_{w_i^{-1}i} \circ \quo{m}(g^{-1})^*\  \conj{w_{ig}}{\mu_A(w_{ig}^{-1})}_{ig}\\
    &=& \mu_A(w_i)_i \circ \quo{m}(w_i^{-1})^*\ \mu_A(k_{i,g})_{i_0} \circ \quo{m}(g^{-1})^*\  \conj{w_{ig}}{\mu_A(w_{ig}^{-1})}_{ig}\\
    &=& \mu_A(w_i)_i \circ \quo{m}(w_i^{-1})^*\ \mu_{F}(k_{i,g})\circ \quo{m}(g^{-1})^*\
    \conj{w_{ig}}{\mu_A(w_{ig}^{-1})}_{ig}.
\end{eqnarray*}
Likewise,
\begin{eqnarray*}
\mu_B(g)_i
    &=& \mu_B(w_i)_i \circ \quo{m}(w_i^{-1})^*\ \mu_{F}(k_{i,g})\circ \quo{m}(g^{-1})^*\
    \conj{w_{ig}}{\mu_B(w_{ig}^{-1})}_{ig}.
\end{eqnarray*}
Using Lemma~\ref{lemma: iso} and the definition of $\phi^i$ we have
\begin{eqnarray*}
&&\hskip-20pt
 \phi_i \circ \mu_A(g)_i\\
 &=& \mu_B(w_i)_i \circ \conj{w_i}{\mu_A(w_i^{-1})}_i \circ \mu_A(w_i)_i\\
 &&  \circ\ \quo{m}(w_i)^* \mu_{F}(\bar{k}_{i,g}) \circ \quo{m}(g^{-1})^*\ \conj{w_{ig}^{-1}}{\mu_A(w_{ig}^{-1})}_{ig}\\
 &=& \mu_B(w_i)_i \circ \quo{m}(w_i)^*\ \mu_{F}(\bar{k}_{i,g}) \circ \quo{m}(g^{-1})^*\
 \conj{w_{ig}^{-1}}{\mu_A(w_{ig}^{-1})}_{ig}.
\end{eqnarray*}
On the other hand, using the work above and the definition of $\conj{g}{\phi}_i$ we have
\begin{eqnarray*}
&&\hskip-20pt
 \mu_B(g)_i \circ \conj{g}{\phi}_i\\
    &=& \mu_B(w_i)_i \circ \quo{m}(w_i^{-1})^*\
    \mu_{F}(k_{i,g})\circ \quo{m}(g^{-1})^*\
    \conj{w_{ig}}{\mu_B(w_{ig}^{-1})}_{ig} \\
    && \circ\ \quo{m}(g^{-1})^*\left( \mu_B(w_{ig})_{ig}
    \circ \conj{w_{ig}}{\mu_A(w_{ig}^{-1})}_{ig}\right)\\
    &=& \mu_B(w_i)_i \circ \quo{m}(w_i^{-1})^*\
    \mu_{F}(k_{i,g})\\
    && \circ\ \quo{m}(g^{-1})^* \left( \conj{w_{ig}}{\mu_B(w_{ig}^{-1})}_{ig}
    \circ \mu_B(w_{ig})_{ig} \circ
    \conj{w_{ig}}{\mu_A(w_{ig}^{-1})}_{ig}\right)\\
    &=& \mu_B(w_i)_i \circ \quo{m}(w_i^{-1})^*\
    \mu_{F}(k_{i,g}) \circ \quo{m}(g^{-1})^*
    \conj{w_{ig}}{\mu_A(w_{ig}^{-1})}_{ig}.
\end{eqnarray*}
Thus,
\begin{equation}
\phi_i \circ \mu_A(g)_i = \mu_B(g)_i \circ \conj{g}{\phi}_i,
\end{equation}
for all facets $i$ of $I(\Ksch{G},\KK)$. This concludes the proof that $\phi : A\to B$ is a weakly-equivariant morphism in $\catqC\Ksch{G}$. Since $\phi$ is clearly an isomorphism by Lemma~\ref{lemma: iso}, this concludes the proof of Proposition~\ref{proposition: cuspidal}.	
\end{proof}

\begin{definition}\label{definition: compact induction}
For any vertex $i$ of $I(\Ksch{G},\KK)$ and any cuspidal character sheaf $F$ on $\quo{G}_i$, let $\cind_i F$
denote the cuspidal object of $\catqC\Ksch{G}$ given by Proposition~\ref{proposition: cuspidal}. 
\end{definition}

\begin{theorem}\label{theorem: cuspidal}
If $C$ is a cuspidal coefficient system for $\Ksch{G}$ then there is a vertex $i_0$ of $I(\Ksch{G},\KK)$ and a cuspidal character sheaf $F$ such that $C_{i_0} = F$ and $C_i=0$ unless $i$ is in the $\Ksch{G}(\KK)$-orbit of $i_0$. Moreover, this vertex $i_0$ is unique up to $\Ksch{G}(\KK)$-conjugation and $F$ is unique up to isomorphism in $D^b_c(\quo{G}_{i_0};\EE)$.
\end{theorem}

\begin{proof}
If $C_i =0$ for all vertices $i$ of $I(\Ksch{G},\KK)$, then $C=0$ because each $C_{i\leq j} : \res_{i\leq j} C_i \to C_j$ is an isomorphism for each $i \leq j$ (see Definition~\ref{definition: catqD}). Since cuspidal coefficient systems are non-zero (see Definition~\ref{definition: cuspidal}), that is not the case. Thus, there is some vertex $i$ of $I(\Ksch{G},\KK)$ such that $C_i \neq 0$. Using Definition~\ref{definition: cuspidal}(a) we write $C_i = \oplus_m C_{i,m}$ (finite direct sum) where each $C_{i,m}$ is a character sheaf. Let $\ksch{M}$ be a proper Levi subgroup of $\quo{G}_i$. Using \cite{Lan1} we can identify the star of $i$ in $I(\Ksch{G},\KK)$ with the building for $\quo{G}_i$, so there is some facet $j$ in the star of $i$ such that $\ksch{M} = \quo{G}_j$. Since $\ksch{M}\subset \quo{G}_i$ is proper, $i \lneq j$. Now, there is a parabolic subgroup $\Ksch{P}$ with levi component $\Ksch{L}$ and a facet $k$ of $I(\Ksch{L},\KK)$ such that $k_\Ksch{G} = i$ and $k_\Ksch{P} = j$. Note also that $k_\Ksch{G} \lneq k_\Ksch{P}$. By Definition~\ref{definition: cuspidal}(d), $\res^{\Ksch{G}}_{\Ksch{P}} C =0$. Thus, $(\res^{\Ksch{G}}_{\Ksch{P}} C)_k=0$ so $\res_{k_\Ksch{G}\leq k_\Ksch{P}} C_{k_\Ksch{G}} =0$ so $\res_{i\leq j} C_i =0$ so $\res^{\quo{G}_i}_{\ksch{M}} C_i =0$. Thus, $\res^{\quo{G}_i}_{\ksch{M}} C_{i,m} =0$, for each $m$ above. Since this argument applies to any proper levi subgroup $\ksch{M}$ of $\quo{G}_i$, and since character sheaves are equivariant, it follows that $C_{i,m}$ is a cuspidal character sheaf. Thus we have shown that if $i$ is a vertex then $C_i$ is $0$ or is a finite direct sum of cuspidal character sheaves for $\quo{G}_i$. Note that it follows immediately that $C_i=0$ unless $i$ is a vertex.

	Now, let $\{i_0, i_1, \ldots, i_d\}$ be the vertices of a fundamental $\Ksch{G}(\KK)$-domain in $I(\Ksch{G},\KK)$; thus, the convex hull of $\{i_0, i_1, \ldots, i_d\}$ is a chamber in $I(\Ksch{G},\KK)$. For each such vertex $i_n$, we write $C_{i_n} = \oplus_{m} C_{i_n,m}$, where $C_{i_n,m}$ is a cuspidal character sheaf for $\quo{G}_{i_n}$.  Consider 
\[
A \ceq \mathop{\oplus}\limits_{n,m} \cind_{i_n} C_{i_n,m},
\]
where $0\leq n \leq d$ and $0 \leq m \leq d_n$ runs over an index set corresponding to the irreducible summands of $C_{i_n}$ in the category of perverse sheaves on $\quo{G}_{i_n}$. Notice that $A_{i_n} = \oplus_m C_{i_n,m} = C_{i_n} = C_{i_n}$ for each vertex $i_n$ above. Thus, $A_{i_n} = C_{i_n}$ for each vertex $i_n$ above and $A_i = C_i = 0$ unless $i$ is a vertex. Note also that $A$ and $C$ are both weakly-equivariant.

We now show that $A\iso C$. Observe that $A$ and $C$ are weakly-equivariant and $A_{i_n} = C_{i_n}$ for every vertex $i_n \in \{ i_0, \ldots, i_n\}$. For each such vertex $i_n$, let $i \mapsto w(n)_i$ denote a normalized section of $g \mapsto g i_n$, as in the proof of Proposition~\ref{proposition: cuspidal}; thus, $i = w(n)_i i_n$ for each vertex $i$ in the $\Ksch{G}(\KK)$-orbit of $i_n$. Since the set of vertices in $I(\Ksch{G},\KK)$ is partitioned into $\Ksch{G}(\KK)$-orbits, for each vertex $i$ there is a unique vertex $i_n$ such that $i = w(n)_i i_n$. Let $i$ be any vertex and define $\phi_i : A_i \to C_i$ by 
\[
\phi_i 
= 
\mu_C(w(n)_i)_i \circ \mu_A(w(n)_i)_i^{-1}.
\]
This composition is defined since the codomain of $\mu_A(w(n)_i)_i^{-1}$ is the domain of $\mu_A(w(n)_i)_i$ which is $\conj{w(n)_i}{A}_i$, and by Proposition~\ref{proposition: conjugation}, 
\begin{eqnarray*}
\conj{w(n)_i}{A}_i 
	&=& \quo{m}(w(n)_i^{-1})_i^* A_{w(n)_i^{-1} i}\\
 	&=& \quo{m}(w(n)_i^{-1})_i^* A_{i_n}\\
 	&=& \quo{m}(w(n)_i^{-1})_i^* B_{i_n}\\
	&=& \quo{m}(w(n)_i^{-1})_i^* B_{w(n)_i^{-1} i}\\
	&=& \conj{w(n)_i}{B}_i,
\end{eqnarray*} 
which is the domain of $\mu_C(w(n)_i)_i$. Observe also that the domain of $\phi_i$ is indeed $A_i$ since the codomain of $\mu_A(w(n)_i)_i$ is $A_i$; likewise, the codomain of $\phi_i$ is indeed $B_i$ since the codomain of $\mu_C(w(n)_i)_i$ is $C_i$. If $i$ is not a vertex, define $\phi_i = 0$ (Recall that $A_i = C_i =0$ unless $i$ is a vertex.) To see that this defines a morphism in $\catqD\Ksch{G}$ is it necessary to see that the diagramme appearing in Definition~\ref{definition: catqD}(mor) commutes for all $i \leq j$. But the case $i=j$ is trivial since $A_{i \leq j} = C_{i \leq j}  = 1$ if $i = j$; the case $i \lneq j$ is equally trivial, since $A_{i\leq j} = C_{i\leq j} = 0$ if $i \lneq j$; Thus, $\phi$ is a morphism in $\catqD\Ksch{G}$. In fact, from the definition of each $\phi_i$ is is also clear that $\phi$ is an isomorphism in $\catqD\Ksch{G}$. Since $\catqC\Ksch{G}$ is a full subcategory of $\catqD\Ksch{G}$ and since the domain and codomain of $\phi$ are objects of $\catqC\Ksch{G}$ it follows that $\phi$ is an isomorphism in $\catqC\Ksch{G}$. (In fact, $\phi$ is weakly-equivariant.)

Since $A\iso C$ and $C$ is cuspidal, it follows that $A$ is also cuspidal. Since each $\cind_{i_n} C_{i_n,m}$ is weakly-equivariant by Proposition~\ref{proposition: cuspidal}, it follows from Definition~\ref{definition: cuspidal}(c) that $\cind_{i_n} C_{i_n,m}=0$ for all but one vertex $i_n$ and for all but one index $m$. Let $i_0$ be that vertex, set $m=0$ and let $F = C_{i_0,0}$. Then $A = \cind_{i_0} F$. Now, $A\iso C$ implies $C \iso  \cind_{i_0} F$, which completes the proof of Theorem~\ref{theorem: cuspidal}.
\end{proof}

\begin{corollary}\label{corollary: cuspidal}
Let $i_0, i_1, \ldots , i_d$ be a set of representatives for the $\Ksch{G}(\KK)$-orbits of vertices in $I(\Ksch{G},\KK)$. The isomorphism classes in $\catqA^{(0)}\Ksch{G}$ (see Definition~\ref{definition: compact induction}) are parameterized by
\[
\left\{
(i_n,F) \tq 0\leq n \leq d,\ F \in \hat{\quo{G}}_{i_n}^{(0)} \right\},
\]
where $\hat{\quo{G}}_{i_n}^{(0)}$ denotes a set of representative for the isomorphism classes of cuspidal character sheaves for $\quo{G}_{i_n}$ (\cf \cite[3.10]{CS}).
\end{corollary}

\begin{definition}\label{definition: cuspidal Levi}
A Levi subgroup $\Ksch{L}\subseteq \Ksch{G}$ is said to be a \emph{cuspidal Levi subgroup} if there is a cuspidal coefficient system for $\Ksch{L}$; in other words, $\Ksch{L}$ is a cuspidal Levi subgroup of $\Ksch{G}$ if $\catqA^{(0)}\Ksch{L}$ is non-empty (\cf Definition~\ref{definition: cuspidal}).
\end{definition}

\begin{example}\label{example: cuspidal Levi}
The algebraic group $GL(n)_\KK$ admits only one cuspidal Levi subgroup, up to conjugacy, and that is the split torus.
The split torus $GL(1)_\KK$ in $SL(2)_\KK$ is also a cuspidal Levi subgroup, as is $SL(2)_\KK$ itself, as we shall see in Section~\ref{section: examples}. The cuspidal Levi subgroups of $Sp(4)$ are described in Section~\ref{section: examples}.
\end{example}


\section{Admissible coefficient systems}\label{section: admissible}

Throughout Section~\ref{section: admissible} we assume $\KK$ is a maximal unramified closure of a local field with finite residue field. We note that such a field is strictly henselian and that the residue field is an algebraic closure of the finite field. Let $\Ksch{G}$ be a connected reductive linear algebraic group over $\KK$ satisfying the conditions of Section~\ref{subsection: stabilizers}.  

Let $\Ksch{\sigma} : \Ksch{P} \to \Ksch{G}$ be a parabolic subgroup with Levi component $\Ksch{L}$. In Section~\ref{subsection: parabolic subgroups} we define a parabolic induction function $\ind^{\Ksch{G}}_{\Ksch{P}}$ taking weakly-equivariant objects of $\catqC\Ksch{L}$ to weakly-equivariant objects of $\catqC\Ksch{G}$. It should be noted that, in constrast to parabolic restriction $\ind^{\Ksch{G}}_{\Ksch{P}}$ (Section~\ref{subsection: parabolic restriction}), conjugation $\quo{m}(g)^*$ (Section~\ref{subsection: weakly-equivariant objects}), and Frobenius $\frob^*$ (Section~\ref{section: frobenius} below), parabolic induction $\ind^{\Ksch{G}}_{\Ksch{P}}$ is a function, not a functor. This is because we can only apply our definition to weakly-equivariant objects in $\catqC\Ksch{L}$, and we have not constructed a category of weakly-equivariant objects in $\catqC\Ksch{L}$. This is not to say that such a construction is not possible. Note that we took pains in Section~\ref{subsection: parabolic induction on reductive quotients} to treat $\ind_{i\leq j}$ as a functor on the category of equivariant perverse sheaves, not just a function, and we will use that improvement below. However, the definition of weakly-equivariant objects, while built upon that of equivariant perverse sheaves, is considerably less sophisticated and almost certainly not the correct definition upon which to try to build a well-behaved category. Nevertheless, our parabolic induction $\ind^{\Ksch{G}}_{\Ksch{P}}$, as a function, is all that is needed to define admissible coefficient systems (see Definition~\ref{definition: admissible} below) for the same reason that Lusztig's parabolic induction -- also a function, not a functor -- suffices to define character sheaves.


\subsection{Parabolic subgroups}\label{subsection: parabolic subgroups}

Let $\Ksch{G} \to \Ksch{G/P}$ be the cokernel of $\sigma : \Ksch{P} \to \Ksch{G}$, let $\Ksch{\pi} : \Ksch{P} \to \Ksch{L}$ be canonical quotient map and let $\Ksch{U} \hookrightarrow \Ksch{P}$ be the unipotent radical of $\Ksch{P}$. For any facet $j$ of $I(\Ksch{G},\KK)$, let $\Rsch{P}_j$ (resp. $\Rsch{L}^j$, $\Rsch{U}_j$) be the schematic closure of $\Ksch{P}$ (resp. $\Ksch{L}$, $\Ksch{U}$) in $\Rsch{G}_j$. Then $\Rsch{P}_{j}$ (resp. $\Rsch{L}^{j}$, $\Rsch{U}_{j}$) is a smooth integral model for $\Ksch{P}$ (resp. $\Rsch{L}^{j}$, $\Rsch{U}_{j}$). Since $\Ksch{\sigma}(\Rsch{P}_j(\RK)) \subseteq\Rsch{G}_j(\RK)$ and $\Ksch{\pi}(\Rsch{P}_j(\RK)) \subseteq\Rsch{L}_j(\RK)$, it follows from the Extension Principle (\cf\cite[1.7]{BT2}) that $\Ksch{\sigma} : \Ksch{P} \to\Ksch{G}$ and $\Ksch{\pi} : \Ksch{P} \to \Ksch{L}$ and
$\Ksch{U}\hookrightarrow \Ksch{P}$ extend uniquely to $\RK$-scheme morphisms $\Rsch{\sigma}_j : \Rsch{P}_j \to
\Rsch{G}_j$ and $\Rsch{\pi}_j : \Rsch{P}_j \to \Rsch{L}^j$ and $\Rsch{U}_j \to \Rsch{P}_j$ such that the squares commute in the following diagramme.
\begin{equation}
	\xymatrix{
\ar[d] \Ksch{G} &  \ar[l]_{\sigma}  \Ksch{P} \ar[d] \ar[r]^{\pi} & \Ksch{L} \ar[d]\\ 
\Rsch{L}^j   & \ar[l]_{\sigma_j} \Rsch{P}_j \ar[r]^{\pi_j} & \Rsch{G}_j 
	}
\end{equation}

\begin{definition}\label{definition: Di}
Let $i$ be an arbitrary facet of $I(\Ksch{G},\KK)$. Let $D_i(\Ksch{G},\Ksch{P},\KK)$ denote the elements of the double coset space $\Rsch{G}_i(\RK) \backslash \Ksch{G}(\KK) / \Ksch{L}(\KK)$ represented by $g\in \Ksch{G}(\KK)$ with the following properties, where $\Rsch{L}^{ig}$ denotes the integral closure of $\Ksch{L}$ in $\Rsch{G}_{ig}$: $\Rsch{L}^{ig}$ is a smooth integral model for $\Ksch{L}$, $\Rsch{L}^{ig}(\RK)$ is a parahoric subgroup of $\Ksch{L}(\KK)$, and $\quo{L}^{ig}$ is a levi subgroup of $\quo{G}_{ig}$. It follows from \cite[9.22]{Lan1} that there is a unique facet $i^g_\Ksch{P}$ in the star of $ig$ such that $\quo{G}_{i^g_\Ksch{P}} = \quo{L}^{ig}$ and $\quo{G}_{ig\leq i^g_\Ksch{P}}$ coincides with the schematic intersection of $\Ksch{P}$ with $\quo{G}_{ig}$. Let $i^g_\Ksch{L}$ denote the facet of $I(\Ksch{L},\KK)$ such that $\Rsch{L}^{ig} = \Rsch{L}_{i^g_\Ksch{L}}$.
Notice that if $g$ satisfies this condition then so does $hgl$ for all $h\in \Rsch{G}_i(\RK)$ and $l\in \Ksch{L}(\KK)$.
\end{definition}

We begin by remarking that $D_i(\Ksch{G},\Ksch{P},\KK)$ is finite.
First, notice that $\Rsch{G}_i(\RK)\backslash \Ksch{G}(\KK)/\Ksch{P}(\KK)$ is finite. Thus, the image of $D_i(\Ksch{G},\Ksch{P},\KK)$ under the surjection
\begin{eqnarray*}
\Rsch{G}_i(\RK)\backslash \Ksch{G}(\KK)/\Ksch{L}(\KK) &\to&
\Rsch{G}_i(\RK)\backslash \Ksch{G}(\KK)/\Ksch{P}(\KK)\\
\Rsch{G}_i(\RK) g \Ksch{L}(\KK) &\mapsto& \Rsch{G}_i(\RK) g
\Ksch{P}(\KK)
\end{eqnarray*}
is finite. Suppose now that $g$ represents a double coset in $D_i(\Ksch{G},\Ksch{P},\KK)$. Then each element in the pre-image of $\Rsch{G}_i(\RK) g \Ksch{P}(\KK)$ under the map above is represented by $gu$, for some $u \in \Ksch{U}(\KK)$. If $gu$ represents an element of $D_i(\Ksch{G},\Ksch{P},\KK)$ then $g$ and $gu$ have the properties listed in Definition~\ref{definition: Di}.  This implies $u \in \Rsch{U}_{ig}(\RK)$, in which case $gu=hg$ for some $h\in \Rsch{G}_i(\RK)$.  Thus, the intersection of the pre-image of $\Rsch{G}_i(\RK) g \Ksch{P}(\KK)$ with $D_i(\Ksch{G},\Ksch{P},\KK)$ is a singleton. It follows that $D_i(\Ksch{G},\Ksch{P},\KK)$ is finite.

In order to define parabolic induction functors, we now begin the process of picking specific representatives in $\Ksch{G}(\KK)$ for the double cosets appearing in Definition~\ref{definition: Di}. We begin by recalling some basic notions and establishing some notation. Let $\Ksch{T}$ be a maximal $\KK$-split torus of $\Ksch{G}$ and let $A(\Ksch{G},\Ksch{T},\KK)$ be the apartment for $\Ksch{T}$. Let $\Rsch{T}_0$ be the N\'eron model for
$\Ksch{T}$ and let $W(\Ksch{G},\Rsch{T},\KK)$ be the associated affine Weyl group associated to the pair
$(\Ksch{G}(\KK),\Rsch{T}_0(\RK))$; that is, $W(\Ksch{G},\Ksch{T},\KK) = \Ksch{N}(\KK)/\Rsch{T}_0(\RK)$,
where $\Ksch{N}$ is the normalizer of $\Ksch{T}$ in $\Ksch{G}$. If $i$ is a facet of $A(\Ksch{G},\Ksch{T},\KK)$ let
$W_i(\Ksch{G},\Ksch{T},\KK)$ denote the stabiliser of $i$ in $W(\Ksch{G},\Ksch{T},\KK)$. Let $\dot{W}(\Ksch{G},\Ksch{T},\KK)$ be a set of representatives for $W(\Ksch{G},\Ksch{T},\KK)$ contained in $\Ksch{N}(\KK)$ and chosen so that $vw$ is represented by $\dot{v} \dot{w}$ when the length of $vw$ equals the length of $v$ plus the length of $w$ (\cf \cite[5.2]{Morris}). Let $\dot{W}_i(\Ksch{G},\Ksch{T},\KK)$ denote the subset corresponding to $W_i(\Ksch{G},\Ksch{T},\KK)$.

Suppose $i$ and $j$ are facets of $I(\Ksch{G},\KK)$ such that $i \leq j$. Let $d_i(\Ksch{G},\Ksch{P},\KK)$ (resp.
$d_j(\Ksch{G},\Ksch{P},\KK)$) be a set of representatives for the double coset space $D_i(\Ksch{G},\Ksch{P},\KK)$ (resp.
$D_i(\Ksch{G},\Ksch{P},\KK)$). The map $\Rsch{f}_{i\leq j} : \Rsch{G}_j \to \Rsch{G}_i$
(\cf Section~\ref{subsection: restriction between reductive quotients}) 
defines an inclusion of points $\Rsch{G}_j(\RK) \subseteq \Rsch{G}_i(\RK)$ which in turn defines a surjection
\begin{eqnarray*}
\Rsch{G}_j(\RK) \backslash \Ksch{G}(\KK) / \Ksch{L}(\KK)
&\to&
\Rsch{G}_i(\RK) \backslash \Ksch{G}(\KK) / \Ksch{L}(\KK)\\
\Rsch{G}_j(\RK) y \Ksch{L}(\KK) &\mapsto& \Rsch{G}_i(\RK) y
\Ksch{L}(\KK).
\end{eqnarray*}
If $y\in d_j(\Ksch{G},\Ksch{P},\KK)$ then $y$ satisfies the conditions of Definition~\ref{definition: Di}; in particular,
$\Rsch{P}_{jy} \hookrightarrow \Rsch{G}_{jy}$ is projective. Now $\Rsch{f}_{iy\leq  jy}$ induces $\Rsch{P}_{jy} \hookrightarrow \Rsch{P}_{iy}$ so $\Rsch{P}_{iy} \hookrightarrow \Rsch{G}_{iy}$ is projective. The other conditions appearing in Definition~\ref{definition: Di} are also satisfied, so $y$ represents a double coset in $D_i(\Ksch{G},\Ksch{P},\KK)$. Therefore, the surjection above restricts to a surjection
\begin{equation}\label{equation: Dij}
D_j(\Ksch{G},\Ksch{P},\KK) \to D_i(\Ksch{G},\Ksch{P},\KK)
\end{equation}
which in turn defines
\begin{equation}\label{equation: dij}
d_j(\Ksch{G},\Ksch{P},\KK) \to d_i(\Ksch{G},\Ksch{P},\KK).
\end{equation}
Let $x$ be the image of $y$ under Equation~\ref{equation: dij}. The fibre of Equation~\ref{equation: Dij} above the double coset represented by $x$ is
\begin{equation}
    \left\{ \Rsch{G}_j(\RK) ux \Ksch{L}(\KK)
    \tq u\in \Rsch{G}_{jx}(\RK) \backslash \Rsch{G}_{ix}(\RK)
    /(\Rsch{G}_{ix}(\RK) \cap \Ksch{L}(\KK)\right\},
\end{equation}
which we now denote $D^x_{i\leq j}(\Ksch{G},\Ksch{P},\KK)$. Observe that
\begin{equation}
\Rsch{G}_{ix}(\RK) \cap \Ksch{L}(\KK) = \Rsch{L}^{ix}(\RK) = \Rsch{L}_{i^x_\Ksch{L}}(\RK)
\end{equation}
and this is a parahoric subgroup of $\Ksch{L}(\KK)$. Recall also that $ix \leq i^x_\Ksch{P}$ (see Definition~\ref{definition: Di}). Using the affine Bruhat Decomposition (\cf\cite[3.22]{Morris}) we find $D^x_{i\leq j}(\Ksch{G},\Ksch{P},\KK)$ is in bijection with
\begin{equation}
W_{jx}(\Ksch{G},\Ksch{T},\KK) \backslash
W_{ix}(\Ksch{G},\Ksch{T},\KK) /
W_{i^x_\Ksch{L}}(\Ksch{L},\Ksch{T},\KK).
\end{equation}
Following \cite[C.1.1]{V2} we find a unique element of minimal length, called a distinguished element, in each double coset
above, and let $d^x_{i\leq j}(\Ksch{G},\Ksch{P},\KK)$ be the elements of $\dot{W}(\Ksch{G},\Ksch{T},\KK)$ corresponding to
distinguished elements representing the double coset space above. Thus, there is a unique $u \in d_{i\leq j}(\Ksch{G},\Ksch{P},\KK)_x$ such that $u x  = h y l$, for some $l\in \Ksch{L}(\KK)$ and some $h \in \Rsch{G}_j(\RK)$. In fact, we can say much more, as Lemma~\ref{lemma: di} shows.

\begin{lemma}\label{lemma: di}
Keep notation as above. There exists a collection of sets $d_i(\Ksch{G},\Ksch{P},\KK)$, as $i$ ranges over all facets of
$I(\ksch{G},\KK)$, consisting of representatives for the double coset spaces $D_i(\Ksch{G},\Ksch{P},\KK)$, with the following properties:
\begin{enumerate}
\item[(a)]
Let $i$ and $j$ be facets of $I(\Ksch{G},\KK)$ with $i \leq j$. If $g\in d_j(\Ksch{G},\Ksch{P},\KK)$ then $ux = gl$, where $x$ is the image of $g$ under the map of Equation~\ref{equation: dij} and $u$ is an element of $d^x_{i\leq j}(\Ksch{G},\Ksch{P},\KK)$. Moreover, $u$ and $l$ (and $x$) are determined uniquely by $g$.
\item[(b)]
Let $i$ be a facet of $I(\Ksch{G},\KK)$ and let $g$ be an element of $\Ksch{G}(\KK)$. If $y$ is an element of
$d_{ig}(\Ksch{G},\Ksch{P},\KK)$ then there is unique $x$ in $d_i(\Ksch{G},\Ksch{P},\KK)$ and $h$ in $\Rsch{G}_i(\RK)$ such that $g y = h x$.
\end{enumerate}
\end{lemma}

\begin{proof}
Let $i$ be a facet of $A(\Ksch{G},\Ksch{T},\KK)$.
Let $d_i(\Ksch{G},\Ksch{P},\KK)$ be a set representatives for $D_i(\Ksch{G},\Ksch{P},\KK)$ chosen from
$\dot{W}(\Ksch{G},\Ksch{T},\KK)$. Suppose $i$ and $j$ are facets
of $A(\Ksch{G},\Ksch{T},\KK)$ such that $i \leq j$. Now, then
$x$, $g$ and $u$ are elements of $\dot{W}(\Ksch{G},\Ksch{T},\KK)$
and $g^{-1}ux$ is contained in $\dot{W}(\Ksch{L},\Ksch{T},\KK)$.
Set $l = g^{-1}ux$; we now have
\begin{equation}\label{equation: dij 1}
    ux = gl,
\end{equation}
as promised. Now, let $i$ be any facet of $I(\Ksch{G},\KK)$.
There is some $z$ in $\Ksch{G}(\KK)/\Ksch{N}(\KK)$ such that $i$
is contained in the apartment $A(\Ksch{G},\Ksch{T}^z,\KK)$. Since
$zi$ is a facet of $A(\Ksch{G},\Ksch{T},\KK)$ we set
\begin{equation}
d_{i}(\Ksch{G},\Ksch{P},\KK) = z^{-1}
d_{zi}(\Ksch{G},\Ksch{P},\KK).
\end{equation}
If $i' \leq j'$ and $j'$ is a facet of
$A(\Ksch{G},\Ksch{T}^z,\KK)$ then $i'$ is a facet of
$A(\Ksch{G},\Ksch{T}^z,\KK)$ also. Suppose $x'\in
d_{j'}(\Ksch{G},\Ksch{P},\KK)$ and ${g'}$ is the image of ${x'}$
under the map to $d_{i'}(\Ksch{G},\Ksch{P},\KK)$. Set $i = zi'$,
$j = zj'$, $x=zx'$ and $g=zg'$. Then $z^{-1}u'z$ is an element of
$d_{i\leq j}(\Ksch{G},\Ksch{P},\KK)$ so set $u = zu'z^{-1}$. Then $ux
= gl$ by Equation~\ref{equation: dij 1}, so $(z u' z^{-1}) (zx') =
(zg') l$, so $u' x'=  g' l$, as desired.

The proof of part (b) is omitted.
\end{proof}


\subsection{Local parabolic induction}\label{subsection: local parabolic induction}

\begin{definition}\label{definition: local parabolic induction}
Let $\Ksch{P}$ be a parabolic subgroup of $\Ksch{G}$ with Levi component $\Ksch{L}$. Let $i$ be any facet of $I(\Ksch{G}, \KK)$. Let $g$ be an element of $\Ksch{G}(\KK)$ such that $ig$ satisfies the conditions appearing in Definition~\ref{definition: Di}. Let $\ind^{\quo{G}_{ig}}_{\quo{L}^{ig}}$ denote the functor from the category of equivariant perverse sheaves whose objects are finite direct sums of character sheaves on $\quo{L}^{ig}$ to the category of equivariant perverse sheaves whose objects are finite direct sums of character sheaves on $\quo{G}_{ig}$ given by $\ind_{ig\leq i^g_\Ksch{P}}$. For any $B\in \obj \catqC\Ksch{L}$ we will write $B^{ig}$ for $B_{i^g_\Ksch{L}}$. Thus, $\ind^{\quo{G}_{ig}}_{\quo{L}^{ig}} B^{ig} = \ind_{ig \leq i^g_\Ksch{P}} B_{i^g_\Ksch{L}}$. (This is not a typo; observe that $\quo{G}_{i^g_\Ksch{P}} = \quo{L}_{i^g_\Ksch{L}}$.) Likewise for morphisms in $\catqC\Ksch{L}$.
\end{definition}


For $j$ a facet of $I(\Ksch{G},\KK)$, let $\Sigma_j$ be the set of affine
roots vanishing on $j$. Let $l$ be a second facet of $I(\Ksch{G},\KK)$. We will 
denote by $j\wedge l$ the facet of $I(\Ksch{G},\KK)$ defined by
$\Sigma_{j\wedge l}=\Sigma_j\cap\Sigma_l$.
 
\begin{lemma}\label{lemma: restriction and induction}
Let $i$, $j$ and $l$ be facets of $I(\Ksch{G},\KK)$ such that $i\leq j$ 
and $i \leq l$.  There is a choice of representatives $d^i_{j,l}$ for 
$\Rsch{G}_{i\leq j}(\RK) \backslash \Rsch{G}_i(\RK) / \Rsch{G}_{i\leq l}(\RK)$ 
such that
\[
\res_{i\leq j} \ind_{i\leq l} F 
\iso 
\sum_{x\in d^i_{j,l}} \quo{m}(x^{-1})^*\ \ind_{jx\leq jx\wedge l} \res_{l\leq l\wedge jx} F
\]
and the diagramme
\[
	\xymatrix{
	\ar[d]^{\res_{i\leq j\leq k} \ind_{i\leq l} F}
\res_{j\leq k} \res_{i\leq j} \ind_{i\leq l} F 
	\ar[r]
& 
\res_{j\leq k} \sum_{x\in d^i_{j,l}} \quo{m}(x^{-1})^*\ \ind_{jx\leq jx\wedge l} \res_{l\leq l\wedge jx} F
	\ar[d]
\\
\res_{i\leq k} \ind_{i\leq l} F
	\ar[r]
&
\sum_{x\in d^i_{k,l}} \quo{m}(x^{-1})^*\ \ind_{kx\leq kx\wedge l} \res_{l\leq l\wedge kx} F
	}
\]
commutes, for any character sheaf $F$ on $\quo{G}_l$. 
\end{lemma}

\begin{proof}
Recall that $\KK$ is an unramified closure of a $p$-adic field.
We choose a uniformizer and, as in \cite[C.1.1]{V2}, we find a unique element of minimal length, called a
distinguished element, in each double coset above, and let $d^i_{k,l}$ be the 
set of distinguished elements representing the double coset space above.
Then use \cite[Prop~10.1.2]{MS} (which gives considerably more information concerning the isomorphism above than \cite[Prop~15.2]{CS}). 
Note that the proofs of \cite[Prop~10.1.2]{MS} and \cite[Prop~15.2]{CS}, while different, both depend on the fact that $\kK$ is the algebraic closure of a finite field.
\end{proof}

\subsection{Parabolic induction}\label{subsection: parabolic induction}

\begin{proposition}\label{proposition: parabolic induction}
Let $\Ksch{P}$ be a parabolic subgroup $\Ksch{G}$ with Levi component $\Ksch{L}$. If $B$ is a weakly-equivariant admissible object of $\catqC\Ksch{L}$ then there is a weakly-equivariant object $\ind^{\Ksch{G}}_{\Ksch{P}} B$ of $\catqC\Ksch{G}$ such that 
\[
    (\ind^{\Ksch{G}}_{\Ksch{P}} B)_i
 \ceq
    \sum\limits_{g\in d_i(\Ksch{G},\Ksch{P},\KK)}
    \quo{m}(g^{-1})_i^*\
    \ind^{\quo{G}_{ig}}_{\quo{L}^{ig}}
    \ B^{ig}.
\]
for each facet $i$ of $I(\Ksch{G},\KK)$.
\end{proposition}

\begin{proof}
We must begin by defining $(\ind^{\Ksch{G}}_{\Ksch{P}} B)_{i\leq j}$ for every $i\leq j$ in $I(\Ksch{G},\KK)$.To that end, we fix one such pair of facets and consider $\res_{i\leq j}(\ind^{\Ksch{G}}_{\Ksch{P}} B)_i$. 
\begin{eqnarray*}
\res_{i\leq j}(\ind^{\Ksch{G}}_{\Ksch{P}} B)_i
	&=& \res_{i\leq j}  \sum\limits_{y\in d_i(\Ksch{G},\Ksch{P},\KK)} 
		\quo{m}(g^{-1})_i^*\
  		\ind^{\quo{G}_{ig}}_{\quo{L}^{ig}} \ B^{ig}\\
    	&=& \sum\limits_{g\in d_i(\Ksch{G},\Ksch{P},\KK)} 
		\res^{\quo{G}_i}_{\quo{G}_j} \ \quo{m}(g^{-1})_i^*\
    		\ind^{\quo{G}_{ig}}_{\quo{L}^{ig}} \ B^{ig}.
\end{eqnarray*}
Let
\begin{equation}\label{equation: pi.0}
\xymatrix{
\sum\limits_{g\in d_i(\Ksch{G},\Ksch{P},\KK)} 
		\res^{\quo{G}_i}_{\quo{G}_j} \ \quo{m}(g^{-1})_i^*\
    		\ind^{\quo{G}_{ig}}_{\quo{L}^{ig}} \ B^{ig}
\ar[d]^{1.}
\\
\sum\limits_{g\in d_i(\Ksch{G},\Ksch{P},\KK)} 
		\quo{m}(g^{-1})_i^*\ \res^{\quo{G}_{ig}}_{\quo{G}_{jg}}
    		\ind^{\quo{G}_{ig}}_{\quo{L}^{ig}} \ B^{ig}
}
\end{equation}
be the isomorphism in $\catM_{\quo{G}_j}\quo{G}_j$ determined by the natural isomorphism
\begin{equation}\label{equation: pi.1}
		\res^{\quo{G}_i}_{\quo{G}_j}\ \quo{m}(g^{-1})_i^*
\iso
		\quo{m}(g^{-1})_i^*\ \res^{\quo{G}_{ig}}_{\quo{G}_{jg}}
\end{equation}
of Lemma~\ref{lemma: conjugation and restriction}. Now, since $B$ is an object in $\catqC\Ksch{L}$ it follows that $B^{ig} = B_{i^g_\Ksch{L}}$ is a finite direct sum of character sheaves. Recall that $\quo{L}^{ig} = \quo{G}_{i^g_\Ksch{P}}$. By Lemma~\ref{lemma: restriction and induction} (MacKey's formula for character sheaves) we have an isomorphism 
\begin{equation}\label{equation: pi.2}
\res^{\quo{G}_{ig}}_{\quo{G}_{jg}}
    		\ind^{\quo{G}_{ig}}_{\quo{G}_{i^g_\Ksch{P}}} \ B^{ig}
\iso
\sum_{u\in d^g_{i\leq j}(\Ksch{G},\Ksch{P},\KK)}
    \quo{m}(u^{-1})^*
    \ind^{\quo{G}_{jgu}}_{\quo{G}_{j^{gu}_\Ksch{P}}}\
    \res^{\quo{G}_{i^g_\Ksch{P}}}_{\quo{G}_{j^{gu}_\Ksch{P}}}\ B^{ig}.
\end{equation}
which defines
\begin{equation}\label{equation: pi.3}
\xymatrix{
\sum\limits_{g\in d_i(\Ksch{G},\Ksch{P},\KK)} 
		\quo{m}(g^{-1})_i^*\ \res_{ig\leq jg} \
    		\ind^{\quo{G}_{ig}}_{\quo{L}^{ig}} \ B^{ig}
\ar[d]^{2.}
\\
  \sum\limits_{g\in d_i(\Ksch{G},\Ksch{P},\KK)} 
		\quo{m}(g^{-1})_i^*\ \sum_{u\in d^g_{i\leq j}(\Ksch{G},\Ksch{P},\KK)}
		\quo{m}(u^{-1})_j^* \ind^{\quo{G}_{jgu}}_{\quo{G}_{j^{gu}_\Ksch{P}}}\
		\res^{\quo{G}_{i^g_\Ksch{P}}}_{\quo{G}_{j^{gu}_\Ksch{P}}}\ B^{ig}.
}
\end{equation}
(Here we have used the notation of Lemma~\ref{lemma: di}.)
Let
\begin{equation}\label{equation: pi.4}
\xymatrix{
\sum\limits_{g\in d_i(\Ksch{G},\Ksch{P},\KK)} 
		\quo{m}(g^{-1})_i^*\ \sum_{u\in d^g_{i\leq j}(\Ksch{G},\Ksch{P},\KK)}
		\quo{m}(u^{-1})_j^* \ind^{\quo{G}_{jgu}}_{\quo{G}_{j^{gu}_\Ksch{P}}}\
		\res^{\quo{G}_{i^g_\Ksch{P}}}_{\quo{G}_{j^{gu}_\Ksch{P}}}\ B^{ig}
\ar[d]^{3.}
\\
\sum\limits_{g\in d_i(\Ksch{G},\Ksch{P},\KK) \atop u\in d^g_{i\leq j}(\Ksch{G},\Ksch{P},\KK)}
		\quo{m}((gu)^{-1})_j^* \ind^{\quo{G}_{jgu}}_{\quo{G}_{j^{gu}_\Ksch{P}}}\
		\res^{\quo{G}_{i^g_\Ksch{P}}}_{\quo{G}_{j^{gu}_\Ksch{P}}}\ B^{ig}.
}
\end{equation}
be the isomorphism in determined by the natural isomorphism
\begin{equation}\label{equation: pi.5}
\quo{m}(g^{-1})^*\ \quo{m}(u^{-1})^*\iso \quo{m}((gu)^{-1})^* .
\end{equation}
Since $u\in \Rsch{G}_{jg}(\RK)$ we have $\Rsch{G}_{jg} = \Rsch{G}_{jgu}$; it follows that the integral closure of
$\Ksch{L}$ in $\Rsch{G}_{jg}$ is the integral closure of $\Ksch{L}$ in $\Rsch{G}_{jgu}$, which we denote
$\Rsch{L}_{j^{gu}_P}$. Moreover, since $i \leq j$ we have $\Rsch{G}_i \supseteq \Rsch{G}_j$ and $\Rsch{G}_{ig} \supseteq
\Rsch{G}_{jg}$. Also, $\Rsch{G}_{ig} = \Rsch{G}_j^{gu}$, so it follows that the integral closure of $\Ksch{L}$ in $\Rsch{G}_i^g$ is the integral closure of $\Ksch{L}$ in $\Rsch{G}_i^{gu}$. Therefore,
\begin{equation}\label{equation: pi.6}
\sum\limits_{g\in d_i(\Ksch{G},\Ksch{P},\KK) \atop u\in d^g_{i\leq j}(\Ksch{G},\Ksch{P},\KK)}
		\quo{m}((gu)^{-1})_j^* \ind^{\quo{G}_{jgu}}_{\quo{G}_{j^{gu}_\Ksch{P}}}\
		\res^{\quo{G}_{i^g_\Ksch{P}}}_{\quo{G}_{j^{gu}_\Ksch{P}}}\ B^{ig}
=
\sum\limits_{g\in d_i(\Ksch{G},\Ksch{P},\KK) \atop u\in d^g_{i\leq j}(\Ksch{G},\Ksch{P},\KK)}
		\quo{m}((gu)^{-1})_j^* \ind^{\quo{G}_{jgu}}_{\quo{L}^{jgu}}\
		\res^{\quo{L}^{igu}}_{\quo{L}^{jgu}} \ B^{igu}.
\end{equation}
By Lemma~\ref{lemma: di}, for each $g$ and $u$ as above there is a unique $h \in d_j(\Ksch{G},\Ksch{P},\KK)$ and $l = l^i_j(h)$ such that $g u = h l$. Thus,
\begin{equation}\label{equation: pi.7}
     \sum\limits_{g\in d_i(\Ksch{G},\Ksch{P},\KK) \atop u\in d^g_{i\leq j}(\Ksch{G},\Ksch{P},\KK)}
		\quo{m}((gu)^{-1})_j^* \ind^{\quo{G}_{jgu}}_{\quo{L}^{jgu}}\
		\res^{\quo{L}^{igu}}_{\quo{L}^{jgu}} \ B^{igu}
=
\sum\limits_{h\in d_{j}(\Ksch{G},\Ksch{P},\KK)}
    \quo{m}((h l)^{-1})_j^*
    \ind^{\quo{G}_{jhl}}_{\quo{L}^{j h l}}\
     \res^{\quo{L}^{i h l}}_{\quo{L}_{j h l}}\
    B^{ihl}.
\end{equation}
Now, the natural isomorphism
\begin{equation}\label{equation: pi.8}
\quo{m}((hl)^{-1})^* \iso \quo{m}(h^{-1})^*\ \quo{m}(l^{-1})^* 
\end{equation}
defines
\begin{equation}\label{equation: pi.9}
\xymatrix{
\sum\limits_{h\in d_{j}(\Ksch{G},\Ksch{P},\KK)}
    \quo{m}((h l)^{-1})_j^*
    \ind^{\quo{G}_{jhl}}_{\quo{L}^{j h l}}\
     \res^{\quo{L}^{i h l}}_{\quo{L}_{j h l}}\
    B^{ihl}
\ar[d]^{4.}
\\
\sum\limits_{h\in d_{j}(\Ksch{G},\Ksch{P},\KK)}
    \quo{m}(h^{-1})^*\ \quo{m}(l^{-1})^*\
    \ind^{\quo{G}_{jhl}}_{\quo{L}^{jhl}}\
    \res^{\quo{L}^{ihl}}_{\quo{L}^{jhl}}\
    B^{ihl}.
}
\end{equation}
Lemma~\ref{lemma: conjugation and induction} gives
\begin{equation}
\quo{m}(l^{-1})^*\
    \ind^{\quo{G}_{jhl}}_{\quo{L}^{jhl}}
   \iso
    \ind^{\quo{G}_{jh}}_{\quo{L}^{jh}}\
   \quo{m}(l^{-1})^*\
\end{equation}
and therefore defines
\begin{equation}\label{equation: pi.10}
\xymatrix{
\sum\limits_{h\in d_{j}(\Ksch{G},\Ksch{P},\KK)}
    \quo{m}(h^{-1})^*\ \quo{m}(l^{-1})^*\
    \ind^{\quo{G}_{jhl}}_{\quo{L}^{jhl}}\
    \res^{\quo{L}^{ihl}}_{\quo{L}^{jhl}}\
    B^{ihl}
\ar[d]^{5.}
\\
\sum\limits_{h\in d_{j}(\Ksch{G},\Ksch{P},\KK)}
    \quo{m}(h^{-1})^*\ 
    \ind^{\quo{G}_{jh}}_{\quo{L}^{jh}}\
    \quo{m}(l^{-1})^*\ 
    \res^{\quo{L}^{ihl}}_{\quo{L}^{jhl}}\ B^{ihl}.
}
\end{equation}
Since $B$ is an object of $\catqD\Ksch{L}$ (see Definition~\ref{definition: catqD}(obj)) we have the isomorphism
\begin{equation}\label{equation: pi.11}
\res^{\quo{L}^{ihl}}_{\quo{L}^{jhl}}\ B^{ihl} \iso B^{jhl}
\end{equation}
which defines
\begin{equation}\label{equation: pi.12}
\xymatrix{
\sum\limits_{h\in d_{j}(\Ksch{G},\Ksch{P},\KK)}
    \quo{m}(h^{-1})^*\ 
    \ind^{\quo{G}_{jh}}_{\quo{L}^{jh}}\
    \quo{m}(l^{-1})^*\ 
    \res^{\quo{L}^{ihl}}_{\quo{L}^{jhl}}\ B^{ihl}
\ar[d]^{6.}
\\
\sum\limits_{h\in d_{j}(\Ksch{G},\Ksch{P},\KK)}
    \quo{m}(h^{-1})^*\ 
    \ind^{\quo{G}_{jh}}_{\quo{L}^{jh}}\
    \quo{m}(l^{-1})^*\ B^{jhl}
}
\end{equation}
Since $B$ is weakly-equivariant (see Definition~\ref{definition: weakly-equivariant}) and $l\in \Ksch{L}(\KK)$, we have the isomorphism
\begin{equation}\label{equation: pi.13}
B^{ihl} \iso B^{ih}
\end{equation}
which defines
\begin{equation}\label{equation: pi.14}
\xymatrix{
\sum\limits_{h\in d_{j}(\Ksch{G},\Ksch{P},\KK)}
    \quo{m}(h^{-1})^*\ 
    \ind^{\quo{G}_{jh}}_{\quo{L}^{jh}}\
    \quo{m}(l^{-1})^*\ 
    \res^{\quo{L}^{ihl}}_{\quo{L}^{jhl}}\ B^{jhl}
\ar[d]^{7.}
\\
\sum\limits_{h\in d_{j}(\Ksch{G},\Ksch{P},\KK)}
    \quo{m}(h^{-1})^*\ 
    \ind^{\quo{G}_{jh}}_{\quo{L}^{jh}}\
    B^{jh}.
}
\end{equation}
Since this last expression is precisely, $(\ind^{\Ksch{G}}_{\Ksch{P}} B)_j$, composing the isomorphisms of Equations~\ref{equation: pi.0}, \ref{equation: pi.3}, \ref{equation: pi.4}, \ref{equation: pi.9}, \ref{equation: pi.10}, \ref{equation: pi.12}, \ref{equation: pi.14} defines an isomorphism
\begin{equation}\label{equation: pi.15}
(\ind^{\Ksch{G}}_{\Ksch{P}} B)_{i\leq j} : \res_{i\leq j} (\ind^{\Ksch{G}}_{\Ksch{P}} B)_i \to (\ind^{\Ksch{G}}_{\Ksch{P}} B)_j.
\end{equation}

Now we must show that the family $((\ind^{\Ksch{G}}_{\Ksch{P}} B)_i, (\ind^{\Ksch{G}}_{\Ksch{P}} B)_{i\leq j} )_{i,j}$ satisfies the condition of Definition~\ref{definition: catqD}(obj). Inspecting the isomorphisms appearing in the definition of $(\ind^{\Ksch{G}}_{\Ksch{P}} B)_{i\leq j}$, is clear that $(\ind^{\Ksch{G}}_{\Ksch{P}} B)_{i\leq i} = \id_{(\ind^{\Ksch{G}}_{\Ksch{P}} B)_i}$, so we turn now to the second condition in Definition~\ref{definition: catqD}(obj).
Let $i$, $j$ and $k$ be facets of $I(\Ksch{G},\KK)$ and consider the following diagramme.
\begin{equation}
	\xymatrix{
	\ar[d]_{\res_{i\leq j\leq k} (\ind^{\Ksch{G}}_{\Ksch{P}} B)_i}
\res_{j\leq k}\res_{i\leq j}(\ind^{\Ksch{G}}_{\Ksch{P}} B)_i
	\ar[rrr]^{\res_{j\leq k} (\ind^{\Ksch{G}}_{\Ksch{P}} B)_{i\leq j}}
&&&
\res_{j\leq k} (\ind^{\Ksch{G}}_{\Ksch{P}} B)_j
	\ar[d]^{(\ind^{\Ksch{G}}_{\Ksch{P}} B)_{j\leq k}}
\\
\res_{i\leq k}(\ind^{\Ksch{G}}_{\Ksch{P}} B)_k
	\ar[rrr]^{(\ind^{\Ksch{G}}_{\Ksch{P}} B)_{i\leq k}}
&&&
(\ind^{\Ksch{G}}_{\Ksch{P}} B)_k
	}
\end{equation}
In order to see that this diagramme commutes, it is sufficient to observe that each isomorphism appearing in the definition of 
$(\ind^{\Ksch{G}}_{\Ksch{P}} B)_{i\leq j}$ (and $(\ind^{\Ksch{G}}_{\Ksch{P}} B)_{j\leq k}$ and $(\ind^{\Ksch{G}}_{\Ksch{P}} B)_{i\leq k}$) is compatible with the restriction functors of Section~\ref{subsection: restriction between reductive quotients}. Fortunately, this work has already been done in various lemmata, in anticipation of this need:
for isomorphism $1.$ of Equation~\ref{equation: pi.0} use Lemma~\ref{lemma: conjugation and restriction}; 
for isomorphism $2.$ of Equation~\ref{equation: pi.3} use Lemma~\ref{lemma: restriction and induction} and Lemma~\ref{lemma: di} (to see that $d^g_{i\leq j} = d_{j^{g}_\Ksch{P} \geq jg  \leq i^{g}_\Ksch{P}}$);
for isomorphism $3.$ of Equation~\ref{equation: pi.4} use Lemma~\ref{lemma: conjugation and restriction};
for isomorphism $4.$ of Equation~\ref{equation: pi.9} use Lemma~\ref{lemma: conjugation and restriction}; 
for isomorphism $5.$ of Equation~\ref{equation: pi.10} use Lemma~\ref{lemma: conjugation and induction};
for isomorphism $6.$ of Equation~\ref{equation: pi.12} use Definition~\ref{definition: catqD}(obj);
for isomorphism $7.$ of Equation~\ref{equation: pi.14} use Definition~\ref{definition: weakly-equivariant}.
Each of these arguments makes use of Proposition~\ref{proposition: transitive restriction on quotients}!
\end{proof}


\subsection{Transitive parabolic induction}\label{subsection: transind} 

Recall the notation of Proposition~\ref{proposition: transitive restriction}.

\begin{proposition}\label{proposition: transitive induction}
Let $\Ksch{P}$ be a parabolic subgroup of $\Ksch{G}$ with Levi component $\Ksch{L}$. Let $\Ksch{Q}$ be a parabolic subgroup of $\Ksch{L}$ with Levi component $\Ksch{M}$. Let $\Ksch{R}$ be a parabolic subgroup of $\Ksch{G}$ with Levi
component $\Ksch{M}$. If $C\in \obj\catqC\Ksch{G}$ is weakly-equivariant then
\begin{equation}
\ind^{\Ksch{G}}_{\Ksch{P}}\ \ind^{\Ksch{L}}_{\Ksch{Q}} C \iso
\ind^{\Ksch{G}}_{\Ksch{R}} C.
\end{equation}
\end{proposition}

\begin{proof}
To simplify notation, let $B = \ind^{\Ksch{L}}_{\Ksch{Q}} C$ and let $A = \ind^{\Ksch{G}}_{\Ksch{P}} B$. We
fix a facet $i$ of $I(\Ksch{G},\KK)$. Then
\begin{eqnarray*}
A_i 
    &=& (\ind^{\Ksch{G}}_{\Ksch{P}} B)_i\\
    &=& \sum_{x\in d_i(\Ksch{G},\Ksch{P},\KK)}
    \quo{m}(x^{-1})_i^*\
    \ind^{\quo{G}_{ix}}_{\quo{P}_{i^x_\Ksch{P}}}\ B_{i^x_\Ksch{P}}\\
    &=& \sum_{x\in d_i(\Ksch{G},\Ksch{P},\KK)}
    \quo{m}(x^{-1})_i^*\
    \ind^{\quo{G}_{ix}}_{\quo{P}_{i^x_\Ksch{P}}}\
    \sum_{y\in d_{i^x_\Ksch{P}}(\Ksch{L},\Ksch{Q},\KK)}
    \quo{m}(y^{-1})_{i^x_\Ksch{P}}^*\
    \ind^{\quo{L}_{i^x_\Ksch{P}y}}_{\quo{Q}_{(i^x_\Ksch{P})^y_\Ksch{Q}}}
    C_{(i^x_\Ksch{P})^y_\Ksch{Q}}.
\end{eqnarray*}
Since $\Rsch{L}_{i^x_\Ksch{P}} = \Rsch{L}^{ix}$ and $y \in \Ksch{L}(\KK)$, it follows that $\Rsch{L}_{i^x_\Ksch{P} y} =
\Rsch{L}^{ixy} = \Rsch{L}_{i^{xy}_\Ksch{P}}$. Since $\Ksch{T} \subseteq \Ksch{M}\subseteq \Ksch{L}$, the schematic closure of $\Ksch{M}$ in $\Rsch{L}_{i^{xy}_\Ksch{P}}$ is equal to the schematic closure  of $\Ksch{M}$ in $\Rsch{G}_{ixy}$, and we have $\Rsch{L}_{i^x_\Ksch{P} y} = \Rsch{L}_{i^{xy}_\Ksch{P}}$ and $\Rsch{M}_{(i^x_\Ksch{P})^y_\Ksch{Q}} =
\Rsch{M}_{i^{xy}_\Ksch{Q}}$. Using Lemma~\ref{lemma: conjugation and induction} we now have
\begin{eqnarray*}
&&\hskip-20pt
    \sum_{x\in d_i(\Ksch{G},\Ksch{P},\KK)} \quo{m}(x^{-1})_i^*\
    \ind^{\quo{G}_{ix}}_{\quo{P}_{i^x_\Ksch{P}}}\
    \sum_{y\in d_{i^x_\Ksch{P}}(\Ksch{L},\Ksch{Q},\KK)}
    \quo{m}(y^{-1})_{i^x_\Ksch{P}}^*\
    \ind^{\quo{L}_{i^x_\Ksch{P} y}}_{\quo{Q}_{(i^x_\Ksch{P})^y_\Ksch{Q}}}
    C_{(i^x_\Ksch{P})^y_\Ksch{Q}}\\
 &=&\sum_{x\in d_i(\Ksch{G},\Ksch{P},\KK) \atop y\in d_{i^x_\Ksch{P}}(\Ksch{L},\Ksch{Q},\KK)}
    \quo{m}(x^{-1})_i^*\
    \ind^{\quo{G}_{ix}}_{\quo{P}_{i^x_\Ksch{P}}}\
    \quo{m}(y^{-1})_{i^x_\Ksch{P}}^*\
    \ind^{\quo{L}_{i^{xy}_\Ksch{P}}}_{\quo{Q}_{i^{xy}_\Ksch{Q}}}
    C_{i^{xy}_\Ksch{Q}}\\
 &=&\sum_{x\in d_i(\Ksch{G},\Ksch{P},\KK) \atop y\in d_{i^x_\Ksch{L}}(\Ksch{L},\Ksch{Q},\KK)}
    \quo{m}(x^{-1})_i^*\
    \quo{m}(y^{-1})_{i^x_L}^*\
    \ind^{\quo{G}_{ixy}}_{\quo{P}_{i^{xy}_\Ksch{P}}}\
    \ind^{\quo{L}_{i^{xy}_\Ksch{P}}}_{\quo{Q}_{i^{xy}_\Ksch{Q}}}
    C_{i^{xy}_\Ksch{Q}}\\
 &=&\sum_{x\in d_i(\Ksch{G},\Ksch{P},\KK) \atop y\in d_{i^x_\Ksch{L}}(\Ksch{L},\Ksch{Q},\KK)}
    \quo{m}((xy)^{-1})_i^*\
    \ind^{\quo{G}_{ixy}}_{\quo{P}_{i^{xy}_\Ksch{P}}}\
    \ind^{\quo{L}_{i^{xy}_\Ksch{P}}}_{\quo{Q}_{i^{xy}_\Ksch{Q}}}
    C_{i^{xy}_\Ksch{Q}}.
\end{eqnarray*}
Using Propositions~\cite[4.2]{CS}, \cite[4.8(b)]{CS} and \cite[2.18(a)]{CS} we have
\begin{equation}
\ind^{\quo{G}_{ixy}}_{\quo{P}_{i^{xy}_\Ksch{P}}}\
    \ind^{\quo{L}_{{i}^{xy}_\Ksch{P}}}_{\quo{Q}_{i^{xy}_\Ksch{Q}}}
    C_{i^{xy}_\Ksch{Q}}
 = \ind^{\quo{G}_{ixy}}_{\quo{R}_{i^{xy}_\Ksch{Q}}}
    C_{i^{xy}_\Ksch{Q}}.
\end{equation}
Since $i^{xy}_\Ksch{Q} = i^{xy}_\Ksch{R}$, we have
\begin{equation}\label{equation: tr 1}
(\ind^{\Ksch{G}}_{\Ksch{P}}\ \ind^{\Ksch{L}}_{\Ksch{Q}} C)_i
 =\sum_{x\in d_i(\Ksch{G},\Ksch{P},\KK) \atop y\in d_{i^x_\Ksch{L}}(\Ksch{L},\Ksch{Q},\KK)}
    \quo{m}((xy)^{-1})_i^*\
    \ind^{\quo{G}_{i}^{xy}}_{\quo{R}_{i^{xy}_\Ksch{R}}}
    C_{i^{xy}_\Ksch{R}}.
\end{equation}
Now, the canonical surjection from $D_i(\Ksch{G},\Ksch{R},\KK) \to D_i(\Ksch{G},\Ksch{P},\KK)$ defines a surjection
$d_i(\Ksch{G},\Ksch{R},\KK) \to d_i(\Ksch{G},\Ksch{P},\KK)$; the pre-image of $x$ in $d_i(\Ksch{G},\Ksch{P},\KK)$ is exactly $d_{i^x_\Ksch{P}}(\Ksch{L},\Ksch{Q},\KK)$. For each pair $(x,y)$ in Equation~\ref{equation: tr 1} there is a unique $z$ in $d_i(\Ksch{G},\Ksch{R},\KK)$ such that $xy = h z m$ for $h$ in $\Rsch{G}_i(\RK)$ and $m\in \Ksch{M}(\KK)$. Thus,
$i^{xy}_\Ksch{R} = i^z_\Ksch{R} m$ and
\begin{eqnarray*}
 &&\hskip-20pt
    \sum_{x\in d_i(\Ksch{G},\Ksch{P},\KK)
    \atop y\in d_{i^x_L}(\Ksch{L},\Ksch{Q},\KK)}
    \quo{m}((xy)^{-1})_i^*\
    \ind^{\quo{G}_{ixy}}_{\quo{R}_{i^{xy}_\Ksch{R}}}
    C_{i^{xy}_\Ksch{R}}\\
 &=&
    \sum_{z\in d_i(\Ksch{G},\Ksch{R},\KK)}
    \quo{m}((hzm)^{-1})_i^*\
    \ind^{\quo{G}_{izm}}_{\quo{R}_{i^{z}_\Ksch{R} m}}
    C_{i^{z}_\Ksch{R} m}\\
 &=&
    \sum_{z\in d_i(\Ksch{G},\Ksch{R},\KK)}
    \quo{m}(h^{-1})^* \quo{m}(z^{-1})^* \quo{m}(m^{-1})^*
    \ind^{\quo{G}_{izm}}_{\quo{R}_{i^{z}_\Ksch{R} m}}
    C_{i^{z}_\Ksch{R} m}\\
  &=&
    \sum_{z\in d_i(\Ksch{G},\Ksch{R},\KK)}
    \quo{m}(h^{-1})^* \quo{m}(z^{-1})^*
    \ind^{\quo{G}_{iz}}_{\quo{R}_{i^{z}_\Ksch{R}}}
    \conj{m}{C}_{i^{z}_\Ksch{R}}.
\end{eqnarray*}
Since $C\in \obj\catqD\Ksch{M}$ and $h\in \Rsch{G}_i(\RK)$ we have
\begin{eqnarray*}
 &&\hskip-20pt
    \sum_{z\in d_i(\Ksch{G},\Ksch{R},\KK)}
    \quo{m}(h^{-1})^* \quo{m}(z^{-1})^*
    \ind^{\quo{G}_{iz}}_{\quo{R}_{i^{z}_\Ksch{R}}}
    \conj{m}{C}_{i^{z}_\Ksch{R}}\\
 &\iso& \sum_{z\in d_i(\Ksch{G},\Ksch{R},\KK)}
    \quo{m}(z^{-1})^*
    \ind^{\quo{G}_{iz}}_{\quo{R}_{i^{z}_\Ksch{R}}}
    {C}_{i^{z}_\Ksch{R}}\\
 &=& (\ind^{\Ksch{G}}_{\Ksch{R}} C)_i.
\end{eqnarray*}
In summary, we have shown that
\begin{equation}
(\ind^{\Ksch{G}}_{\Ksch{P}}\ \ind^{\Ksch{L}}_{\Ksch{Q}} C)_i \iso
(\ind^{\Ksch{G}}_{\Ksch{R}} C)_i.
\end{equation}
Now, as $i$ ranges over all facets of $I(\Ksch{G},\KK)$ these isomorphisms define an isomorphism in $\catqC\Ksch{G}$ (details omitted).
\end{proof}


\subsection{Admissible coefficient systems}\label{subsection: admissible}

We now come to the main definition of Section~\ref{section: admissible}.

\begin{definition}\label{definition: admissible}
An \emph{irreducible admissible coefficient system for $\Ksch{G}$} is a simple object of $\catqC\Ksch{G}$ which is a summand of $\ind^{\Ksch{G}}_{\Ksch{P}} C$ in $\catqC\Ksch{G}$ for some parabolic subgroup $\Ksch{P}\subseteq\Ksch{G}$ with Levi component $\Ksch{L}$ and some cuspidal coefficient system $C$ for $\Ksch{L}$ (see Definition~\ref{definition: cuspidal}). We write $\catqA\Ksch{G}$ for the set of irreducible admissible objects in $\catqC\Ksch{G}$. An \emph{admissible coefficient system for $\Ksch{G}$} is an object of $\catqC\Ksch{G}$ which is a finite direct sum of irreducible admissible coefficient systems.
\end{definition}

\begin{remark}
Observe that any cuspidal coefficient system (see Definition~\ref{definition: cuspidal}) is an irreducible admissible coefficient system; thus, $\catqA^{(0)}\Ksch{G} \subseteq \catqA\Ksch{G}$.
\end{remark}

The adjective `admissible' is surely one of the most over-used in mathematics, and our use of it here suggests a lack of imagination. However, as we shall show in Section~\ref{section: representations}, since our admissible coefficient systems form a bridge between certain admissible perverse sheaves and certain admissible representations, we have elected to use the adjective here also. Nevertheless, a word of caution is in order: admissible perverse sheaves are, by definition, irreducible, while our admissible coefficient systems are not. 

Examples of admissible coefficient systems are provided in Section~\ref{section: examples}.


\section{Enters Frobenius}\label{section: frobenius}

Let $\Kq$ be a $p$-adic field, let $\Rq$ be the ring of integers of $\Kq$ and let $\kq$ denote the residue field for $\Kq$. Let
$\Knr$ be a maximal unramified extension of $\Kq$ and let $\knr$ be the residue field for $\Knr$. As in Section~\ref{section: admissible} we note that $\Knr$ is strictly henselian and that $\knr$ is an algebraic closure of $\kq$, which is a finite field. Fix an isomorphism $\Gal{\Knr/\Kq} \iso \Gal{\knr/\kq}$; this determines a `lift' $\frob_\Kq \in \Gal{\Knr/\Kq}$ of the geometric Frobenius $\Frob_\kq \in \Gal{\knr/\kq}$.

Let $\Ksch{G}_{\Kq}$ be a connected reductive linear algebraic group over $\Kq$ such that
\begin{equation}
\Ksch{G} \ceq \Ksch{G}_\Kq \times_{\Spec{\Kq}} \Spec{\Knr}
\end{equation}
satisfies the conditions of Section~\ref{section: fundamental notions}. Then $\Ksch{G}$ is a connected reductive split linear algebraic group over $\Knr$ which is defined over $\Kq$.  In particular, since $\Ksch{G}$ is split it follows that $\Ksch{G}_\Kq$ splits over an unramified extension of $\Kq$. Let $\GKq$ denote the group of $\Kq$-rational points on $\Ksch{G}$. 

Since $\Ksch{G}$ is defined over $\Kq$, the Galois group $\Gal{\Knr/\Kq}$ acts on $I(\Ksch{G},\Knr)$, and
\begin{equation}
I(\Ksch{G},\Kq) = I(\Ksch{G},\Knr)^{\Gal{\Knr/\Kq}},
\end{equation}
where $I(\Ksch{G},\Kq)$ is the Bruhat-Tits building for $\GKq$.


\subsection{Frobenius-stable coefficient systems}

Let $i$ be any facet of $I(\Ksch{G},\Knr)$ and let $\Rsch{G}_i$ be the associated $\Rnr$-scheme, as in Section~\ref{subsection: integral models}, and let $\Rsch{G}_{\frob(i)}$ be the $\Rnr$-scheme associated to $\frob(i)$. The geometric Frobenius $\frob_\Kq : \Knr \to \Knr$ defines a isomorphism $\frob_i : \Rsch{G}_i \to \Rsch{G}_{\frob(i)}$ of $\Knr$-schemes which restricts to an isomorphism of special fibres, and factors through the reduction quotient maps to an isomorphism of reductive quotients $\Frob_i : \quo{G}_i \to \quo{G}_{\frob(i)}$.
Let $\Frob_i^* : D^b_c(\quo{G}_{\frob(i)};\EE) \to D^b_c(\quo{G}_i;\EE)$ be the derived functor.

\begin{lemma}\label{lemma: frobenius and restriction}
Let $i$, $j$ and $k$ be facets of $I(\Ksch{G},\Knr)$ with $i \leq j\leq k$. Then there is an isomorphism of functors
\[
\res^{\frob}_{i\leq j} : \res_{i\leq j}\ \Frob_{i}^* \to \Frob_{j}^* \res_{\frob(i)\leq \frob(j)}
\]
such that 
\[
	\xymatrix{
	\ar[rrrr]^{\res_{i\leq j\leq k} \Frob_i^*}
\res_{j\leq k} \res_{i\leq j} \Frob_i^* 
	\ar[d]^{\res_{j\leq k} \res^{\frob}_{i\leq j}}
&&&&
\res_{i\leq k} \Frob_i^* 
	\ar[dd]_{\res^{\frob}_{i\leq k}}
\\
\res_{j\leq k} \Frob_j^*\ \res_{\frob(i)\leq \frob(j)}
	\ar[d]^{\res^{\frob}_{j\leq k}\ \res_{\frob(i)\leq \frob(j)}}
&&&&
\\
\Frob_k^*\ \res_{\frob(j)\leq \frob(k)} \res_{\frob(i)\leq \frob(j)}
	\ar[rrrr]_{\Frob_k^*\ \res_{\frob(i)\leq \frob(j)\leq \frob(k)}}
&&&&
\Frob_k^*\ \res_{\frob(i)\leq \frob(k)}
	}
\]
commutes.
\end{lemma}

\begin{proof}
The proof of Lemma~\ref{lemma: frobenius and restriction} follows the lines of the proofs of Lemmas~\ref{lemma: parabolic restriction on quotients} and \ref{lemma: conjugation and restriction}. In particular, the isomorphism $\res^{\frob}_{i\leq j}$ is defined by the base-change homomorphism and the natural isomorphisms resulting from the construction of the derived functors $\Frob_i^*$, $\ksch{s}_{i\leq j}^*$ and ${\ksch{r}_{i\leq j}}_!$. The proof of the lemma then follows from \cite[Expos\'e XVII, \S 5.2]{SGA4}.
\end{proof}

\begin{proposition}\label{proposition: frobenius}
There is a unique canonical functor $\frob^* :\catqD\Ksch{G} \to \catqD\Ksch{G}$ such that $(\frob^* A)_i = \Frob_i^*\ A_{\frob(i)}$ for each object $A$ of $\catqC\Ksch{G}$ and for each facet $i$ of $I(\Ksch{G},\KK)$.
\end{proposition}

\begin{proof}
The promised functor is given as follows. 
For $A\in \obj\catqD\Ksch{G}$ and $\phi\in \mor\catqD\Ksch{G}$, and for any facet $i$ and $j$ of $I(\Ksch{G},\Knr)$ such that $i\leq j$, define \[(\frob^* A)_i\ceq \Frob_i^*\ A_{\frob(i)},\] define \[(\frob^* A)_{i\leq j}\ceq \Frob_j^*\ A_{\frob(i)\leq \frob(j)} \circ \res^{\frob}_{i\leq j} A_{\frob(i)}\] and define \[(\frob^* \phi)_i\ceq \Frob_i^*\ \phi_{\frob(i)}.\]

To show that $\frob^* A$ is an object of $\catqD\Ksch{G}$ when $A$ is an object of $\catqD\Ksch{G}$, we must turn again to Definition~\ref{definition: catqD}(obj). From the definition of $\res^{\frob}_{i\leq j}$ in the proof of Lemma~\ref{lemma: frobenius and restriction} (and the fact that $\res_{i\leq j} = \id$) it is clear that $\res^{\frob}_{i\leq i} A_{\frob(i)} = \id_{\Frob(i)^* A_{\frob(i)}}$. Thus,
\begin{eqnarray*}
(\frob^*A)_{i\leq i}
    &=& \Frob_{i}^*\ A_{\frob(i)\leq \frob(i)} \circ \res^{\frob}_{i\leq i} \Frob_i^* A_{\frob(i)}\\
    &=& \Frob_{i}^*\ \id_{A_{\frob(i)}} \circ \id_{\Frob_i^* A_{\frob(i)}} \\
    &=& \id_{\Frob_{i}^* A_{\frob(i)}} \circ \id_{\Frob_i^* A_{\frob(i)}}\\
    &=& \id_{(\frob^*A)_i}.
\end{eqnarray*}
Having shown that $\frob^*A$ satisfies the first condition set out in Definition~\ref{definition: catqD}(obj), we now turn to the second part of Definition~\ref{definition: catqD}(obj). Suppose $i$, $j$ and $k$ are facets of $I(\Ksch{G},\KK)$ with $i \leq j \leq k$. We must now show that the following diagramme commutes.
\begin{equation}\label{equation: f.0}	
\xymatrix{
	\ar[d]_{\res_{i\leq j\leq k} (\frob^* A)_{i\leq j}} 
\res_{j\leq k} \res_{i\leq j} (\frob^* A)_{i}
	\ar[rrr]^{\res_{j\leq k} (\frob^* A)_{i\leq j} } 
&&& 
	\ar[d]^{(\frob^* A)_{j\leq k}}
\res_{j\leq k} (\frob^* A)_{j}
\\
\res_{i\leq k} (\frob^* A)_{i}	
	\ar[rrr]_{(\frob^* A)_{i\leq k}}
&&&  
(\frob^* A)_{k}
	}
\end{equation}
To that end, consider the diagramme below, in which the outer square is Diagramme~\ref{equation: f.0}. (To save space we have written $\res^{\frob(i)}_{\frob(j)}$ for $\res_{\frob(i)\leq \frob(j)}$ and $A^{\frob(i)}_{\frob(j)}$ for $A_{\frob(i)\leq \frob(j)}$, etc...)
\[
	\xymatrix{
	\ar[ddd]
\res^{j}_{k} \res^{i}_{j} \Frob_i^* A_{\frob(i)}
	\ar[rrr] 
	\ar[dr]_{4.}
&&& 
	\ar[ddd]
\res^{j}_{k} \frob_j^* A_{\frob(j)}
	\ar[dl]_{2.} 
\\
&	
\Frob_k^* \res^{\frob(j)}_{\frob(k)} \res^{\frob(i)}_{\frob(j)} A_{\frob(i)}
	\ar[r]
	\ar[d]
& 
\Frob_k^* \res^{\frob(j)}_{\frob(k)} A_{\frob(j)} 
	\ar[d]
& 
\\
&	
\Frob_k^* \res^{\frob(i)}_{\frob(k)} A_{\frob(i)}
	\ar[r] 
& 
\Frob_k^* A_{\frob(k)}
& 
\\
	\ar[ur]_{3.} 
\res^{i}_{k} \Frob_i^* A_{\frob(i)}	
	\ar[rrr]
&&&  
\Frob_k^* A_{\frob(k)}
	\ar[ul]^{1.} 
}
\]
To show that $\frob^* A$ satisfies the second condition appearing in Definition~\ref{definition: catqD}(obj) we must show that the outer square commutes in Diagramme~\ref{equation: f.0}. The inner square is the result of applying the functor $\Frob_k^*\ $ to the commuting square appearing in Definition~\ref{definition: catqD}(obj) applied to $A_{\frob(i)}$, and is therefore commutative. The arrow marked $1.$ is the identity. The arrow marked $2.$ is $\res^{\frob}_{j\leq k} A_{\frob(j)}$, so the right-hand square commutes by virtue of the definition of $(\frob^* A)_{j\leq k}$; likewise, the arrow marked $3.$ is $\res^{\frob}_{i \leq k} A_{\frob(i)}$, so the bottom square commutes by virtue of the definition of $(\frob^* A)_{j\leq k}$. The arrow marked $4.$ is $\res^{\frob}_{j \leq k} \res_{\frob(i)\leq \frob(j)} A_{\frob(i)} \circ \res_{j\leq k} \res^{\frob}_{i \leq j} A_{\frob(i)}$ and the top and left-hand squares commute by Lemma~\ref{lemma: frobenius}. This concludes the demonstration that $\frob^* A$ is an object in $\catqD\Ksch{G}$.

Suppose $\phi : A \to B$ is a morphism in $\catqD\Ksch{G}$. In order to show that $\frob^* \phi$ is a morphism in $\catqD\Ksch{G}$ we must show that the following diagramme commutes.
\begin{equation}\label{equation: f.1}	
	\xymatrix{
	\ar[rrr]^{\res_{i\leq j} (\frob^*\phi)_{i}} 
\res_{i\leq j} (\frob^*A)_{i} 
	\ar[d]_{(\frob^*A)_{i\leq j}} 
&&&
\res_{i\leq j} (\frob^*B)_{i} 
	\ar[d]^{(\frob^*B)_{i\leq j}} 
\\
(\frob^*\ A)_{j} 
	\ar[rrr]_{ (\frob^*\phi)_{j}} 
&&& 
(\frob^*B)_{j}
	}
\end{equation}
Consider the following diagramme, in which the outer square is Diagramme~\ref{equation: f.1}. (To save space we have written $\res^{\frob(i)}_{\frob(j)}$ for $\res_{\frob(i)\leq \frob(j)}$ and $A^{\frob(i)}_{\frob(j)}$ for $A_{\frob(i)\leq \frob(j)}$, etc... , as above.)
\[
	\xymatrix{
	\ar[rrr]
\res_{i\leq j} \Frob_i^*A_{\frob(i)} 
	\ar[dd]
	\ar[dr]^{3.}
&&&
\res_{i\leq j} \Frob_i^*B_{\frob(i)} 
	\ar[dd] 
	\ar[dl]^{4.}
\\
&
\Frob_j^*\res_{\frob(i)\leq \frob(j)} A_{\frob(i)}
	\ar[r]^{0.}
	\ar[dl]^{2.}
&
\Frob_j^*\res_{\frob(i)\leq \frob(j)} B_{\frob(i)}
	\ar[dr]^{1.}
&
\\
\Frob_j^*A_{\frob(j)} 
	\ar[rrr]
&&& 
\Frob_j^*B_{\frob(j)}
	}
\]
The arrow marked $0.$ is $\Frob_j^*\res_{\frob(i)\leq \frob(j)} \phi_{\frob(i)}$, the arrow marked $1.$ is $\Frob_j^* B_{\frob(i)\leq \frob(j)}$ and the arrow marked $2.$ is $\Frob_j^*A_{\frob(i)\leq \frob(j)}$; thus, the bottom square is the result of applying the functor $\Frob_j^*$ to the relevant form of the commuting square appearing in Definition~\ref{definition: catqD}(mor), and therefore commutes since $\phi$ is a morphism in $\catqD\Ksch{G}$. The arrow marked $3.$ is $\res^{\frob}_{i\leq j} A_{\frob(i)}$ and the arrow marked $3.$ is $\res^{\frob}_{i\leq j} B_{\frob(i)}$, so the upper square commutes because $\res^{\frob}_{i\leq j} $ is a natural transformation. The left-hand triangle commutes by virtue of the definition of $(\frob^* A)_{i\leq j}$ and likewise the right-hand triangle commutes by virtue of the definition of $(\frob^* B)_{i\leq j}$. Therefore, the outer square commutes. This concludes the demonstration that $\frob^*\phi$ is a morphism in $\catqD\Ksch{G}$.
\end{proof}

\begin{definition}\label{definition: frobenius stable}
An object $A$ of the category $\catqD\Ksch{G}$ is \emph{frobenius-stable} if there is an isomorphism $\alpha : \frob^*A \to A$ in category $\catqD\Ksch{G}$.
\end{definition}

\begin{lemma}\label{lemma: frobenius}
Let $i$ and $j$ be facets of $I(\Ksch{G},\Knr)$ with $i \leq j$. Let $g$ be an element of $\Ksch{G}(\Knr)$. Then
\[
\ind_{i\leq j}\ \Frob_{j}^* \iso
\Frob_{i}^*\
    \ind_{\frob(i)\leq \frob(j)},
\]
and
\[
\quo{m}(g)_i^*\ \Frob_{gi}^* \iso \Frob_{i}^*\
\quo{m}(\frob(g))_{\frob(i)}^*.
\]
\end{lemma}


\subsection{Parabolic restriction and frobenius}\label{subsection: parabolic restriction and frobenius}

\begin{definition}\label{definition: twisted-Levi}
Let $\Ksch{L}_\Kq\subseteq \Ksch{G}_\Kq$ be a subgroup. Here, $\Ksch{L}_\Kq$, $\Ksch{G}_\Kq$ and the inclusion are all over $\Kq$. We say that $\Ksch{L}_\Kq\subseteq \Ksch{G}_\Kq$ is an \emph{unramified twisted-Levi subgroup} if there is a finite extension $\Kq' : \Kq$ contained in $\Knr$ such that $\Ksch{L}_{\Kq'}\subseteq \Ksch{G}_{\Kq'}$ is a Levi subgroup. Note that this implies that there is a parabolic subgroup $\Ksch{P}$ defined over $\Kq'$ such that $\Ksch{L}_{\Kq'}$ is the maximal reductive quotient of $\Ksch{P}_{\Kq'}$.
\end{definition}

Note that an unramified twisted-Levi subgroup is not, in general, a Levi subgroup; rather, an unramified twisted-Levi subgroup is a \emph{form} of a Levi subgroup.

\begin{proposition}\label{proposition: parabolic restriction and frobenius}
Let $\Ksch{P}$ be a parabolic subgroup of $\Ksch{G}$ with Levi component $\Ksch{L}$. Suppose $\Ksch{L}$ is defined over $\Kq$ and $\Ksch{L}_\Kq \subseteq \Ksch{G}_\Kq$ is an unramified twisted-Levi subgroup. Let $A$ be an object of $\catqD\Ksch{G}$. If $A$ is frobenius-stable then $\res^{\Ksch{G}}_{\Ksch{P}}A$ is also frobenius-stable.
\end{proposition}

\begin{proof}
Let $\alpha : \frob^* A \to A$ be an isomorphism in category $\catqD\Ksch{G}$. Recall the definition of $\res^{\Ksch{G}}_{\Ksch{P}}\phi$ (see Proposition~\ref{proposition: parabolic restriction}). Lemma~\ref{lemma: frobenius} shows that $\res^{\Ksch{G}}_{\Ksch{P}}\alpha$ is an isomorphism and therefore that $\res^{\Ksch{G}}_{\Ksch{P}} A$ is frobenius-stable.
\end{proof}


\subsection{Cuspidal coefficient systems revisited}\label{subsection: cuspidal and frobenius}

In this section we briefly revisit Section~\ref{subsection: cuspidal} and adapt Definition~\ref{definition: cuspidal} and Theorem~\ref{theorem: cuspidal} to the present context.

Let $i_0$ be a vertex of $I(\Ksch{G},\Kq)$ and let $F$ be a cuspidal character sheaf for $\quo{G}_{i_0}$ equipped with an
isomorphism $\varphi : \Frob_{i_0}^* F \to F$ in $D^b_c(\quo{G}_{i_0};\EE)$. Then there is an object $B$ and isomorphism $\beta : \frob_{\Ksch{G}}^* B \to B$ in $\catqC\Ksch{G}$ such that $B_{i_0} = F$, $\beta_{i_0} = \varphi$.
The proof follows the lines of the proof of Proposition~\ref{proposition: cuspidal}.
For any vertex $i$ of $I(\ksch{G},\Kq)$, cuspidal character sheaf $F$ on $\quo{G}_i$ and isomorphism
$\varphi : \Frob_{i}^* F \to F$, we will write $\cind^{\Ksch{G}}_{\quo{G}_i} F$ and $\cind^{\Ksch{G}}_{\quo{G}_i}\varphi$ for the object and morphism in $\catqC\Ksch{G}$ promised above.
Suppose $C$ is a cuspidal coefficient system for $\Ksch{G}$ and frobenius-stable. Then there is a vertex
$i_0$ of $I(\Ksch{G},\Kq)$ and a frobenius-stable cuspidal character sheaf $F$ for $\quo{G}_{i_0}$ such that
$C \iso \cind^{\Ksch{G}}_{\quo{G}_{i_0}} F$.
The proof follows the lines of the proof of Theorem~\ref{theorem: cuspidal}.


\subsection{Parabolic induction and frobenius}\label{subsection: parabolic induction and frobenius}

\begin{proposition}\label{proposition: induction and frobenius}
Let $\Ksch{P}$ be a parabolic subgroup of $\Ksch{G}$ with Levi component $\Ksch{L}$. Suppose $\Ksch{L}$ is defined over $\Kq$ and $\Ksch{L}_\Kq\subseteq \Ksch{G}_\Kq$ is an unramified-twisted Levi subgroup. If $B$ is a cuspidal coefficient system for  $\Ksch{L}$ then
\begin{equation}
\ind^\Ksch{G}_\Ksch{P}\ \frob_\Ksch{L}^* B \iso \frob_\Ksch{G}^*\
\ind^\Ksch{G}_\Ksch{P} B.
\end{equation}
\end{proposition}

\begin{proof}
Let $i$ be a facet of $I(\Ksch{G},\Knr)$ and consider $(\ind^\Ksch{G}_\Ksch{P}\ \frob_\Ksch{L}^* B)_i$. Using
Proposition~\ref{proposition: parabolic induction} and Proposition~\ref{proposition: frobenius} we have
\begin{eqnarray*}
(\ind^\Ksch{G}_\Ksch{P}\ \frob_\Ksch{L}^* B)_i
 &=& \sum_{g\in d_i(\Ksch{G},\Ksch{P},\Knr)} \quo{m}(g^{-1})^*\
    \ind^{\quo{G}_{ig}}_{\quo{P}_{i^g_\Ksch{P}}}\
    (\frob_\Ksch{L}^* B)_{i^g_\Ksch{P}}\\
 &=& \sum_{g\in d_i(\Ksch{G},\Ksch{P},\Knr)} \quo{m}(g^{-1})^*\
    \ind^{\quo{G}_{ig}}_{\quo{P}_{i^g_\Ksch{P}}}\
    (\frob_\Ksch{L}^* B)_{\frob(i^g_\Ksch{P})}\\
 &=& \sum_{g\in d_i(\Ksch{G},\Ksch{P},\Knr)} \quo{m}(g^{-1})^*\
    \ind^{\quo{G}_{ig}}_{\quo{P}_{i^g_\Ksch{P}}}\
    \Frob_{i^g_\Ksch{P}}^*\ B_{\frob(i)^{\frob(g)}_\Ksch{P}},
\end{eqnarray*}
since $\Ksch{L}$ is defined over $\Kq$. Now,using Lemma~\ref{lemma: frobenius} we have
\begin{eqnarray*}
 &&\hskip-20pt
     \sum_{g\in d_i(\Ksch{G},\Ksch{P},\Knr)} \quo{m}(g^{-1})^*\
    \ind^{\quo{G}_{ig}}_{\quo{P}_{i^g_\Ksch{P}}}\
    \Frob_{i^g_\Ksch{P}}^*\ B_{\frob(i)^{\frob(g)}_\Ksch{P}}\\
 &=& \sum_{g\in d_i(\Ksch{G},\Ksch{P},\Knr)} \quo{m}(g^{-1})^*\
    \Frob_{i^g_\Ksch{P}}^*\
    \ind^{\quo{G}_{\frob(ig)}}_{\quo{P}_{{\frob(i)}^{\frob(g)}_\Ksch{P}}}\
     B_{\frob(i)^{\frob(g)}_\Ksch{P}}\\
 &=& \sum_{g\in d_i(\Ksch{G},\Ksch{P},\Knr)} \Frob_{i}^*\
 \quo{m}(\frob(g)^{-1})^*\
     \ind^{\quo{G}_{\frob(ig)}}_{\quo{P}_{{\frob(i)}^{\frob(g)}_\Ksch{P}}}\
     B_{\frob(i)^{\frob(g)}_\Ksch{P}}\\
 &=& \Frob_{i}^* \sum_{g\in d_i(\Ksch{G},\Ksch{P},\Knr)}
 \quo{m}(\frob(g)^{-1})^*\
     \ind^{\quo{G}_{\frob(ig)}}_{\quo{P}_{{\frob(i)}^{\frob(g)}_\Ksch{P}}}\
     B_{\frob(i)^{\frob(g)}_\Ksch{P}}.
\end{eqnarray*}
If $\Rsch{P}_{ig} \subset \Rsch{G}_{ig}$ is projective then $\Rsch{P}_{\frob(i)\frob(g)} \subset
\Rsch{G}_{\frob(i)\frob(g)}$ is projective, since $\Ksch{L}$ is a $\Kq$-scheme, so $\frob(g)$ represents an element of
$D_{\frob(i)}(\Ksch{G},\Ksch{P},\Knr)$ (see Definition~\ref{definition: Di}). Accordingly, there is some $h\in \Rsch{G}_{\frob(i)}(\Rnr)$ and $l\in \Ksch{L}(\Knr)$ such that $\frob(g) = h g' l$, for a unique $g'\in d_{\frob(i)}(\Ksch{G},\Ksch{P},\Knr)$ (see Definition~\ref{definition: Di} again). Thus,
\begin{eqnarray*}
 && \hskip-20pt
    \Frob_{i}^* \sum_{g\in d_i(\Ksch{G},\Ksch{P},\Knr)}\
    \quo{m}(\frob(g)^{-1})^*\
     \ind^{\quo{G}_{\frob(ig)}}_{\quo{P}_{{\frob(i)}^{\frob(g)}_\Ksch{P}}}\
     B_{\frob(i)^{\frob(g)}_\Ksch{P}}\\
 &=& \Frob_{i}^* \sum_{{g'}\in d_{\frob(i)}(\Ksch{G},\Ksch{P},\Knr)}  \quo{m}((h g' l)^{-1})^*\
     \ind^{\quo{G}_{\frob(i)h g' l}}_{\quo{P}_{{\frob(i)}^{h g' l}_\Ksch{P}}}\
     B_{\frob(i)^{h g' l}_\Ksch{P}}\\
 &=& \Frob_{i}^* \sum_{{g'}\in d_{\frob(i)}(\Ksch{G},\Ksch{P},\Knr)}
    \quo{m}(h^{-1})^* \quo{m}(g'^{-1})^* \quo{m}(l^{-1})^*
     \ind^{\quo{G}_{\frob(i)xl}}_{\quo{P}_{{\frob(i)}^{{g'}}_\Ksch{P} l}}\
     B_{\frob(i)^{{g'}}_\Ksch{P} l}\\
 &=& \Frob_{i}^* \sum_{{g'}\in d_{\frob(i)}(\Ksch{G},\Ksch{P},\Knr)}
    \quo{m}(h^{-1})^* \quo{m}(g'^{-1})^*
     \ind^{\quo{G}_{\frob(i){g'}}}_{\quo{P}_{{\frob(i)}^{{g'}}_\Ksch{P}}}
     \quo{m}(l^{-1})^*\ B_{\frob(i)^{{g'}}_\Ksch{P} l}\\
 &=& \Frob_{i}^* \sum_{{g'}\in d_{\frob(i)}(\Ksch{G},\Ksch{P},\Knr)}
    \quo{m}(h^{-1})^* \quo{m}(g'^{-1})^*
     \ind^{\quo{G}_{\frob(i){g'}}}_{\quo{P}_{{\frob(i)}^{{g'}}_\Ksch{P}}}\
     \conj{l}{B}_{\frob(i)^{{g'}}_\Ksch{P}}.
\end{eqnarray*}
Since $B$ is weakly-equivariant (\cf Definition~\ref{definition: weakly-equivariant}) we have
\begin{eqnarray*}
 &&\hskip -20pt
    \Frob_{i}^* \sum_{{g'}\in d_{\frob(i)}(\Ksch{G},\Ksch{P},\Knr)}
    \quo{m}(h^{-1})^* \quo{m}(g'^{-1})^*\
     \ind^{\quo{G}_{\frob(i){g'}}}_{\quo{L}_{{\frob(i)}^{{g'}}_\Ksch{P}}}\
     \conj{l}{B}_{\frob(i)^{{g'}}_\Ksch{P}}\\
 &\iso&\Frob_{i}^* \sum_{{g'}\in d_{\frob(i)}(\Ksch{G},\Ksch{P},\Knr)}
    \quo{m}(h^{-1})^* \quo{m}(g'^{-1})^*\
     \ind^{\quo{G}_{\frob(i)g'}}_{\quo{P}_{{\frob(i)}^{{g'}}_\Ksch{P}}}\
     {B}_{\frob(i)^{{g'}}_\Ksch{P}}.
\end{eqnarray*}
Since $B$ is cuspidal (\cf Definition~\ref{definition: cuspidal}) we have
\begin{eqnarray*}
 &&\hskip -20pt
    \Frob_{i}^* \sum_{{g'}\in d_{\frob(i)}(\Ksch{G},\Ksch{P},\Knr)}
    \quo{m}(h^{-1})^* \quo{m}(g'^{-1})^*\
     \ind^{\quo{G}_{\frob(i){g'}}}_{\quo{P}_{{\frob(i)}^{{g'}}_\Ksch{P}}}\
     {B}_{\frob(i)^{{g'}}_\Ksch{P}}\\
  &\iso&\Frob_{i}^* \sum_{{g'}\in d_{\frob(i)}(\Ksch{G},\Ksch{P},\Knr)}
    \quo{m}(g'^{-1})^*\
     \ind^{\quo{G}_{\frob(i)}^{{g'}}}_{\quo{P}_{{\frob(i)}^{{g'}}_\Ksch{P}}}\
     {B}_{\frob(i)^{{g'}}_\Ksch{P}}\\
 &=&\Frob_{i}^* (\ind^\Ksch{G}_\Ksch{P} B)_{\frob(i)}\\
 &=& (\frob^* \ind^\Ksch{G}_\Ksch{P} B)_i.
\end{eqnarray*}
The isomorphism of Proposition~\ref{proposition: induction and frobenius} is now found by arguments similar to those employed in the proof of Proposition~\ref{proposition: parabolic induction}.
\end{proof}

\begin{corollary}\label{corollary: induction and frobenius}
Let $\Ksch{P}$ be a parabolic subgroup of $\Ksch{G}$ with Levi component $\Ksch{L}$. Suppose $\Ksch{L}$ is defined over $\Kq$ and $\Ksch{L}_\Kq\subseteq \Ksch{G}_\Kq$ is an unramified twisted-Levi subgroup. Let $B$ be an weakly-equivariant object equipped with an isomorphism $\beta : \frob_\Ksch{L}^* B \to B$ in $\catqC\Ksch{L}$. Then
$\ind^\Ksch{G}_\Ksch{L} \beta$ defines an isomorphism $\frob_\Ksch{G}^*\ind^{\Ksch{G}}_{\Ksch{P}}B \to
\ind^{\Ksch{G}}_{\Ksch{P}}B$ in $\catqC\Ksch{G}$.
\end{corollary}

\begin{proof}
Let $A = \ind^\Ksch{G}_\Ksch{P}B$ and let $\alpha = \ind^\Ksch{G}_\Ksch{P} \beta$. Since $\beta : \frob_\Ksch{L}^* B
\to B$ is an isomorphism in $\catqC\Ksch{L}$, and since parabolic induction is a functor, it follows that $\alpha$ is an
isomorphism. However, the domain of $\alpha$ is $\ind^\Ksch{G}_\Ksch{P} \frob_\Ksch{L}^* B$, rather than
$\frob_\Ksch{G}^* \ind^\Ksch{G}_\Ksch{P} B$. Corollary~\ref{corollary: induction and frobenius} now follows
directly from Proposition~\ref{proposition: induction and frobenius}.
\end{proof}


\section{Supercuspidal depth-zero representations}\label{section: representations}

Let the field extensions $\Knr:\Kq$ and $\knr:\kq$ be as in Section~\ref{section: frobenius}; likewise, let $\Ksch{G}_\Kq$ be a connected unramified linear algebraic group over $\Kq$ and let $\Ksch{G}$ denoted the group scheme over $\Knr$ obtained by extension of scalars, as in Section~\ref{section: frobenius}.

Suppose now that $i$ is fixed by the Galois action on the building. Using the principle of \'etale descent we see that there
is a smooth group scheme ${\Rsch{G}_i}_{/\Rq}$ over $\Rq$ such that: ${\Rsch{G}_i}_{/\Rq}$ is an integral model of $\Ksch{G}_\Kq$ equipped with a $\Rq$-rational structure, compatible with the isomorphism of generic fibres, and such that
${\Rsch{G}_i}_{/\Rq}(\Rq) = \Ksch{G}(\Kq)_i$ (\cf \cite[10.10]{Lan1}). The special fibre $\sfib{G}_{i/\kq}$ of
$\Rsch{G}_{i/\Rq}$ is a linear algebraic group and it defines a $\kq$-rational structure for $\sfib{G}_i$. Moreover, there is a
maximal reductive quotient ${\nu_i}_{/\kq} : \sfib{G}_{i/\kq} \to \quo{G}_{i/\kq}$ (\cf Section~\ref{subsection: integral models}) and it defines a $\kq$-rational structure for $\quo{G}_i$. Let $\rho_{i/\Kq}: \Rsch{G}_i(\Rq) \to
\quo{G}_i(\kq)$ be the canonical quotient (\cf Section~\ref{subsection: integral models}).


\subsection{Characteristic functions}

Let $A$ be a frobenius-stable coefficient system for $\Ksch{G}$ and let $\alpha : \frob_\Ksch{G}^* A \to A$ be an isomorphism in $\catqC\Ksch{G}$. For each facet $i$ of $I(\Ksch{G},\Knr)$, $\alpha$ defines an isomorphism $\alpha_i : \frob_{i}^* A_{\frob(i)} \to A_i$ in category $D^b_c(\quo{G}_i;\EE)$. If $i$ is actually a facet of $I(\Ksch{G},\Kq)$
(so $\frob(i) = i$) then $A_i \in  D^b_c(\quo{G}_i;\EE)$ is frobenius-stable in the usual sense; in that case, let $\chf{A_i,\alpha_i} : \quo{G}_i(\kq) \to \bar\QQ_\ell$ be the characteristic function associated to the pair $(A_i, \alpha_i)$ (\cf \cite[8.4]{CS}, for example).

\begin{proposition}\label{proposition: conjugation and characteristic functions} 
Let $A$ be a weakly-equivariant coefficient system for $\Ksch{G}$ and let $\alpha : \frob^* A \to A$ be an isomorphism in
$\catqC\Ksch{G}$. If $i$ is a facet of $I(\Ksch{G},\Kq)$ then
\begin{equation}
\forall x\in \Rsch{G}_{i}(\Rq), \qquad
\chf{A_{gi},\alpha_{gi}}(\rho_{gi}(gxg^{-1})) =
\chf{A_i,\alpha_i}(\rho_i(x)),
\end{equation}
for all $g\in \Ksch{G}(\Kq)$.
\end{proposition}

Let $\Ksch{L}_\Kq \subseteq \Ksch{G}_\Kq$ be a twisted-Levi subgroup; so $\Ksch{L} \ceq \Ksch{L}_\Kq \times_{\Spec{\Kq}}
\Spec{\Knr}$ is a Levi subgroup of $\Ksch{G} \ceq \Ksch{G}_\Kq \times_{\Spec{\Kq}} \Spec{\Knr}$. Let $\Ksch{P}$ be a parabolic subgroup of $\Ksch{G}$ with Levi component $\Ksch{L}$. Let $B$ be a frobenius-stable weakly-equivariant coefficient system for $\Ksch{L}$ equipped with $\beta : \frob^*B \to B$ and let $A = \ind^\Ksch{G}_\Ksch{P} B$ equipped with $\alpha : \frob^* A \to A$ as in Corollary~\ref{corollary: induction and frobenius}. For any facet $i$ of $I(\Ksch{G},\Knr)$ and $g\in d_i(\Ksch{G},\Ksch{P}, \Knr)$ let $A_i(g)$ denote the summand of $A_i$ in $D^b_c(\quo{G}_i;\EE)$ given by
\begin{equation}
A_i(g) \ceq \quo{m}(g^{-1})_{i}^*\
\ind^{\quo{G}_{ig}}_{\quo{L}_{i^g_L}} \ B_{i^g_L}.
\end{equation}

We will need the following result concerning characteristic functions of induced objects in the proof of Theorem~\ref{theorem:
supercuspidal models}.

\begin{proposition}\label{proposition: diq}
With notation as above, suppose $i$ is a facet of $I(\Ksch{G},\Kq)$ and let $d_i(\Ksch{G},\Ksch{L},\Kq)$ denote the
set of $g\in d_i(\Ksch{G},\Ksch{P},\Knr)$ (\cf Definition~\ref{definition: Di}) such that $\Rsch{G}_{ig}$ is defined over $\Rq$. Then
\begin{equation}
\forall x\in \quo{G}_i(\kq), \qquad
    \chf{A_i,\alpha_i}(x)
    = \sum_{g\in d_i(\Ksch{G},\Ksch{P},\Kq)}
    \chf{A_i(g),\alpha_i(g)}(x),
\end{equation}
where the right-hand side is trivial if $d_i(\Ksch{G},\Ksch{L},\Kq)$ is empty.
\end{proposition}

\begin{proof}
If $g$ is an element of $d_i(\Ksch{G},\Ksch{P},\Kq)$ then $A_i(g)$ is itself frobenius-stable; more precisely, the restriction of
$\alpha_i$ to $\Frob_{i}^* A_i(g)$, which we denote $\alpha_i(g)$, is an isomorphism onto $A_i(g)$. It follows that
\begin{equation}
\sum_{g\in d_i(\Ksch{G},\Ksch{P},\Kq)} A_i(g)
\end{equation} 
is a frobenius-stable summand of $A_i$. For any $g\in d_i(\Ksch{G},\Ksch{P},\Knr)$ we have $\frob(g) = h g' l$ for a
unique $g'\in d_i(\Ksch{G},\Ksch{L},\Knr)$ with $h\in \Rsch{G}_i(\Rnr)$ and $l\in \Ksch{L}(\Knr)$ (\cf proof of
Proposition~\ref{proposition: induction and frobenius}). Suppose $g\in d_i(\Ksch{G},\Ksch{P},\Knr)$ and $g\not\in
d_i(\Ksch{G},\Ksch{P},\Kq)$. Then $\quo{G}_{ig}(\kq) = \emptyset$. There are two cases to consider: either $g'=g$ or $g'\ne g$. In the first case, $A_i(g)$ is a frobenius-stable summand of $A_i$ with $\chf{A_i(g),\alpha_i(g)} =0$, with $\alpha_i(g)$ as above; it follows that the sum of such $A_i(g)$ is a frobenius-stable summand of $A_i$ with trivial characteristic function. The sum of the objects $A_i(g)$ with $g$ in the second case is also a frobenius-stable summand of $A_i$ with trivial characteristic function. However, in this case the summands $A_i(g)$ are not themselves frobenius-stable. In summary, we have
\begin{equation}
\forall x\in \quo{G}_i(\kq), \qquad
    \chf{A_i,\alpha_i}(x)
    = \sum_{g\in d_i(\Ksch{G},\Ksch{P},\Kq)}
    \chf{A_i(g),\alpha_i(g)}(x),
\end{equation}
where the right-hand side is trivial if $d_i(\Ksch{G},\Ksch{P},\Kq)$ is empty, as desired.
\end{proof}


\subsection{Models for representations}\label{subsection: models for representations}

Let $\pi: \GKq \to \End_{\bar\QQ_\ell}(V)$ be an admissible representation. For each facet $i$ of the building
$I(\Ksch{G},\Kq)$, let $V_i$ denote the $\bar\QQ_\ell$-vector space consisting of all $v\in V$ for which $\pi(h)v = v$ for each
$h\in \Rsch{G}_i(\Rnr)$ such that $\rho_{i/\Kq}(h) = 1$. We let
\begin{equation}
{\pi}_{i}: \Rsch{G}_{i}(\Rq) \to \End_{\bar\QQ_\ell}(V_{i}),
\end{equation}
denote the \emph{compact restriction} of $\pi$ at $i$; that is to say, we let $\pi_{i}$ denote the representation of the group
$\Rsch{G}_i(\Rq)$ on $V_i$ defined by $\pi_i(h)v = \pi(h)v$ for $h\in \Rsch{G}_i(\Rq)$ and $v\in V_i$. We note that, for each
$x\in i$, the group $\Ksch{G}(\Kq)_i$ equals $\Ksch{G}(\Kq)_{x,0}$, where the latter is defined in \cite{Moy-Prasad}. Moreover, the group of all $h\in \Rsch{G}_i(\Kq)$ such that $\rho_{i/\Kq}(h) =1$ is exactly the group $\Ksch{G}(\Kq)_{x,0^+}$, for $x\in i$ (\cf \cite{Moy-Prasad} also). By \cite[Propn 1.1]{V1}, $\Ksch{G}(\Kq)_{x,0^+}$ equals $U^{(0)}_x$, where the latter is defined in \cite{SS}. Thus the set of $\bar\QQ_\ell$-vector spaces
\begin{equation}
\gamma_0(V) = \left\{  V_{i} \tq i \text{ facet of }
I(\Ksch{G},\Kq) \right\},
\end{equation}
equipped with inclusions $V_j \hookrightarrow V_i$ for $i \leq j$, forms a $\GKq$-equivariant coefficient system of
$\bar\QQ_\ell$-vector spaces, in the sense of \cite{SS}.

Now suppose $\pi$ is a depth-zero supercuspidal representation. Then, for each facet $i$ of $I(\Ksch{G},\Kq)$, the compact restriction $\pi_i$ factors through $\rho_{i/\Kq}$ to a representation
\begin{equation}
\quo{\pi}_{i} : \quo{G}_{i}(\kq)
\to \End_{\bar\QQ_\ell}(V_{i});
\end{equation}
that is, $\quo{\pi}_i = \pi_i \circ \rho_{i/\Kq}$. If $\pi$ is non-trivial, then, by the definition of depth-zero
representations, the coefficient system $\gamma_0(V)$ is non-trivial; in fact there is some vertex $i_0$ of
$I(\Ksch{G},\Kq)$ such that $V_{i_0} \ne 0$.

\begin{definition}\label{definition: model}
Let $\catqK\Ksch{G}$ denote the subgroup of the Grothendieck group for $\catqD\Ksch{G}$ generated by admissible coefficient systems for $\Ksch{G}$ (\cf Definition~\ref{definition: admissible}). Let $\pi : \Ksch{G}(\Kq) \to \Aut_{\bar\QQ_\ell}(V)$ be a depth-zero admissible representation. A \emph{model for $\pi$} is an element $\sum_n a_n [A^n]$ of $\catqK\Ksch{G}\otimes_\ZZ \bar\QQ_\ell$ where each $A^n$ is equipped with an isomorphism $\alpha^n : \frob^* A^n \to A^n$ (so $A^n$ is frobenius-stable)
such that
\begin{equation}
\forall x\in \quo{G}_i(\kq), \qquad \sum_n a_n
\chf{A^n_i,\alpha^n_ii}(x) = \trace\quo{\pi}_i(x),
\end{equation}
for each facet $i$ of $I(\Ksch{G},\Kq)$.
\end{definition}

\begin{theorem}\label{theorem: supercuspidal models}
Supercuspidal depth-zero representations admit models; that is, for each supercuspidal depth-zero representation $\pi$ of $\Ksch{G}(\Kq)$ there is a some $\sum_n a_n [A^n]\in \catqK\Ksch{G}\otimes_\ZZ \bar\QQ_\ell$ and isomorphisms $\alpha^n : \frob^* A^n \to A^n$  such that
\begin{equation}
\forall x\in \quo{G}_i(\kq), \qquad \sum_n a_n
\chf{A^n_i,\alpha^n_i}(x) = \trace\quo{\pi}_i(x),
\end{equation}
for each facet $i$ of $I(\Ksch{G},\Kq)$.
\end{theorem}

\begin{proof}
Let $\pi : \GKq \to \Aut_{\bar\QQ_\ell}(V)$ be an irreducible supercuspidal representation. Since $\pi$ has depth-zero, there is
a vertex $i_0$ of $I(\Ksch{G},\Kq)$ such that $(\Rsch{G}_{i_0}(\Rq), \pi_{i_0})$ is a type of $\pi$. Moreover,
since $\pi$ is irreducible, we have $\pi_i =0$ unless $i$ is contained in the $\GKq$-orbit of $i_0$ in $I(\Ksch{G},\Kq)$. By
\cite[25.1]{CS} and \cite[(10.4.5)]{CS}, we may write 
\begin{equation}\label{eqn: ms 1}
  \trace\quo{\pi}_{i_0}
  =
  \sum_n a_{n}\,
  \chf{\ind_{\ksch{M}_n}^{\quo{G}_{i_0}} F_n,\ind_{\ksch{M}_n}^{\quo{G}_{i_0}} \varphi_n},
\end{equation}
where the sum is taken over cuspidal pairs $(\ksch{M}_n,F_n)$ for $\quo{G}_{i_0}$  (as defined in \cite[2.4]{L0}) such that $\ksch{M}_n$ is defined over $\kq$ and where $F_n$ is equipped with a fixed isomorphism $\varphi_n : \Frob_{\ksch{M}_n}^* F_n \to F_n$. Since $\quo{\pi}_{i_0}$ is cuspidal, the scalars $a_n$ are zero except when $\ksch{M}_n$ is anisotropic over $\kq$ (\cf \cite[(15.2.1)]{CS}).
We define $E_n = \ind_{\ksch{M}_n}^{\quo{G}_{i_0}} F_n$ and $\varepsilon_n = \ind_{\ksch{M}_n}^{\quo{G}_{i_0}}\varphi_n$.

Let $\Ksch{T}^n$ be a maximal $\Knr$-split torus in $\Ksch{G}$, defined over $\Kq$, such that the associated apartment contains $i_0$, and such that the image of $\Ksch{T}^n(\Knr)\cap \Rsch{G}_{i_0}(\Rq)$ in $\ksch{M}_n$ is the group of
$\kq$-rational points of an elliptic torus of $\ksch{M}_n$ (the existence of $\Ksch{T}^n$ is guaranteed by \cite[end of the proof of prop. 5.1.10]{BT2}). The torus $\Ksch{T}^n$ is elliptic. Thus, the centralizer $\Ksch{L}^n$ of $\Ksch{T}^n$ in $\Ksch{G}$ is an elliptic Levi subgroup defined over $\Kq$ such that the image of $I(\Ksch{L}^n,\Kq)$ under $I(\Ksch{L}^n,\Knr) \to I(\Ksch{G},\Knr)$ is $\{ i_0\}$, and $\quo{L}^{n}_{i_0} = \ksch{M}_n$. Now $\Ksch{L}^n_\Kq$ is an elliptic unramified twisted-Levi subgroup of $\Ksch{G}_\Kq$. Let $\Ksch{P}^n$ be a parabolic subgroup of $\Ksch{G}$ with Levi component $\Ksch{L}^n = \Ksch{L}^n_\Kq \times_\Spec{\Kq} \Spec{\Knr}$.

Using Section~\ref{subsection: cuspidal and frobenius}, set
\begin{equation}
B^n = \cind^{\Ksch{L}^n}_{\quo{L}^n_{i_0}} F_n
 \quad
\text{and}
 \quad
\beta^n = \cind^{\Ksch{L}^n}_{\quo{L}^n_{i_0}}\varphi_n.
\end{equation}
Then
\begin{equation}\label{eqn: sc 3}
B^n_{i_0} = F_n
 \qquad
\text{and}
 \qquad
\beta^n_{i_0} = \varphi_n.
\end{equation}
Define
\begin{equation}
A^n = \ind^\Ksch{G}_\Ksch{P} B^n.
\end{equation}
Let
\begin{equation}
\alpha^n : \frob^* A^n \to A^n
\end{equation}
be the isomorphism given by Corollary~\ref{corollary: induction and frobenius}. In order to prove Theorem~\ref{theorem:
supercuspidal models} we now show that \begin{equation}\label{eqn: sc main}
\forall x\in \quo{G}_i(\kq),\qquad
    \sum_n a_n \chf{A^n_i,\alpha^n_i}(x) = \trace \quo{\pi}_i(x),
\end{equation}
for each facet $i$ of $I(\Ksch{G},\Kq)$.

Consider the case $i = i_0$. Notice that since $i_0$ is contained in the image of the map from $I(\Ksch{L}^n,\Knr)$ to
$I(\Ksch{G},\Knr)$, it follows that the schematic closure of $\Ksch{L}^n$ in $\Rsch{G}_{i_0}$ is $\Rsch{L}^n_{i_0}$. If $g\in
d_{i_0}(\Ksch{G},\Ksch{P}^n,\Kq)$ then $\Rsch{G}_{i_0g}$ is defined over $\Kq$, so the schematic closure of $\Ksch{L}^n$ in $\Rsch{G}_{i_0 g}$ is a parahoric subgroup of $\Ksch{L}^n$ defined over $\Kq$. Since $\Ksch{L}^n_\Kq$ is elliptic there is exactly one such parahoric subgroup, namely $\Ksch{L}^n_{i_0}$, and $d_{i_0}(\Ksch{G},\Ksch{P}^n,\Kq)$ is a singleton. Thus, $g$ and $1$ represent the same double coset in $D_{i_0}(\Ksch{G},\Ksch{L}^n,\Knr)$ (\cf Definition~\ref{definition: Di}), so $g = hl$ for some $h \in \Rsch{G}_i(\Rnr)$ and $l\in \Ksch{L}^n(\Knr)$. Using Lemma~\ref{lemma: frobenius} we have
\begin{eqnarray*}
 A^n_{i_0 g} &=& \quo{m}(g^{-1})_{i_0}^*\
 \ind^{\quo{G}_{i_0 g}}_{(\quo{P}^n)_{(i_0)^g_{\Ksch{P}^n}}}\
 B_{(i_0)^g_{\Ksch{P}^n}}\\
 &=&\quo{m}((hl)^{-1})_{i_0}^*\
 \ind^{\quo{G}_{i_0 hl}}_{(\quo{P}^n)_{(i_0)^{(hl)}_{\Ksch{P}^n}}}\
 B_{(i_0)^{hl}_{\Ksch{P}^n}}\\
 &=&\quo{m}(h^{-1})^* \quo{m}(l^{-1})_{i_0}^*\
 \ind^{\quo{G}_{i_0l}}_{(\quo{P}^n)_{i_0}^l}\
 B_{i_0 l}\\
 &=& \quo{m}(h^{-1})^*\
 \ind^{\quo{G}_{i_0}}_{\quo{P}^n_{i_0}}
 \quo{m}(l^{-1})_{i_0}^*\ B_{i_0 l}\\
 &=& \quo{m}(h^{-1})^*
 \ind^{\quo{G}_{i_0}}_{\quo{P}^n_{i_0}}\
 \conj{l}{B}_{i_0}.
\end{eqnarray*}
Moreover, since $B$ is a weakly-equivariant object of $\catqC\Ksch{L}^n$ and $h \in \Rsch{G}_{i_0}(\Rnr)$ we have canonical isomorphisms
\[
 \quo{m}(h^{-1})^*\
 \ind^{\quo{G}_{i_0}}_{\quo{P}^n_{i_0}}
 \conj{l}{B}_{i_0}
 \iso \quo{m}(h)_{i_0}^*\ \ind^{\quo{G}_{i_0}}_{\quo{P}^n_{i_0}}\
 B^n_{i_0}\\
 \iso \ind^{\quo{G}_{i_0}}_{\quo{P}^n_{i_0}}\ B^n_{i_0}.
\]
Finally, since
\begin{equation}
\quo{P}^n_{i_0} = \quo{L}^n_{i_0} = \ksch{M}_n,
\end{equation}
we have
\[
   \ind^{\quo{G}_{i_0}}_{\quo{L}^n_{i_0}}\ B^n_{i_0}
 =\ind^{\quo{G}_{i_0}}_{\ksch{M}_n}\ F_n =E_n,
\]
by Equation~\ref{eqn: sc 3}. Since the isomorphisms just used are exactly those appearing in the definition of $\alpha^n_{i_0}$ (\cf Corollary~\ref{corollary: induction and frobenius}), together with Proposition~\ref{proposition: diq} we have
\begin{equation}
\forall x\in \quo{G}_{i_0}(\kq),\qquad
\chf{A^n_{i_0},\alpha^n_{i_0}}(x) =
\chf{E_n,\varepsilon_n}(x).
\end{equation}
By Equation~\ref{eqn: ms 1} we now have
\begin{equation}\label{eqn: sc 7}
\forall x\in \quo{G}_{i_0}(\kq), \qquad \sum_{n} a_n
\chf{A^n_{i_0},\alpha^n_{i_0}}(x) = \trace\quo{\pi}_{i_0}(x).
\end{equation}
This verifies Equation~\ref{eqn: sc main} in this case. Proposition~\ref{proposition: conjugation and characteristic
functions} extends Equation~\ref{eqn: sc 7}, \mutmut, to all $i$ in the $\Ksch{G}(\Kq)$-orbit of $i_0$.

If $i$ is a facet of $I(\Ksch{G},\Kq)$ which does not lie in the $\Ksch{G}(\Kq)$-orbit of $i_0$ then $d_i(\Ksch{G},\Ksch{P}^n,\Kq)$ is empty, again since $\Ksch{L}^n_\Kq$ is an elliptic twisted-Levi
subgroup of $\Ksch{G}_\Kq$, so Proposition~\ref{proposition: diq}
gives
\begin{equation}\label{eqn: sc 8}
\forall x \in \quo{G}_i(\kq), \qquad \chf{A^n_i,\alpha^n_i}(x) = 0.
\end{equation}
Gathering Equations~\ref{eqn: ms 1}, \ref{eqn: sc 7} and \ref{eqn: sc 8} gives Equation~\ref{eqn: sc main} in this case also.

Finally, consider the $\bar\QQ_\ell$-vector space formed by taking the tensor product of the subgroup $\catK_0(\quo{G}_{i_0})$ of the Grothendieck group of perverse sheaves on $\quo{G}_{i_0}$ generated by character sheaves of $\quo{G}_{i_0}$ with $\bar\QQ_\ell$ (\cf \cite[14.10]{CS}). Then $\sum_{n} a_n [E_n]$ (summation as in Equation~\ref{eqn: ms 1}) is an element of $\catK_0(\quo{G}_{i_0}) \otimes_{\ZZ} \bar\QQ_\ell$. Now, Equation~\ref{eqn: sc main} shows that $\sum_n a_n [A^n]$, equipped with the isomorphisms $\alpha^n$, is a model for $\pi$, which completes the proof of Theorem~\ref{theorem: supercuspidal models}.
\end{proof}


\subsection{Distributions associated to admissible coefficient systems}\label{subsection: characters}

\begin{definition}\label{definition: character}
Let $\Ksch{G}(\Kq)_\er$ denote the set of elliptic regular elements of $\Ksch{G}(\Kq)$. Let $A$ be an admissible coefficient system for $\Ksch{G}$ which is frobenius-stable with respect to $\alpha$ and weakly-equivariant. Let $\chf{A,\alpha} : \Ksch{G}(\Kq)_\er \to \bar\QQ_\ell$ be the function defined by
\begin{equation}
\chf{A,\alpha}(g) \ceq \sum_{\{i \in I(\Ksch{G},\Kq)\tq g\in
\Rsch{G}_i(\Rq)\}} (-1)^{\dim i} \chf{A_i,\alpha_i}(\rho_i(g)).
\end{equation}
We will refer to $\chf{A,\alpha}$ as the \emph{character of $A$}.
\end{definition}

\begin{theorem}\label{theorem: character formula}
Let $\pi : \Ksch{G}(\Kq) \to \Aut_{\bar\QQ_\ell}(V)$ be a depth-zero admissible representation and let $\Theta_\pi$ be the character of $\pi$ in the sense of Harish-Chandra. Let $A$ be a model for $\pi$. Then
\begin{equation}
\Theta_\pi(g) = \chf{A,\alpha}(g),
\end{equation}
for all $g\in \Ksch{G}(\Kq)_\er$.
\end{theorem}

\begin{proof}
Since $\pi$ is a depth-zero representation, there is a vertex $i_0$ of $I(\Ksch{G},\Kq)$ such that the $\Ksch{G}(\Kq)$-algebra $V$ is generated by $V_{i_0}$. Thus, $V$ is an object in category $\operatorname{Alg}^0 \Ksch{G}(\Knr)$ (\cf \cite{SS}). Now \cite[III.4.10]{SS} and \cite[III.4.16]{SS} extend to any quasi-split reductive linear algebraic group  (\cf
\cite[Cor.3.33]{Co}) so
\begin{equation}
\Theta_\pi(g) = \sum_{i\in I(\Ksch{G},\Kq),\ g \in
\Rsch{G}_i(\Rq)} (-1)^{\dim i} \trace(g, V_i).
\end{equation}
For $g\in \Rsch{G}_i(\Rq)$, the trace $\trace(g,V_i)$ is the character of $\quo{\pi}_i$. By Definition~\ref{definition: model}
we have $\trace(g,V_i) = \chf{A_i,\alpha_i}(\rho_i(g))$, which concludes the demonstration.
\end{proof}


\section{Examples}\label{section: examples}

Let $\Kq$ and $\Knr$ be as in Section~\ref{section: representations}.


\subsection{SL(2)}\label{example: SL(2)}

Let $\Ksch{G}_\Kq = \SL(2)_\Kq$. We now describe all frobenius-stable cuspidal coefficient systems for each
unramified twisted-Levi subgroup of $\Ksch{G}_\Kq$. Table~\ref{table: models for SL2} records models for all irreducible depth-zero supercuspidal representations of $\SL(2,\Kq)$; in this section we describe the terms used in Table~\ref{table: models for SL2}.

\begin{table}[htdp]
\caption{Models for depth-zero supercuspidal representations of $\SL(2,\Kq)$}
\begin{center}
\begin{tabular}{|c|c|}
\hline 
{Representation} & {Model}\\
\hline
$\pi_\theta$ & $-B_\theta$\\
\hline
$\pi_{+}$ & $-\frac{1}{2} B_\sgn + \frac{1}{2} \left( C^{+} - C^{-} \right)$\\
\hline
$\pi_{-}$ & $-\frac{1}{2} B_\sgn - \frac{1}{2} \left( C^{+} - C^{-} \right)$\\
\hline
$\pi_{+}'$ & $-\frac{1}{2} B_\sgn' + \frac{1}{2} \left( D^{+} - D^{-} \right)$\\
\hline
$\pi_{-}'$ & $-\frac{1}{2} B_\sgn' - \frac{1}{2} \left( D^{+} - D^{-} \right)$\\
\hline
$\pi_\theta'$ & $-B_\theta'$\\
\hline
\end{tabular}
\end{center}
\label{table: models for SL2}
\end{table}

Up to conjugation over $\Kq$, there are four unramified twisted-Levi subgroups of $\Ksch{G}_\Kq$: the split torus $\Ksch{S}_\Kq$, two unramified elliptic tori $\Ksch{T}_\Kq$ and $\Ksch{T}_\Kq'$, and the group $\Ksch{G}_\Kq$ itself.

\begin{itemize}
\item
Consider the split torus
\[
\Ksch{S}_\Kq
    =
    \left\{ \begin{pmatrix} s_1 & 0 \\ 0 & s_2 \end{pmatrix}
    \ \Big\vert\
    s_1 s_2 =1\right\}.
\]
Set $\Ksch{S} = \Ksch{S}_\Kq \otimes \Knr$. All frobenius-stable cuspidal coefficient systems for $\Ksch{S}$ take the form $\cind^{\Ksch{S}}_{\quo{S}_{(0)}} \fais{L}_\theta[1]$, where $\fais{L}_\theta$ indicates a Kummer system on
$\quo{S}_{(0)} = \GL(1)_\knr$ equipped with an isomorphism $\Frob^* \fais{L}_\theta \to \fais{L}_\theta$ such that the
characteristic function of $\fais{L}_\theta$ equals the character $\theta$ of $\GL(1,\kq)$. Define 
\[
A(\theta) \ceq \ind^{\Ksch{G}}_{\Ksch{S}} \cind^{\Ksch{S}}_{\quo{S}_{(0)}} \fais{L}_\theta[1].
\]
Using Definition~\ref{definition: cuspidal} and Proposition~\ref{proposition: parabolic induction} we see that
\begin{eqnarray*}
A(\theta)_{(0)}
    &=& \ind^{\quo{G}_{(0)}}_{\quo{T}_{(0)}} \fais{L}_\theta[1]\\
A(\theta)_{(01)} &=&
     \fais{L}_\theta[1]
    \oplus
    \quo{m}(s_{(1)})_{(01')}^*
    \fais{L}_\theta[1]\\
A(\theta)_{(1)}
    &=& \ind^{\quo{G}_{(1)}}_{\quo{T}_{(0)}} \fais{L}_\theta[1]
\end{eqnarray*}
Observe that $\quo{G}_{(01)} = \quo{G}_{(01')} = \quo{T}_{(0)}$. The object $A(\theta)$ is simple in $\catqC\Ksch{G}$ unless $\theta^2 =1$; thus, in that case, $A(\theta)$ is an admissible coefficient system.

\item
All frobenius-stable cuspidal coefficient systems for
\[
\Ksch{T}_{\Kq}
    =
    \left\{ \begin{pmatrix} x & y \\ \varepsilon y & x \end{pmatrix}
    \ \Big\vert\
    x^2-\varepsilon y^2 =1 \right\},
\]
where $\varepsilon\in \Rq$ is a fixed quadratic residue and a unit in $\Rq$, take the form $\cind^{\Ksch{T}}_{\quo{T}_{(0)}} \fais{L}_\theta[1]$, where $\theta$ is a character of the finite group $\quo{T}_{(0)}(\kq)$. Define
\[
B(\theta) \ceq \ind^{\Ksch{G}}_{\Ksch{T}} \cind^{\Ksch{T}}_{\quo{T}_{(0)}} \fais{L}_\theta[1].
\]
The object $B(\theta)$, equipped with an isomorphism $\beta(\theta) : \frob^* B(\theta) \to B(\theta)$ induced from $\phi_\theta: \Frob^*\fais{L}_\theta \to \fais{L}_\theta$, provides our first non-trivial example of the phenomenon described in Proposition~\ref{proposition: diq}. From Definition~\ref{definition: cuspidal} and the definition of parabolic induction on objects (\cf Proposition~\ref{proposition: parabolic induction}) we have
\begin{eqnarray*}
B(\theta)_{(0)}
    &=& \ind^{\quo{G}_{(0)}}_{\quo{T}_{(0)}} \fais{L}_{\theta}[1] \\
B(\theta)_{(01)}
    &=& \quo{m}(z^{-1})_{(01)}^* \fais{L}_{\theta}[1] \oplus
    \quo{m}(z^{-1}w)_{(01)}^* \fais{L}_{\theta}[1] \\
B(\theta)_{(1)}
    &=& \quo{m}(z^{-1})_{(1)}^* \ind^{\quo{G}_{(1)z}}_{\quo{T}_{(1)z}}
    \fais{L}_{\theta}[1],
\end{eqnarray*}
where $\Ksch{T}_\Kq = \Ksch{S}^{z}_\Kq$ and where $w$ is a representative for the non-trivial element in the Weyl group for $\Ksch{S}$. (Recall that $\Ksch{S}$ is split.) Here, $z$ is an element of $\Ksch{G}(\Kq')$, where $\Kq'$ is a quadratic extension of $\Kq$, {\em not} an element of $\Ksch{G}(\Kq)$, which represents a class in $H^1(\Ksch{S}_\Kq,\Gal{\Knr/\Kq})$; in other words, we view $\Ksch{T}_\Kq$ as a twist of $\Ksch{S}$ and represent that twist by conjugation by $z^{-1}$. Now, from Proposition~\ref{proposition: frobenius} and the definition of parabolic induction on maps (\cf Proposition~\ref{proposition: parabolic induction}) we have
\begin{eqnarray*}
\beta(\theta)_{(0)} &=& \ind^{\quo{G}_{(0)}}_{\quo{T}_{(0)}} \phi_{\theta}[1] \\
\beta(\theta)_{(01)} &=& \quo{m}(z^{-1})_{(01)}^* \phi_{\theta}[1] \oplus
    \quo{m}(z^{-1}w)_{(01)}^* \phi_{\theta}[1] \\
\beta(\theta)_{(1)} &=& \quo{m}(z^{-1})_{(1)}^* \ind^{\quo{G}_{(0)}}_{\quo{T}_{(0)}}
    \phi_{\theta}[1].
\end{eqnarray*}
Therefore
\begin{eqnarray*}
    \chf{B(\theta)_{(0)},\beta(\theta)_{(0)}} &=& \chf{\ind^{\quo{G}_{(0)}}_{\quo{T}_{(0)}}
    \fais{L}_{\theta}[1],\ind^{\quo{G}_{(0)}}_{\quo{T}_{(0)}}
    \phi_{\theta}[1]}\\
    \chf{B(\theta)_{(01)},\beta(\theta)_{(01)}} &=& 0 \\
    \chf{B(\theta)_{(1)},\beta(\theta)_{(1)}} &=& 0.
\end{eqnarray*}

\item
Likewise, all frobenius-stable cuspidal coefficient systems for
\[
\Ksch{T}_{\Kq}'
    =
    \left\{ \begin{pmatrix} x & \varpi^{-1} y \\ \varepsilon \varpi y & x \end{pmatrix}
    \ \Big\vert\
    x^2-\varepsilon y^2 =1 \right\},
\]
where $\varpi\in \Rq$ is a fixed uniformizer for $\Kq$, take the form $\cind^{\Ksch{T}}_{\quo{T}_{(1)}'} \fais{L}_\theta[1]$, where $\theta$ is a character of the finite group $\quo{T}_{(1)}'(\kq)$.  Define
\[
B(\theta)' \ceq \ind^{\Ksch{G}}_{\Ksch{T}'} \cind^{\Ksch{T}'}_{\quo{T}_{(1)}'} \fais{L}_\theta[1].
\]
The object $B(\theta)'$, equipped with an isomorphism $\beta(\theta)' : \frob^* B(\theta) \to B(\theta)'$ induced from $\phi_\theta: \Frob^*\fais{L}_\theta \to \fais{L}_\theta$ provides another example of the phenomenon described in Proposition~\ref{proposition: diq}. From Definition~\ref{definition: cuspidal} and the definition of parabolic induction on objects (\cf Proposition~\ref{proposition: parabolic induction}) we have
\begin{eqnarray*}
B(\theta)_{(0)}'
    &=& \quo{m}(v^{-1})_{(0)}^* \ind^{\quo{G}_{(0)v}}_{\quo{T}_{(0)v}'}
    \fais{L}_{\theta}[1]\\
B(\theta)_{(01)}'
    &=& \quo{m}(v^{-1})_{(01)}^* \fais{L}_{\theta}[1] \oplus
    \quo{m}(v^{-1}w)_{(01)}^* \fais{L}_{\theta}[1] \\
B(\theta)_{(1)}'
    &=& \ind^{\quo{G}_{(1)}}_{\quo{T}_{(1)}'} \fais{L}_{\theta}[1] ,
\end{eqnarray*}
where $\Ksch{T}_\Kq' = \Ksch{S}^{v}_\Kq$ and where $w$ is a representative for the non-trivial element in the Weyl group for $\Ksch{S}$ as above. Note that $v$ is an element of $\Ksch{G}(\Knr_2)$, {\em not} an element of $\Ksch{G}(\Kq)$, which represents a class in $H^1(\Ksch{S}_\Kq,\Gal{\Knr/\Kq})$. As above, from Proposition~\ref{proposition: frobenius} and the definition of parabolic induction on maps (\cf Proposition~\ref{proposition: parabolic induction}) it follows that
\begin{eqnarray*}
\beta(\theta)_{(0)}' &=& \quo{m}(v^{-1})_{(0)}^* \ind^{\quo{G}_{(1)}}_{\quo{T}_{(1)}'}
    \phi_{\theta}[1]\\
\beta(\theta)_{(01)}' &=& \quo{m}(v^{-1})_{(01)}^* \phi_{\theta}[1] \oplus
    \quo{m}(v^{-1}w)_{(01)}^* \phi_{\theta}[1] \\
\beta(\theta)_{(1)}' &=& \ind^{\quo{G}_{(1)}}_{\quo{T}_{(1)}'} \phi_{\theta}[1] .
\end{eqnarray*}
Therefore,
\begin{eqnarray*}
    \chf{B(\theta)_{(0)}',\beta(\theta)_{(0)}'} &=& 0\\
    \chf{B(\theta)_{(01)}',\beta(\theta)_{(01)}'} &=& 0 \\
    \chf{B(\theta)_{(1)}',\beta(\theta)_{(1)}'} &=& \chf{\ind^{\quo{G}_{(1)}}_{\quo{T}_{(1)}'}
    \fais{L}_{\theta}[1],\ind^{\quo{G}_{(1)}}_{\quo{T}_{(1)}'}
    \phi_{\theta}[1]}.
\end{eqnarray*}

\item
Finally, there are exactly four frobenius-stable cuspidal coefficient systems for
\[
\Ksch{G}_\Kq
    =
    \left\{ \begin{pmatrix} a & b \\ c & d \end{pmatrix}
    \ \Big\vert\
    ad-bc =1\right\};
\]
they are
\begin{eqnarray*}
        C^{\pm} &\ceq& \cind_{\quo{G}_{(0)}}^{\Ksch{G}} F_\pm\\
    D^{\pm} &\ceq& \cind_{\quo{G}_{(1)}}^{\Ksch{G}} F_\pm.
\end{eqnarray*}
The following facts are simple consequences of Definition~\ref{definition: cuspidal}:
\begin{eqnarray*}
C^{\pm}_{(0)} &=& F_\pm\\ 
C^{\pm}_{(01)} &=& 0\\
C^{\pm}_{(1)} &=& 0,
\end{eqnarray*}
and
\begin{eqnarray*}
D^{\pm}_{(0)} &=& 0\\
D^{\pm}_{(01)} &=& 0 \\
D^{\pm}_{(1)} &=& F_\pm.
\end{eqnarray*}
\end{itemize}

We now describe all depth-zero supercuspidal representations of $\SL(2,\Kq)$. Let $\pi_\theta$ denote the
representation obtained by compact induction from the Deligne-Lusztig representation
$-R^{\quo{G}_{(0)}}_{\quo{T}_{(0)}}(\theta)$ of $\quo{G}_{(0)}(\kq)$, where $\theta$ is a character of
$\quo{T}_{(0)}(\kq)$ in general position; likewise, let $\pi_\theta'$ denote the representation obtained by compact
induction from the Deligne-Lusztig representation $-R^{\quo{G}_{(1)}}_{\quo{T}_{(1)}'}(\theta)$ of
$\quo{G}_{(1)}(\kq)$, where $\theta$ is a character of $\quo{T}_{(1)'}(\kq)$ in general position. Next, let $\chi_0^\pm$
denote the two cuspidal representations appearing in the Lusztig series for $R^{\SL(2)}_{\SU(1)}(\sgn)$, where $\sgn$ is the sign character of $\SU(1,\kq)$; let $\pi_\pm$ denote the representation obtained by compact induction from $\chi_0^\pm$ on
$\quo{G}_{(0)}(\kq)$ and likewise let $\pi_\pm'$ denote the representation obtained by compact induction from $\chi_0^\pm$ on $\quo{G}_{(1)}(\kq)$. 

Although $-B(\theta)$ is a model for the depth-zero supercuspidal representation $\pi_\theta$ (see Table~\ref{table: models for SL2}) this example illustrates how {\em not} to find models for representations. Consider the coefficient system $\underbar{V}$ for $\pi_\theta$ in the sense of \cite{SS}. Even though $V_{(01)}$  and $V_{(1)}$ are both zero, it does not follow that $B(\theta)_{(01)}$ and $B(\theta)_{(1)}$ are zero. As one sees from the proof of Theorem~\ref{theorem: supercuspidal models}, to make a model for a representation it is necessary to identify a
type for that representation. Observe that $(\Rsch{G}_{(0)}(\Rq), -R^{\quo{G}_{(0)}}_{\quo{T}_{(0)}}(\theta))$ is a type for
$\pi_\theta$ and that
\begin{equation}
-\trace R^{\quo{G}_{(0)}}_{\quo{T}_{(0)}}(\theta) =
\chf{\ind^{\quo{G}_{(0)}}_{\quo{T}_{(0)},
    \fais{L}_{\theta}[1],}{\ind^{\quo{G}_{(0)}}_{\quo{T}_{(0)}}
    \phi_{\theta}[1]}}.
\end{equation}

We remark that $\{ \pi_{+}, \pi_{+}', \pi_{-}, \pi_{-}' \}$ is an L-packet and that $\{ \pi(\theta), \pi'(\theta) \}$ is also an L-packet when $\theta$ is in general position.


\subsection{Sp(4)}\label{example: Sp(4)}

Let $\Ksch{G}_\Kq = \Sp(4)_\Kq$. We now describe all frobenius-stable cuspidal coefficient systems for each
unramified twisted-Levi subgroup of $\Ksch{G}_\Kq$, as in Example~\ref{example: SL(2)}. Using \cite{S1} and \cite{S2} we produce models for all supercuspidal depth-zero representations $\Sp(4,\Kq)$ and present the results in Tables~\ref{table: models for Sp4.1}, \ref{table: models for Sp4.0} and \ref{table: models for Sp4.2}, in which $\zeta\in \bar\QQ_\ell$ is a primitive fourth root of unity such that $\zeta^2$ equals the quadratic residue of $-1$ in $\kq$. We now explain the other terms appearing in these tables.

\begin{table}[htdp]
\caption{Models for representations induced from $\Rsch{G}_{(1)}(\Rq)$}
\begin{center}
\begin{tabular}{|c|c|}
\hline
{Representation} & {Model}\\
\hline
 $\pi_0$ &
    $A_4(\theta_1,\theta_2)$\\
\hline
 $\pi_1$ &
    $\begin{array}{c}
    \frac{1}{2}A_4(\theta,\sgn)\\
    -\frac{1}{2}\left(B_2^+(\theta)-B_2^-(\theta)\right)
    \end{array}$\\
\hline
 $\pi_2$ &
    $\begin{array}{c}
    \frac{1}{2}A_4(\theta,\sgn)\\
    +\frac{1}{2}\left(B_2^+(\theta)-B_2^-(\theta)\right)
    \end{array}$\\
\hline
 $\pi_3$ &
    $\begin{array}{c}
    \frac{1}{4}A_4(\sgn,\sgn)\\
    -\frac{1}{4}\left(B_2^+(\sgn)-B_2^-(\sgn)\right)\\
    -\frac{1}{4}\left(B_1^+(\sgn)-B_1^-(\sgn)\right)\\
    + \frac{1}{4}\left(D^{+\, +} - D^{+\, -} - D^{-\, +} + D^{-\, -}\right)
    \end{array}$ \\
\hline
 $\pi_4$ &
    $\begin{array}{c}
    \frac{1}{4}A_4(\sgn,\sgn)\\
    +\frac{1}{4}\left(B_2^+(\sgn)-B_2^-(\sgn)\right)\\
    -\frac{1}{4}\left(B_1^+(\sgn)-B_1^-(\sgn)\right)\\
    - \frac{1}{4}\left(D^{+\, +} - D^{+\, -} - D^{-\, +} + D^{-\, -}\right)
    \end{array}$ \\
\hline
 $\pi_5$ &
    $\begin{array}{c}
    \frac{1}{4}A_4(\sgn,\sgn)\\
    +\frac{1}{4}\left(B_2^+(\sgn)-B_2^-(\sgn)\right)\\
    +\frac{1}{4}\left(B_1^+(\sgn)-B_1^-(\sgn)\right)\\
    +\frac{1}{4}\left(C_{+\,+}-C_{+\,-}-C_{-\,+}+C_{-\,-}\right)
    \end{array}$\\
\hline
\end{tabular}
\end{center}
\label{table: models for Sp4.1}
\end{table}

\begin{table}[htdp]
\caption{Models for representations induced from $\Rsch{G}_{(0)}(\Rq)$}
\begin{center}
\begin{tabular}{|c|c|}
\hline
{Representation} & {Model}\\
\hline
 $\pi_6$ &
    $A_7(\theta_1,\theta_2)$\\
\hline
 $\pi_7$ &
    $A_3(\theta_1,\theta_2)$\\
\hline
 $\pi_8$ &
    $\frac{1}{2}A_3(\theta,\sgn) -\frac{1}{2}\zeta\left( B_1^+(\theta)-B_1^-(\theta)\right)$\\
\hline
 $\pi_9$ &
    $\frac{1}{2}A_3(\theta,\sgn) +\frac{1}{2}\zeta\left( B_1^+(\theta)-B_1^-(\theta)\right)$\\
\hline
 $\pi_{10}$ & $-\frac{1}{4}A_7(1) -\frac{1}{4}A_3(1) -\zeta\frac{1}{2}C$\\
\hline
\end{tabular}
\end{center}
\label{table: models for Sp4.0}
\end{table}

\begin{table}[htdp]
\caption{Models for representations induced from $\Rsch{G}_{(2)}(\Rq)$}
\begin{center}
\begin{tabular}{|c|c|}
\hline
{Representation} & {Model}\\
\hline
 $\pi'_6$ &
    $A_8(\theta_1,\theta_2)$\\
\hline
 $\pi'_7$ &
    $A_5(\theta_1,\theta_2)$\\
\hline
 $\pi'_8$ &
    $\frac{1}{2}A_5(\theta,\sgn)
    -\frac{1}{2}\zeta\left( B_2^+(\theta)'-B_2^-(\theta)'\right)$\\
\hline
 $\pi'_9$ &
    $\frac{1}{2}A_5(\theta,\sgn)
    +\frac{1}{2}\zeta\left( B_2^+(\theta)'-B_2^-(\theta)'\right)$\\
\hline
 $\pi'_{10}$ &
    $-\frac{1}{4}A_8(1)
    -\frac{1}{4}A_5(1)
    -\zeta\frac{1}{2}C'$\\
\hline
\end{tabular}
\end{center}
\label{table: models for Sp4.2}
\end{table}

To begin, we fix a representation of $\Sp(4)$. Let $\theta : \GL(4)\to \GL(4)$ be the involution defined by $\theta(g) = J^{-1}\ ^tg^{-1} J$, where
\begin{equation} 
J=
\begin{pmatrix}
        0 & 0 & 0 & 1\\
        0 & 0 & -1 & 0\\
        0 & 1 & 0 & 0\\
        -1 & 0 & 0 & 0
    \end{pmatrix},
\end{equation}
and let $\Ksch{G} = \Sp(4)$ be the subvariety fixed by $\theta$. 

We now consider cuspidal unramified twisted-Levi subgroups of $\Ksch{G}_\Kq$. Up to conjugacy, the group $\Ksch{G}_\Kq$ contains thirteen cuspidal unramified twisted-Levi subgroups: nine inner forms of the torus $\GL(1)_\Kq \times \GL(1)_\Kq$; six inner forms of the Levi subgroup $\GL(1)_\Kq \times \SL(2)_\Kq$; and $\Ksch{G}_\Kq$ itself. The Levi subgroup $\GL(2)_\Kq$ is not cuspidal since it does not admit any cuspidal coefficient systems.

\begin{itemize}
\item
We begin with tori in $\Ksch{G}$. Consider the split torus
\[
\Ksch{T}_\Kq = \left\{ \begin{pmatrix}
  t_1 & 0 & 0 & 0 \\
  0 & t_2 & 0 & 0 \\
  0 & 0 & t_3 & 0 \\
  0 & 0 & 0 & t_4
\end{pmatrix} 
\Big\vert 
\begin{array}{c}t_1t_4 =1\\ t_2t_3 =1\end{array}\right\}.
\]
Let $\Ksch{T} = \Ksch{T}_\Kq \otimes \Knr$. Each frobenius-stable cuspidal coefficient system for
$\Ksch{T}_\Kq$ takes the form
\[
C_0(\theta_1,\theta_2) \ceq \cind_{\quo{T}_{(0)}}^{\Ksch{S}} \fais{L}_{\theta_1}\boxtimes\fais{L}_{\theta_2}[2],
\]
where $\theta_1$ and $\theta_2$ are characters of $\quo{T}_{(0)}(\kq)$. Here,  $\fais{L}_{\theta}$ is the
Frobenius-stable Kummer local system equipped with an isomorphism $\Frob^* \fais{L}_\theta \to \fais{L}_\theta$ such that the characteristic function $\chf{\fais{L}_\theta}{}$ of $\fais{L}$ equals $\theta$. Define
\[
A_0(\theta_1,\theta_2) \ceq \ind^{\Ksch{G}}_{\Ksch{T}} C_0(\theta_1,\theta_2).
\]
Thus, $A_0(\theta_1,\theta_2)$ is an object $\catqC\Ksch{G}(\Kq)$. If $\theta_1$ and $\theta_2$ are in general position, then $A_0(\theta_1,\theta_2)$ is simple in $\catqC\Ksch{G}$, in which case $A_0(\theta_1,\theta_2)$ is admissible (see Definition~\ref{definition: admissible}); in other words, if $\theta_1$ and $\theta_2$ are in general position then $A_0(\theta_1,\theta_2)$ is an admissible coefficient system.

\item
Next, consider torus
\[
\Ksch{T}^1_\Kq = \left\{ \begin{pmatrix}
  x & 0 & 0 & y \\
  0 & t_2 & 0 & 0 \\
  0 & 0 & t_3 & 0 \\
  \varepsilon y & 0 & 0 & x
\end{pmatrix}
\Big\vert \begin{array}{c}x^2-\varepsilon y^2 =1\\ t_2t_3 =1\end{array}\right\}.
\]
Clearly, $\Ksch{T}^1_\Kq$ splits over a quadratic unramified extension of $\Kq$.  Let $\Ksch{T}^1 = \Ksch{T}^1_\Kq \otimes \Knr$. Each frobenius-stable cuspidal coefficient system for $\Ksch{T}^1_\Kq$ takes the form
\[
C_1(\theta_1,\theta_2) \ceq \cind^{\Ksch{T}^1}_{\quo{T}^1_{(0)}}
\fais{L}_{\theta_1}\boxtimes\fais{L}_{\theta_2}[2],
\]
where $\theta_1$ is a character of $\SU(1,\kq)$ and $\theta_2$ is a character of $\GL(1,\kq)$. Define
\[
A_1(\theta_1,\theta_2) = \ind^{\Ksch{G}}_{\Ksch{T}^1} C_1(\theta_1,\theta_2).
\]
If $\theta_1$ and $\theta_2$ are in general position then $A_1(\theta_1,\theta_2)$ is an admissible coefficient system.

\item
Next, consider the torus
\[
\Ksch{T}^2_\Kq = \left\{ \begin{pmatrix}
  x & 0 & 0 & \varpi^{-1} y \\
  0 & t_2 & 0 & 0 \\
  0 & 0 & t_3 & 0 \\
  \varpi \varepsilon y & 0 & 0 & x
\end{pmatrix} \Big\vert \begin{array}{c}x^2-\varepsilon y^2 =1 \\ t_2t_3 =1\end{array}\right\}.
\]
This also splits over a quadratic unramified extension. Let $\Ksch{T}^2 = \Ksch{T}^2_\Kq \otimes \Knr$. Frobenius-stable
depth-zero cuspidal character sheaves for $\Ksch{T}^2_\Kq$ take the form
\[
C_2(\theta_1,\theta_2) \ceq \cind_{\quo{T}^2_{(1)}}^{\Ksch{T}^2}
\fais{L}_{\theta_1}\boxtimes\fais{L}_{\theta_2}[2],
\]
where $\theta_1$ is a character of $\SU(1,\kq)$ and $\theta_2$ is a character of $\GL(1,\kq)$.  Define
\[
A_2(\theta_1,\theta_2) = \ind^{\Ksch{G}}_{\Ksch{T}^2} C_2(\theta_1,\theta_2).
\]
If $\theta_1$ and $\theta_2$ are in general position then $A_2(\theta_1,\theta_2)$ is an admissible coefficient system.

\item
Next, consider the elliptic torus
\[
\Ksch{T}^3_\Kq = \left\{ \begin{pmatrix}
  x_1 & 0 & 0 &  y_1 \\
  0 & x_2 & y_2 & 0 \\
  0 & \varepsilon y_2 & x_2 & 0 \\
  \varepsilon y_1 & 0 & 0 & x_1
\end{pmatrix} \Big\vert
\begin{array}{c}
x_1^2-\varepsilon y_1^2 =1\\
x_2^2-\varepsilon y_2^2 =1
\end{array}
\right\}.
\]
Let $\Ksch{T}^3 = \Ksch{T}^3_\Kq \otimes \Knr$. The building for $\Ksch{T}^3(\Knr)$ in $\Ksch{G}(\Knr)$ is $\{ (0) \}$ and
the reductive quotient is $\quo{T}^3_{(0)} \iso \SU(1)\times \SU(1)$. Accordingly, each frobenius-stable depth-zero cuspidal character sheaf for $\Ksch{T}^3$ takes the form
\[
C_3(\theta_1,\theta_2) \ceq \cind^{\Ksch{T}^3}_{\quo{T}^3_{(0)}}
\fais{L}_{\theta_1}\boxtimes\fais{L}_{\theta_2}[2],
\]
where $\theta_1$ and $\theta_2$ are characters of $\SU(1,\kq)$. Define
\[
A_3(\theta_1,\theta_2) = \ind^{\Ksch{G}}_{\Ksch{T}^3} C_3(\theta_1,\theta_2).
\]
If $\theta_1$ and $\theta_2$ are in general position then $A_3(\theta_1,\theta_2)$ is an admissible coefficient system.

\item
Likewise, consider the elliptic torus
\[
\Ksch{T}^4_\Kq = \left\{ \begin{pmatrix}
  x_1 & 0 & 0 &  \varpi^{-1} y_1 \\
  0 & x_2 & y_2 & 0 \\
  0 & \varepsilon y_2 & x_2 & 0 \\
  \varpi \varepsilon y_1 & 0 & 0 & x_1
\end{pmatrix} \Big\vert
\begin{array}{c}
x_1^2-\varepsilon y_1^2 =1\\
x_2^2-\varepsilon y_2^2 =1
\end{array}
\right\}.
\]
The building for $\Ksch{T}^4(\Kq)$ in $\Ksch{G}(\Kq)$ is $\{ (1) \}$. Each frobenius-stable depth-zero cuspidal character sheaf for $\Ksch{T}^4_\Kq$ takes the form
\[
C_4(\theta_1,\theta_2) \ceq \cind^{\Ksch{T}^4}_{\quo{T}^4_{(1)}}
\fais{L}_{\theta_1}\boxtimes\fais{L}_{\theta_2}[2],
\]
where $\theta_1$ and $\theta_2$ are characters of $\SU(1,\kq)$. Define
\[
A_4(\theta_1,\theta_2) = \ind^{\Ksch{G}}_{\Ksch{T}^4} C_4(\theta_1,\theta_2).
\]
If $\theta_1$ and $\theta_2$ are in general position then $A_4(\theta_1,\theta_2)$ is an admissible coefficient system.

\item
We also have the elliptic torus
\[
\Ksch{T}^5_\Kq = \left\{ \begin{pmatrix}
  x_1 & 0 & 0 &  \varpi^{-1} y_1 \\
  0 & x_2 & \varpi^{-1} y_2 & 0 \\
  0 & \varpi \varepsilon y_2 & x_2 & 0 \\
  \varpi \varepsilon y_1 & 0 & 0 & x_1
\end{pmatrix} \Big\vert
\begin{array}{c}
x_1^2-\varepsilon y_1^2 =1\\
x_2^2-\varepsilon y_2^2 =1
\end{array}
\right\}.
\]
The building for $\Ksch{T}^5(\Kq)$ in $\Ksch{G}(\Kq)$ is $\{ (2) \}$ and $\quo{T}^5_{(2)/\kq} \iso \SU(1)_\kq\times \SU(1)_\kq$. Accordingly, each frobenius-stable depth-zero cuspidal character
sheaf for $\Ksch{T}^5_\Kq$ takes the form
\[
C_5(\theta_1,\theta_2) \ceq \cind^{\Ksch{T}^5}_{\quo{T}^5_{(2)}}
\fais{L}_{\theta_1}\boxtimes\fais{L}_{\theta_2}[2],
\]
where $\theta_1$ and $\theta_2$ are characters of $\SU(1,\kq)$. Define
\[
A_5(\theta_1,\theta_2) = \ind^{\Ksch{G}}_{\Ksch{T}^5} C_5(\theta_1,\theta_2).
\]
If $\theta_1$ and $\theta_2$ are in general position then $A_5(\theta_1,\theta_2)$ is an admissible coefficient system.

\item
Next, consider the torus
\[
\Ksch{T}^6_\Kq = \left\{ \begin{pmatrix}
  x & y & 0 &  0 \\
  \varepsilon y & x & 0 & 0 \\
  0 & 0 & u & v \\
  0 & 0 & \varepsilon v & u
\end{pmatrix} \Big\vert
\begin{array}{c}
xu - \varepsilon yv = 1\\
xv - yu = 0
\end{array}
\right\}.
\]
Frobenius-stable depth-zero cuspidal character sheaves for $\Ksch{T}^6$ take the form
\[
C_6(\theta) \ceq \cind_{\quo{T}^6_{(1)}}^{\Ksch{T}^6}
\fais{L}_\theta[2],
\]
where $\theta$ is a character of $\quo{T}^6_{(0)}(\kq)$. Define
\[
A_6(\theta) = \ind^{\Ksch{G}}_{\Ksch{T}^6} C_6(\theta).
\]
If $\theta$ is in general position then $A_6(\theta)$ is an admissible coefficient system.

\item
Next consider the unramified elliptic torus
\[
\Ksch{T}^7_\Kq = \left\{
\begin{pmatrix}
  x_1 & x_2 & x_3 & x_4 \\
  \varepsilon x_4 & x_1 & x_2 & x_3 \\
  \varepsilon x_3 & \varepsilon x_4 & x_1 & x_2 \\
  \varepsilon x_2& \varepsilon x_3 & \varepsilon x_4 & x_1
\end{pmatrix} \Big\vert
\begin{array}{c}
x_1^2 - 2\varepsilon x_2 x_4 + \varepsilon x_3^2 = 1\\
x_2^2 - 2\varepsilon x_1 x_3 + \varepsilon x_4^2 = 0
\end{array} \right\}.
\]
The image of $I(\Ksch{T}^7,\Knr) \hookrightarrow I(\Ksch{G},\Knr)$ is $\{ (0)\}$. Frobenius-stable admissible coefficient systems for $\Ksch{T}^7_\Kq$ take the form
\[
C_7(\theta) \ceq \cind^{\Ksch{T}^7}_{\quo{T}^7_{(0)}}
\fais{L}_\theta[2],
\]
where $\theta$ is a character of $\quo{T}^7_{(0)}(\kq)$. Define
\[
A_7(\theta) = \ind^{\Ksch{G}}_{\Ksch{T}^7} C_7(\theta).
\]
If $\theta$ is in general position then $A_7(\theta)$ is an admissible coefficient system.

\item
Finally, consider the unramified elliptic torus
\[
\Ksch{T}^8_\Kq = \left\{
\begin{pmatrix}
  x_1 & x_2 & {x_3}{\varpi^{-1}} & {x_4}{\varpi^{-1}} \\
  \varepsilon x_4 & x_1 & {x_2}{\varpi^{-1}} & {x_3}{\varpi^{-1}} \\
  \varpi\varepsilon x_3 & \varpi\varepsilon x_4 & x_1 & x_2 \\
  \varpi\varepsilon x_2& \varpi\varepsilon x_3 & \varepsilon x_4 & x_1
\end{pmatrix} \Big\vert
\begin{array}{c}
x_1^2 - 2\varepsilon x_2 x_4 + \varepsilon x_3^2 = 1\\
x_2^2 - 2\varepsilon x_1 x_3 + \varepsilon x_4^2 = 0
\end{array} \right\}.
\]
The image of $I(\Ksch{T}^8,\Knr) \hookrightarrow I(\Ksch{G},\Knr)$ is $\{ (1)\}$. Frobenius-stable admissible coefficient systems for $\Ksch{T}^8_\Kq$ take the form
\[
C_8(\theta) \ceq \cind^{\Ksch{T}^8}_{\quo{T}^8_{(1)}}
\fais{L}_\theta[2],
\]
where $\theta$ is a character of $\quo{T}^8_{(1)}(\kq)$. Define
\[
A_8(\theta) = \ind^{\Ksch{G}}_{\Ksch{T}^8} C_8(\theta).
\]
If $\theta$ is in general position then $A_8(\theta)$ is an admissible coefficient system.

\item
Let 
\[
\Ksch{L}_\Kq = \left\{ \begin{pmatrix}
  t_1 & 0 & 0 & 0 \\
  0 & a & b & 0 \\
  0 & c & d & 0 \\
  0 & 0 & 0 & t_4
\end{pmatrix} \Big\vert
\begin{array}{c}
t_1t_4 = 1\\
ad-bc =1
\end{array} \right\}.
\]
Up to conjugation, the building for $\Ksch{L}(\Knr)$ has two polyvertices; these may be identified with $(0)$ and $(1)$ in $I(\Ksch{G},\Kq)$. Now $\quo{L}_{(0)}$ and $\quo{L}_{(1)}$ are each isomorphic to $\GL(1)\times \SL(2)$, and the Frobenius-stable cuspidal character sheaves on $\GL(1)\times \SL(2)$ take the form $\fais{L}_\theta\boxtimes K_\pm[1]$, where $\theta$ is a character of $\GL(1,\kq)$ and $K_\pm$ is a cuspidal character sheaf on $\SL(2)_\knr$. (See Example~\ref{example: SL(2)}.) Frobenius-stable cuspidal character sheaves for $\Ksch{L}_\Kq$ take the form
\begin{eqnarray*}
C_0^\pm(\theta)
    &\ceq& \cind^{\Ksch{L}}_{\quo{L}_{(0)}}
    \fais{L}_\theta\boxtimes K_\pm[1]\\
C_0^\pm(\theta)'
    &\ceq& \cind^{\Ksch{L}}_{\quo{L}_{(1)}}
    \fais{L}_\theta\boxtimes K_\pm[1],
\end{eqnarray*}
where $\theta$ is a character of $\GL(1,\kq)$. Define 
\begin{eqnarray*}
B_0^\pm(\theta) &\ceq& \ind^{\Ksch{G}}_{\Ksch{L}} C_0^\pm(\theta)\\
B_0^\pm(\theta)' &\ceq& \ind^{\Ksch{G}}_{\Ksch{L}} C_0^\pm(\theta)'.
\end{eqnarray*}
If $\theta$ is in general position, then $B_0^\pm(\theta)$ and $B_0^\pm(\theta)'$ are admissible coefficient systems.

\item
Next, consider the inner form
\[
\Ksch{L}^1_\Kq = \left\{ \begin{pmatrix}
  a & 0 & 0 & b \\
  0 & x & y & 0 \\
  0 & \varepsilon y & x & 0 \\
  c & 0 & 0 & d
\end{pmatrix} \Big\vert
\begin{array}{c}
x^2-\varepsilon y^2 = 1\\
ad-bc =1
\end{array} \right\}.
\]
Observe that $\Ksch{L}^1(\Kq) = \SU(1,\Kq) \times \SL(2,\Kq)$. Up to $\Ksch{L}^1(\Kq)$ conjugation, the building for
$\Ksch{L}^1(\Kq)$ contains two polyvertices, which may be identified with $(0)$ and $(1)$ in $I(\Ksch{G},\Kq)$. Now
$\quo{L}^1_{(0)/\kq}$ and $\quo{L}^1_{(1)/\kq}$ are each isomorphic to $\SU(1)_\kq\times \SL(2)_\kq$, so we define
\begin{eqnarray*}
C_1^\pm(\theta)
    &\ceq&
    \cind^{\Ksch{L}^1}_{\quo{L}^1_{(0)}} \fais{L}_{\theta}\boxtimes F_\pm[1]\\
C_1^\pm(\theta)'
    &\ceq&
    \cind^{\Ksch{L}^1}_{\quo{L}^1_{(1)}}\fais{L}_{\alpha}\boxtimes F_\pm[1],
\end{eqnarray*}
where $\alpha$ is a character of $\SU(1,\kq)\times \SL(2,\kq)$. Define 
\begin{eqnarray*}
B_1^\pm(\theta) &\ceq& \ind^{\Ksch{G}}_{\Ksch{L}^1} C_1^\pm(\theta)\\
B_1^\pm(\theta)' &\ceq&= \ind^{\Ksch{G}}_{\Ksch{L}^1} C_1^\pm(\theta)'.
\end{eqnarray*}
If $\theta$ is in general position, then $B_1^\pm(\theta)$ and $B_1^\pm(\theta)'$ are admissible coefficient systems.

\item
Likewise,
consider the elliptic unramified-Levi subgroup
\[
\Ksch{L}^2_\Kq = \left\{ \begin{pmatrix}
  x & 0 & 0 & y \varpi^{-1} \\
  0 & a & b & 0 \\
  0 & c & d & 0 \\
  \varepsilon y \varpi & 0 & 0 & x
\end{pmatrix} \Big\vert
\begin{array}{c}
x^2-\varepsilon y^2 = 1\\
ad-bc =1
\end{array} \right\}.
\]
Up to $\Ksch{L}^2(\Kq)$ conjugation, the building for $\Ksch{L}^2(\Kq)$ contains two polyvertices, which may be
identified with $(1)$ and $(2)$ in $I(\Ksch{G},\Kq)$. Define
\begin{eqnarray*}
C_2^\pm(\theta)
    &\ceq&
    \cind^{\Ksch{L}^2}_{\quo{L}^2_{(1)}} \fais{L}_{\theta}\boxtimes F_\pm[1]\\
C_2^\pm(\theta)'
    &\ceq&
    \cind^{\Ksch{L}^2}_{\quo{L}^2_{(2)}} \fais{L}_{\theta}\boxtimes \fais{F}_\pm[1],
\end{eqnarray*}
where $\theta$ is a character of $\quo{L}^2_{(2)}(\kq)$. Define 
\begin{eqnarray*}
B_2^\pm(\theta) &\ceq& \ind^{\Ksch{G}}_{\Ksch{L}^2} C_2^\pm(\theta)\\
B_2^\pm(\theta)' &\ceq&= \ind^{\Ksch{G}}_{\Ksch{L}^2} C_2^\pm(\theta)'.
\end{eqnarray*}
If $\theta$ is in general position, then $B_2^\pm(\theta)$ and $B_2^\pm(\theta)'$ are admissible coefficient systems.

\item
Finally, we turn to the most interesting cuspidal Levi subgroup of $\Ksch{G}$, which is $\Ksch{G}$ itself. A fundamental $\Ksch{G}(\Knr)$-domain for $I(\Ksch{G},\Knr)$ has three polyvertices, which we denote $(0)$, $(1)$ and $(2)$, with $\quo{G}_{(0)/\kq} = \Sp(4)_\kq$, $\quo{G}_{(1)/\kq} = \SL(2)_\kq \times \SL(2)_\kq$ and $\quo{G}_{(2)/\kq} = \Sp(4)_\kq$. There is exactly one cuspidal character sheaf $F_0$ on $\Sp(4)$ while $\SL(2)_\knr\times\SL(2)_\knr$ admits four cuspidal character sheaves, being $F_\pm \boxtimes F_\pm$ in the notation from Example~\ref{example: SL(2)}; of these, only $F_+ \boxtimes F_+$ is a cuspidal unipotent character sheaf. Thus, by Corollary~\ref{corollary: cuspidal}, the elements of $\catqA^{(0)}\Ksch{G}_\Kq$ are
\begin{eqnarray*}
C^0 &\ceq& \cind^{\Ksch{G}}_{\quo{G}_{(0)}} F_0\\
D^{+\, +} &\ceq& \cind^{\Ksch{G}}_{\quo{G}_{(1)}} F_+\boxtimes F_+\\
D^{+\, -} &\ceq& \cind^{\Ksch{G}}_{\quo{G}_{(1)}} F_+\boxtimes F_-\\
D^{-\, +} &\ceq& \cind^{\Ksch{G}}_{\quo{G}_{(1)}} F_-\boxtimes F_+\\
D^{-\, -}&\ceq&\cind^{\Ksch{G}}_{\quo{G}_{(1)}} F_-\boxtimes F_-\\
C^2 &\ceq& \cind^{\Ksch{G}}_{\quo{G}_{(2)}} F_0.
\end{eqnarray*}
\end{itemize}

This completes the list of all frobenius-stable cuspidal
coefficient systems for each cuspidal unramified
twisted-Levi subgroup of $\Ksch{G}_\Kq$. Every depth-zero
character sheaf for $\Ksch{G}_\Kq$ is a simple summand of an
object of $\catqC\Ksch{G}_\Kq$ produced by parabolic induction
from of the cuspidal coefficient systems listed above.

Each irreducible depth-zero supercuspidal representation of $\Ksch{G}(\Kq)$ is equivalent to a representation obtained by compact induction from an irreducible cuspidal representation of $\quo{G}_{(0)}(\kq)$ or $\quo{G}_{(1)}(\kq)$ or $\quo{G}_{(2)}(\kq)$, so we begin by listing all irreducible cuspidal representations of these finite groups. Consider the finite group $\Sp(4,\kq)$. Using notation from \cite{S1} and \cite{S2}, each irreducible cuspidal representation of $\Sp(4,\kq)$ the appears in one of the following families: 
\begin{itemize}
\item
the Deligne-Lusztig representation $\chi_1 = R^{\Sp(4)}_{\ksch{T}_1}(\theta)$ where $\theta$ is a character of
$\ksch{T}_1(\kq)$ in general position; 
\item
the Deligne-Lusztig representation $\chi_4 = R^{\Sp(4)}_{\ksch{T}_4}(\theta)$ where
$\theta$ is a character of $\ksch{T}_4(\kq)$ in general position; 
\item
the irreducible constituent $\xi'_{21}$ of the Deligne-Lusztig
virtual representation $R^{\Sp(4)}_{\ksch{T}_4}(\theta\times\sgn)$ where $\theta$ is a character of $\SU(1,\kq)$ in general position and $\sgn$ is the sign character of $\SU(1,\kq)$; 
\item
the other irreducible constituent $\xi'_{22}$ of the Deligne-Lusztig virtual representation $R^{\Sp(4)}_{\ksch{T}_4}(\theta\times\sgn)$;  
\item
the cuspidal unipotent representation $\theta_{10}$. 
\end{itemize}
With these conventions, and notation from Example~\ref{example: SL(2)}, every irreducible depth-zero supercuspidal representation of $\Ksch{G}(\Kq)$ is equivalent to one of the following:
\begin{eqnarray*}
\pi_{0} &=& \cInd^{\Ksch{G}(\Kq)}_{\Rsch{G}_{(1)}(\Rq)}(\chi_0\times \chi_0)\\
\pi_{1} &=& \cInd^{\Ksch{G}(\Kq)}_{\Rsch{G}_{(1)}(\Rq)}(\chi_0\times \chi_0^+)\\
\pi_{2} &=& \cInd^{\Ksch{G}(\Kq)}_{\Rsch{G}_{(1)}(\Rq)}(\chi_0\times \chi_0^-)\\
\pi_{3} &=& \cInd^{\Ksch{G}(\Kq)}_{\Rsch{G}_{(1)}(\Rq)}(\chi_0^+\times \chi_0^+)\\
\pi_{4} &=& \cInd^{\Ksch{G}(\Kq)}_{\Rsch{G}_{(1)}(\Rq)}(\chi_0^+\times \chi_0^-)\\
\pi_{5}&=&\cInd^{\Ksch{G}(\Kq)}_{\Rsch{G}_{(1)}(\Rq)}(\chi_0^-\times\chi_0^-).
\end{eqnarray*}
and
\begin{equation}
\begin{array}{ccc}
\pi_{6} = \cInd^{\Ksch{G}(\Kq)}_{\Rsch{G}_{(0)}(\Rq)}(\chi_1)
    && \pi'_{6} = \cInd^{\Ksch{G}(\Kq)}_{\Rsch{G}_{(2)}(\Rq)}(\chi_1)\\
\pi_{7} = \cInd^{\Ksch{G}(\Kq)}_{\Rsch{G}_{(0)}(\Rq)}(\chi_4)
    && \pi'_{7} = \cInd^{\Ksch{G}(\Kq)}_{\Rsch{G}_{(2)}(\Rq)}(\chi_4)\\
\pi_{8} = \cInd^{\Ksch{G}(\Kq)}_{\Rsch{G}_{(0)}(\Rq)}(\xi'_{22})
    && \pi'_{8} = \cInd^{\Ksch{G}(\Kq)}_{\Rsch{G}_{(2)}(\Rq)}(\xi'_{22})\\
\pi_{9} = \cInd^{\Ksch{G}(\Kq)}_{\Rsch{G}_{(0)}(\Rq)}(\xi'_{21})
    && \pi'_{9}=\cInd^{\Ksch{G}(\Kq)}_{\Rsch{G}_{(2)}(\Rq)}(\xi'_{21})\\
\pi_{10}=\cInd^{\Ksch{G}(\Kq)}_{\Rsch{G}_{(0)}(\Rq)}(\theta_{10})
    && \pi'_{10} = \cInd^{\Ksch{G}(\Kq)}_{\Rsch{G}_{(2)}(\Rq)}(\theta_{10})\\
\end{array}
\end{equation}


\subsection{GL(n)}\label{example: GL(n)}

Let $\Ksch{G} = \GL(n)$. By Definition~\ref{definition: admissible}, every irreducible admissible coefficient system is a summand of $\ind^{\Ksch{G}}_{\Ksch{L}} C$, where $\Ksch{L}$ is a Levi subgroup of $\Ksch{G}$ and $C$ is a cuspidal coefficient system (\cf Definition~\ref{definition: cuspidal}). Therefore, the first step in describing irreducible admissible coefficient systems for $\Ksch{G}$ is to enumerate all Levi subgroups and all cuspidal coefficient systems on those Levi subgroups. (\cf Definition~\ref{definition: cuspidal Levi}.) Every Levi subgroup of $\Ksch{G}$ is $\Ksch{G}(\Knr)$-conjugate to $\Ksch{L}^{\underline{m}}$, where ${\underline{m}} = [m_1, m_2, \ldots , m_t]$ is a partition of $n$ and $\Ksch{L}^{{\underline{m}}} = \prod_{k=1}^{t} \GL(m_k)_\Knr$. By Corollary~\ref{corollary: cuspidal}, every cuspidal coefficient system on $\Ksch{L}^{\underline{m}}$ is isomorphic to $\cind^{\Ksch{L}^{\underline{m}}}_{\Ksch{L}^{\underline{m}}_i} F$, where $i$ is a vertex of the building
$I(\Ksch{L}^{\underline{m}},\Knr)$ and $F$ is a cuspidal character sheaf for $\quo{L}^{\underline{m}}_i$. Since the building for $\Ksch{L}^{\underline{m}}$ is regular and all vertices are $\Ksch{L}^{\underline{m}}(\Knr)$-conjugate, we have $\quo{L}^{\underline{m}}_i = \prod_{k=1}^{t} \GL(m_k)_\knr$. Thus, cuspidal character sheaves on $\quo{L}^{\underline{m}}_i$ are all of the form $\mathop{\boxtimes}_{k=1}^{t} F_k$, where $F_k$ is a cuspidal character sheaf on $\GL(m_k)_\knr$ (\cf \cite[Lemma~5.4.1]{MS} for example). Since $\GL(m_k)_\knr$ admits cuspidal character sheaves if and only if $m_k =1$, it follows that $\Ksch{L}^{\underline{m}}$ admits cuspidal coefficient systems if and only if ${\underline{m}} = [1^m]$, whence $\Ksch{L}^{\underline{m}}$ is a $\Knr$-split torus. Let $\Ksch{T}$ be a $\Knr$-split torus (there is exactly one in $\Ksch{G}$, up to $\Ksch{G}(\Knr)$-conjugacy). From the discussion above we see that every cuspidal coefficient system for $\Ksch{T}$ takes the form
\begin{equation}\label{equation: GLn 1}
\cind^{\Ksch{T}}_{\Ksch{T}_{(0)}} \fais{L}_1 \boxtimes \fais{L}_2
\boxtimes \cdots \boxtimes \fais{L}_n [n],
\end{equation}
where $(0)$ is \emph{the} vertex of $I(\Ksch{T},\Knr)$ and $\fais{L}_k$ is a Kummer local system on $\GL(1)_\knr$ for each $k=1, \ldots ,n$. Thus, 
\begin{equation}
\ind^{\Ksch{G}}_{\Ksch{T}} C = \ind^{\Ksch{G}}_{\Ksch{T}}
\cind^{\Ksch{T}}_{\Ksch{T}_{(0)}}\fais{L}_1 \boxtimes \fais{L}_2
\boxtimes \cdots \boxtimes \fais{L}_n [n].
\end{equation}
is an admissible coefficient system for $\Ksch{G}$.

We now suppose $\Ksch{G} = \GL(n)_\Kq \times_\Spec{\Kq} \Spec{\Knr}$, so $\Ksch{G}_\Kq = \GL(n)_\Kq$. Let $\Ksch{S}_\Kq$ be the $\Kq$-split torus in $\Ksch{G}_\Kq$ and let $\Ksch{T}_\Kq$ be an inner form of $\Ksch{S}_\Kq$. (Following our conventions, we have $\Ksch{T} = \Ksch{T}_\Kq \times_\Spec{\Kq} \Spec{\Knr}$.)
Then
\begin{equation}
\Ksch{T}^{\underline{n}}_\Kq = \prod_{k=1}^t
\Res_{\mathbb{K}_{n_k}/\Kq}GL(1)_{\mathbb{K}_{n_k}},
\end{equation}
where ${\underline{n}} = [n_1,n_2, \ldots, n_t]$ is a partition of $n$ and $\mathbb{K}_{n_k}$ is an unramified extension of $\Kq$ is degree $n_k$. For each $k=1, \ldots ,t$, let \[\Ksch{T}^{n_k}_\Kq =
\Res_{\mathbb{K}_{n_k}/\Kq} GL(1)_{\mathbb{K}_{n_k}}.\] Then
\[\quo{T}^{n_k}_{(0)/\kq} = \Res_{\Bbbk_{n_k}/\kq} GL(1)_{\Bbbk_{n_k}}.\] Let $\theta_k$ be a character of $\quo{T}^{n_k}_{(0)}(\kq)$ and let $\fais{L}_{\theta_k}$ be the Kummer local system on
$\quo{T}^{n_k}_{(0)}$ equipped with an isomorphism $\Frob^* \fais{L}_{\theta_k} \to \fais{L}_{\theta_k}$ in
$D^b_c(\quo{T}^{n_k}_{(0)};\EE)$ such that $\chf{\fais{L}_{\theta_k}} = \theta_k$. Every frobenius-stable cuspidal coefficient system for $\Ksch{T}^{{\underline{n}}}$ takes the form
\begin{equation}
C^{\underline{n}}_\theta \ceq \cind^{\Ksch{T}^{\underline{n}}}_{\Ksch{T}^{\underline{n}}_{(0)}}
\mathop{\boxtimes}_{k=1}^t \fais{L}_{\theta_k} [n],
\end{equation}
with $\theta = \mathop{\otimes}_{k=1}^t \theta_k$ a character of $\quo{T}^{\underline{n}}_{(0)}(\kq)$, and every frobenius-stable irreducible admissible coefficient system for $\Ksch{G}$ is a summand of
\begin{equation}
A^{\underline{n}}_\theta \ceq \ind^{\Ksch{G}}_{\Ksch{T}^{\underline{n}}} C^{\underline{n}}_\theta.
\end{equation}
If each character $\theta_k$ appearing in $\theta$ is in general position, in the sense of \cite{DL}, then $A^{\underline{n}}_\theta$ is irreducible, and therefore an irreducible admissible coefficient system itself.

One case is particularly important to us. If ${\underline{n}} = [n]$ we denote $\Ksch{T}^{\underline{n}}$ (resp. $C_\theta^{\underline{n}}$, $A_\theta^{\underline{n}}$) by $\Ksch{T}^n$ (resp. $C_\theta^n$, $A_\theta^n$). In this case
\begin{equation}
\Ksch{T}^{n}_\Kq = \Res_{\mathbb{K}_{n}/\Kq} GL(1)_{\mathbb{K}_{n}}
\end{equation}
is elliptic and $\theta$ is a character of $GL(1,\mathbb{K}_{n})$. When $\theta$ is in general position, $A_\theta$ is an irreducible admissible coefficient system. This case is sufficient for a description of the models of all supercuspidal depth-zero representations, as one sees from the proof of Theorem~\ref{theorem: supercuspidal models}.

Let $\pi$ be an irreducible supercuspidal depth-zero representation of $\Ksch{G}(\Kq)$. Then there is a non-trivial
character $\chi$ of $\Ksch{ZG}(\Kq)$ with conductor $\Ksch{ZG}(\Kq)_0$ such that $\pi$ is equivalent to
\begin{equation}
\pi_{\chi,\theta} \ceq \cInd^{\GL(n,\Kq)}_{\Ksch{ZG}(\Kq) \,
\GL(n,\Rq)}(\chi \otimes \sigma),
\end{equation}
where $\sigma$ is a representation of $\GL(n,\Rq) = \Rsch{G}_{(0)}(\Rq)$ produced by inflation from a cuspidal
irreducible representation $\quo{\sigma}$ of $\GL(n,\kq) = \quo{G}_{(0)}(\kq)$. Thus, there is a character $\theta$ of
$\quo{T}^n_{(0)}(\kq)$ in general position such that $\quo{\sigma}$ is equivalent to $(-1)^{n-1} R^{\quo{G}_{(0)}}_{\quo{T}^n_{(0)}}\theta$.

Using this notation, the models for depth-zero supercuspidal representations of $\GL(n,\Kq)$ are presented in Table~\ref{table: models for GLn}.

\begin{table}[htdp]
\caption{Models for depth-zero supercuspidal representations of $\GL(n,\Kq)$}
\begin{center}
\begin{tabular}{|c|c|}
\hline 
{Representation} & {Model}\\
\hline
$\pi_{\chi,\theta}$ & $(-1)^{n-1}
\ind^{\Ksch{G}}_{\Ksch{T}^n}
\cind^{\Ksch{T}^{n}}_{\Ksch{T}^{n}_{(0)}} \fais{L}_{\theta} [n]$ \\
\hline
\end{tabular}
\end{center}
\label{table: models for GLn}
\end{table}


\bibliographystyle{amsplain}

\end{document}